\documentclass[12pt,psamsfonts,leqno,oneside,letterpaper]{amsart}
\usepackage[dvips,text={6.5truein,9truein},left=1truein,top=1truein]{geometry}
\usepackage{amssymb,amsmath,amscd,enumerate,
}
\usepackage[pdftex]{graphicx}
\usepackage{url}

\usepackage[colorlinks,linkcolor=blue,citecolor=blue,pdfstartview=FitH]{hyperref}
\input xy
\xyoption{all}
\SelectTips{cm}{12}

\usepackage{color}
\usepackage{mathrsfs}
\newcommand\redden[1]{{\color{red}#1}}

\newcommand\missingref[1]{\empty}
\newcommand\tinymissingref[1]{\empty}
\newcommand\abstractcomment[1]{\empty}

\theoremstyle{definition}
\newtheorem{para}{}[section]

\newtheorem{remark}[para]{Remark}
\newtheorem{reformulation}[para]{Reformulation}
\newtheorem{remarks}[para]{Remarks}
\newtheorem{notation}[para]{Notation}

\newtheorem{convention}[para]{Convention}
\newtheorem{definition}[para]{Definition}
\newtheorem{definitions}[para]{Definitions}
\newtheorem{definitionnotation}[para]{Definition and Notation}
\newtheorem{remarksnotation}[para]{Remarks and Notation}

\newtheorem{remarknotation}[para]{Remark and Notation}
\newtheorem{notationremark}[para]{Notation and Remark}
\newtheorem{notationreviewremarks}[para]{Notation, Review and Remarks}
\newtheorem{definitionremark}[para]{Definition and Remark}
\newtheorem{definitionsremarks}[para]{Definitions and Remarks}
\newtheorem{notationremarks}[para]{Notation and Remarks}
\newtheorem{definitionsnotation}[para]{Definitions and Notation}

\newtheorem{reviewdefinition}[para]{Review and Definition}

\newtheorem{definitionnotationremarks}[para]{Definition, Notation and Remarks}
\newcommand\Alternatives{\begin{enumerate}[(i)]}
\newcommand\EndAlternatives{\end{enumerate}}
\newcommand\Conditions{\begin{enumerate}[(1)]}
\newcommand\EndConditions{\end{enumerate}}

\theoremstyle{plain}
\newtheorem{theorem}[para]{Theorem}

\newtheorem{lemma}[para]{Lemma}
\newtheorem{proposition}[para]{Proposition}

\newtheorem{corollary}[para]{Corollary}
\newtheorem{conjecture}[para]{Conjecture}
\newtheorem{specialremark}{Remark}
\newtheorem*{propgoodhomology}{Proposition \ref{good homology}}
\newtheorem*{propthreesevensevenhomology}{Proposition \ref{3.77 homology}}
\newtheorem*{fourfreeprop}{Proposition \ref{four-free case}}
\newtheorem*{fivefreeprop}{Proposition \ref{five-free case}}
\newtheorem{claim}[equation]{}
\numberwithin{equation}{para}
\numberwithin{figure}{section}
\numberwithin{specialremark}{para}

\newcommand\Number{\begin{para}}
\newcommand\EndNumber{\end{para}}
\newcommand\Definition{\begin{definition}}
\newcommand\EndDefinition{\end{definition}}
\newcommand\Definitions{\begin{definitions}}
\newcommand\DefinitionsNotation{\begin{definitionsnotation}}
\newcommand\DefinitionNotation{\begin{definitionnotation}}
\newcommand\RemarksNotation{\begin{remarksnotation}}
\newcommand\Reformulation{\begin{reformulation}}
\newcommand\EndRemarksNotation{\end{remarksnotation}}
\newcommand\EndReformulation{\end{reformulation}}
\newcommand\RemarkNotation{\begin{remarknotation}}
\newcommand\EndRemarkNotation{\end{remarknotation}}
\newcommand\NotationRemark{\begin{notationremark}}
\newcommand\EndDefinitionNotationRemarks{\end{definitionnotationremarks}}
\newcommand\NotationReviewRemarks{\begin{notationreviewremarks}}
\newcommand\DefinitionRemark{\begin{definitionremark}}
\newcommand\DefinitionsRemarks{\begin{definitionsremarks}}
\newcommand\DefinitionNotationRemarks{\begin{definitionnotationremarks}}
\newcommand\NotationRemarks{\begin{notationremarks}}
\newcommand\EndNotationRemark{\end{notationremark}}
\newcommand\EndNotationReviewRemarks{\end{notationreviewremarks}}
\newcommand\EndDefinitionRemark{\end{definitionremark}}
\newcommand\EndDefinitionsRemarks{\end{definitionsremarks}}
\newcommand\EndNotationRemarks{\end{notationremarks}}
\newcommand\EndDefinitionsNotation{\end{definitionsnotation}}
\newcommand\EndDefinitionNotation{\end{definitionnotation}}
\newcommand\ReviewDefinition{\begin{reviewdefinition}}
\newcommand\EndReviewDefinition{\end{reviewdefinition}}
\newcommand\EndDefinitions{\end{definitions}}
\newcommand\Theorem{\begin{theorem}}
\newcommand\EndTheorem{\end{theorem}}
\newcommand\Conjecture{\begin{conjecture}}
\newcommand\EndConjecture{\end{conjecture}}
\newcommand\Remark{\begin{remark}}
\newcommand\EndRemark{\end{remark}}
\newcommand\Remarks{\begin{remarks}}
\newcommand\EndRemarks{\end{remarks}}
\newcommand\Convention{\begin{convention}}
\newcommand\EndConvention{\end{convention}}
\newcommand\Notation{\begin{notation}}
\newcommand\EndNotation{\end{notation}}
\newcommand\Lemma{\begin{lemma}}
\newcommand\EndLemma{\end{lemma}}
\newcommand\Proposition{\begin{proposition}}
\newcommand\EndProposition{\end{proposition}}
\newcommand\Corollary{\begin{corollary}}
\newcommand\EndCorollary{\end{corollary}}
\newcommand\Claim{\begin{claim}}
\newcommand\EndClaim{\end{claim}}
\newcommand\Proof{\begin{proof}}
\newcommand\EndProof{\end{proof}}
\newcommand\Equation{\begin{equation}}
\newcommand\EndEquation{\end{equation}}

\newcommand\Bullets{\begin{itemize}}
\newcommand\EndBullets{\end{itemize}}

\newcommand\admissible{admissible}
\newcommand\admissibility{admissibility}

\newcommand\frakN{\mathfrak N}

\newcommand\frakW{{{\mathfrak W}}}

\renewcommand\epsilon{\varepsilon}

\newcommand\kish{\mathop{\rm kish}}

\newcommand\RR{{\bf R}}
\newcommand\fraks{{\mathfrak s}}

\newcommand\waspsi{\scrB}
\newcommand\wasdelta{\epsilon}
\newcommand\waschi{\scrA}

\newcommand\chibar{\overline\chi}

\newcommand\cosech{\mathop{\rm cosech}\nolimits}

\newcommand\Mthick{M_{\rm thick}}
\newcommand\Mthin{M_{\rm thin}}

\newcommand\barOmega{\overline{\Omega}}

\newcommand\tP{\widetilde P}

\newcommand\waslambdaminus{\calL}
\newcommand\maybecalLminus{\calL^-}
\newcommand\maybecalLplus{\calL^+}
\newcommand\wasell{l}

\newcommand\calc{{\mathcal C}}

\newcommand\frakX{{\mathfrak X}}

\newcommand\calf{{\mathcal F}}

\newcommand\calL{{\mathcal L}}

\newcommand\frakG{{\mathfrak G}}

\newcommand\calq{{\mathcal Q}}

\newcommand\calg{{\mathcal G}}

\newcommand\calu{{\mathcal U}}
\newcommand\calU{{\mathcal U}}

\newcommand\calw{{\mathcal W}}

\newcommand\Isom{{\rm Isom}}
\newcommand\arccosh{{\rm arccosh}}

\newcommand\arcsech{{\rm arcsech}}

\newcommand\sech{{\rm sech}}

\newcommand\arcsec{{\rm arcsec}}
\newcommand\ZZ{{\bf Z}}

\newcommand\HH{{\bf H}}

\newcommand\cals{{\mathscr S}}

\newcommand\scrA{{\mathscr A}}
\newcommand\scrB{{\mathscr B}}

\newcommand\scrD{{\mathscr D}}
\newcommand\scrH{{\mathscr H}}
\newcommand\scrM{{\mathscr M}}

\newcommand\iccg{ICC-group}
\newcommand\link{\mathop{\rm link}}

\newcommand\scrv{{\mathscr V}}

\newcommand\dist{\mathop{\rm dist}}

\newcommand\nbhd{\mathop{\rm nbhd}}

\newcommand\density{{\rm density}}

\newcommand\vol{\mathop{\rm vol}}

\newcommand\voct{V_{\rm oct}}

\newcommand\length{\mathop{{\rm length}}}

\newcommand\rank{\mathop{{\rm rank}}}
\newcommand\localrank{\mathop{{\rm local rank}}}

\newcommand\isomplus{\mathop{{\rm Isom}_+}}

\newcommand\restraining{restraining}

\newcommand\clo{{\bold c}}
\newcommand\FF{{\bf F}}
\newcommand\ltr{\langle t\rangle}

\newcommand\iof{{\rm iof}}

\newcommand\shortone{{\mathfrak s}_1}
\newcommand\nextone{{\mathfrak s}_2}

\newcommand\f{f}
\newcommand\fone{f_3}
\newcommand\fthree{f_1}
\newcommand\Wone{V^{\rm near}_{\rm ST}}
\newcommand\Woneminus{V^{-}_{\rm ST}}

\newcommand\g{g}

\newcommand\Vbor{V_{\text{\rm B\"or}}}
\newcommand\tVbor{{\widetilde V}_{\text{\rm B\"or}}}

\DeclareFontFamily{U}{rcjhbltx}{}
\DeclareFontShape{U}{rcjhbltx}{m}{n}{<->rcjhbltx}{}
\DeclareSymbolFont{hebrewletters}{U}{rcjhbltx}{m}{n}

\let\aleph\relax\let\beth\relax
\let\gimel\relax\let\daleth\relax

\DeclareMathSymbol{\aleph}{\mathord}{hebrewletters}{39}
\DeclareMathSymbol{\beth}{\mathord}{hebrewletters}{98}
\DeclareMathSymbol{\gimel}{\mathord}{hebrewletters}{103}
\DeclareMathSymbol{\daleth}{\mathord}{hebrewletters}{100}

\DeclareMathSymbol{\lamed}{\mathord}{hebrewletters}{108}
\DeclareMathSymbol{\mem}{\mathord}{hebrewletters}{109}
\DeclareMathSymbol{\ayin}{\mathord}{hebrewletters}{96}
\DeclareMathSymbol{\tsadi}{\mathord}{hebrewletters}{118}
\DeclareMathSymbol{\qof}{\mathord}{hebrewletters}{114}
\DeclareMathSymbol{\shin}{\mathord}{hebrewletters}{152}

\newcommand\Wtwo{V^{\rm near}_{\rm SG}}
\newcommand\W{W}
\newcommand\Wminus{W^-}
\newcommand\Wthree{V^{\rm far}}
\newcommand\Wfour{V_{\rm VSG}}
\newcommand\Wfive{V_{{\rm D}\text{-}{\rm R}}}
\newcommand\VAD{V_{\rm AD}}

\usepackage[utf8]{inputenc}

\title{Hyperbolic volume, mod $2$ homology, and {$\MakeLowercase{k}$}-freeness}
\author{Rosemary K. Guzman}
\address{Department of Mathematics
\\
University of Illinois\\
1409 W. Green St.\\
Urbana, IL 61801}
\email{rguzma1@illinois.edu}

\author{Peter B. Shalen}
\address{Department of Mathematics, Statistics, and Computer Science
(M/C 249)\\
University of Illinois at Chicago\\
851 S. Morgan St.\\
Chicago, IL 60607-7045}
\email{petershalen@gmail.com}

\begin{document}

\begin{abstract}
We show that if $M$ is any closed, orientable hyperbolic $3$-manifold with $\vol M\le 3.69$, we have $\dim H_1(M;\FF_2)\le7$. This may be regarded as a qualitative improvement of 
one of the main results 
of \cite{fourfree}, because the constant $3.69$ is greater than the ordinal 
corresponding 
to $\omega^2$ in the well-ordered set of finite volumes of hyperbolic $3$-manifolds. We also show that if $\vol M\le 3.77$, we have $\dim H_1(M;\FF_2)\le10$.

These results are applications of a new method for obtaining lower bounds for the volume of  a closed, orientable hyperbolic $3$-manifold such that $\pi_1(M)$ is $k$-free for a given $k\ge4$. Among other  applications we show that if $\pi_1(M)$ is $4$-free we have $\vol M>3.57$ (improving the lower bound of $3.44$ given in \cite{fourfree}), 
and that if $\pi_1(M)$ is $5$-free we have $\vol M>3.77$.
\end{abstract}

\maketitle
    
\section{Introduction}

It is a consequence of the Mostow Rigidity Theorem that a finite-volume hyperbolic $3$-manifold is determined up to isometry by its topological type. Hence any geometrically defined invariant of such a manifold, such as its volume, may be regarded as a topological invariant. The theme of such papers as \cite{paradoxical}, \cite{accs}, \cite{AST}, \cite{cusp}, \cite{acs-surgery}, \cite{CDS}, \cite{ds}, \cite{acs-singular},  \cite{singular-two}, \cite{fourfree}, which will be pursued further in this paper, is to develop explicit quantitative relationships between the volume $\vol M$ of a closed, orientable hyperbolic $3$-manifold $M$ and topologically defined numerical invariants of $M$.

To illustrate the theme of ``quantitative Mostow rigidity,'' we shall first discuss the relationship between the volume and the homology of a closed, orientable hyperbolic $3$-manifold $M$. It has long been known that for any prime $p$, the dimension of $H_1(M;\FF_p)$ (where $\FF_p$ denotes the field of order $p$) is linearly bounded in terms of $\vol M$. According to \cite[Proposition 2.2]{alm}, which builds on results proved in \cite{paradoxical}, \cite{agol} and \cite{cg}, 
we have 
$\dim H_1(M;\FF_p)\le334.08\cdot\vol M$  for every prime $p$. For small values of $\vol M$, 
these results were improved by a couple of orders of magnitude in \cite{acs-surgery}, \cite{singular-two}, and \cite{fourfree}. Theorem 1.1 of \cite{acs-surgery} asserts that if $\vol M\le1.22$ then $\dim H_1(M;\FF_p)\le 2$ for  $p\ne 2,7$, while $\dim H_1(M;\FF_p)\le 3$ if $ p$ is $2$ or $7$; 
this result is sharp, as we have $\dim H_1(M;\FF_5)= 2$ when $M$ is the Weeks manifold.
Theorem 1.2 of \cite{singular-two} asserts that if $\vol M\le3.08$ then $\dim H_1(M;\FF_2)\le 5$. Theorem 1.7 of  \cite{fourfree} asserts that if $\vol M\le3.44$ then $\dim H_1(M;\FF_2)\le 7$.

One of the results of this paper is a strict improvement over Theorem 1.7 of  \cite{fourfree}: 

\begin{propgoodhomology}
If $M$ is any closed, orientable hyperbolic $3$-manifold with $\vol M\le 3.69$, we have $\dim H_1(M;\FF_2)\le7$.
\end{propgoodhomology}

To provide perspective on this improvement, we will now explain some background from \cite{notes} (see especially \cite[Corollaries 6.6.1---6.6.3]{notes}). 
The set $\scrv$ of all finite volumes of orientable hyperbolic $3$-manifolds is a well-ordered set in the ordering inherited from the standard ordering of the real numbers. Each element of $\scrv$ is represented by only finitely many hyperbolic manifolds (up to isometry). The ordinal type of $\scrv$ is $\omega^\omega$, 
so that there is a unique  order-preserving bijection $\alpha\mapsto V_\alpha$ from the set of ordinal numbers less than $\omega^\omega$ to  $\scrv$. For every non-limit ordinal $\alpha$, the number $V_\alpha$ is realized as the volume of a {\it closed} orientable hyperbolic $3$-manifold.

It follows from the main theorem of \cite{two-cusps} that $V_{\omega^2}$  is equal to $\voct=3.66\ldots$, the volume of a regular ideal hyperbolic octahedron; it is realized as the volume of the complement of the Whitehead link, an orientable hyperbolic $3$-manifold. Thus the set $\scrv\cap(0,3.69]$ has ordinal type at least $\omega^2$, whereas the set $\scrv\cap(0,3.44]$ has ordinal type only  $m\omega+n$ for some integers $m,n\ge0$. 
(The ordinal types of the sets of {\it known} elements of $\scrv$ that are at most $3.44$ and $3.69$ respectively are $36\omega+1$ and 
$\omega^2 + \omega + 37$ 
(\cite{nathan}). 
For this reason we regard Proposition \ref{good homology} of the present paper as a {\it qualitative} improvement over Theorem 1.7 of  \cite{fourfree}.

We will also prove:
\begin{propthreesevensevenhomology}
If $M$ is any closed, orientable hyperbolic $3$-manifold with $\vol M\le 3.77$, we have $\dim H_1(M;\FF_2)\le 10$. 
\end{propthreesevensevenhomology}

The ordinal type of the set of  known elements of $\scrv$ that are at most  $3.77$ is $\omega^2+8\omega+24$ 
(\cite{nathan}).

Although they
%, like the  results cited above from \cite{acs-surgery}, \cite{singular-two} and \cite{fourfree}, are vast 
provide qualitative improvements over previously known results, there is no reason to think that Propositions \ref{good homology} and \ref{3.77 homology} are sharp. Indeed, the smallest known closed orientable $3$-manifold $M$ with $\dim H_1(M;\FF_2)\ge 4$ has volume $6.35\ldots$.

In order to explain both the ideas in the proofs of Propositions \ref{good homology} and \ref{3.77 homology}, 
and the statements of the other results in the present paper, we recall some ideas that are used in the papers cited above.
%The main theorem of this paper, Theorem \ref{five-free estimate}, from which Theorem \ref{homology theorem} is derived, is similar in nature to the results  in earlier papers that underlie the results mentioned above about manifolds of volume at most $3.08$ or $3.44$. Recall that
The {\it rank} of a  group $A$ is defined to be the minimal cardinality of a generating set for $A$.
A group $\Gamma$ is said to be {\it $k$-free}, where $k$ is a  given
non-negative 
integer, if every subgroup of $\Gamma$ having rank at most $k$ is a free group. 
A major ingredient in the proofs of the results cited above from \cite{acs-surgery}, \cite{singular-two}, and \cite{fourfree} is that for certain values of $k$, one can establish a stronger lower bound for the volume of a closed, orientable hyperbolic $3$-manifold $M$ under the assumption that $\pi_1(M)$ is $k$-free. For example, according to \cite[Corollary 9.3]{acs-singular}, which is an extension of \cite[Theorem 6.1]{accs}, if  $\pi_1(M)$ is  $3$-free then $\vol M>3.08$.
According to \cite[Theorem 1.6]{fourfree}, if $\pi_1(M)$ is $4$-free then $\vol M>3.44$. (A minor correction to the proof of the latter theorem was given in \cite{correction}.) The transition from these results involving $k$-freeness to the results of \cite[Theorem 1.2]{singular-two} and \cite[Theorem 1.7]{fourfree} about homology involve both the topological results developed in \cite{JS}, \cite{sw}, \cite{bs},
 \cite{acs-singular},  \cite{singular-two}, and the geometric results of \cite{AST}.

One of the main ingredients in the proofs of  Propositions \ref{good homology} and \ref{3.77 homology} 
is a deep and surprising refinement of the method
used in \cite{fourfree} to obtain
lower bounds for $\vol M$ when $\pi_1(M)$ is $4$-free. 
A sketch of this refined method will be given below.

This method is also used to prove another one of the main results of this paper: 
\begin{fourfreeprop}
Every closed, orientable hyperbolic $3$-manifold with $4$-free fundamental group has volume greater than $3.57$.
\end{fourfreeprop}
In addition, the method gives stronger lower bounds for $\vol M$ under the assumption that $\pi_1(M)$ is $k$-free for a given $k>4$; this is illustrated by another result which is proved in the present paper:
\begin{fivefreeprop} 
Every closed, orientable hyperbolic $3$-manifold with $5$-free fundamental group has volume greater than $3.77$.
\end{fivefreeprop}

Preliminary calculations seem to show that by directly applying the methods of this paper under the assumption that a given closed, orientable hyperbolic $3$-manifold has $6$-free fundamental group, one would at best get a lower bound of about $3.80$ for the volume. As this would be a rather modest improvement over Proposition \ref{five-free case}, we have not included it.

Before discussing our refined methods for estimating $\vol M$ when $\pi_1(M)$ is $4$-free, we must of course review the ingredients in the proof of \cite[Theorem 1.6]{fourfree}. 
The key step in the latter proof is \cite[Theorem 1.4]{fourfree}, which involves the notion of a ``semithick point.'' Let $\lambda$ be a positive number. It is standard to define a point $p$ of a hyperbolic $3$-manifold $M$ to be {\it $\lambda$-thick} if every loop based at $p$ having length  less than $\lambda$ is  homotopically trivial. We say that $p$ is {\it $\lambda$-semithick} if there is a cyclic subgroup $C$ of $\pi_1(M,p)$ such that every loop based at $p$ having length  less than $\lambda$ defines an element of $C$. 
(Here and elsewhere in this introduction, we have paraphrased definitions and statements given in the body of the paper so as to minimize the amount of formalism needed to express them.)
Theorem 1.4 of \cite{fourfree} asserts that if $M$ is closed and orientable, and $\pi_1(M)$ is $4$-free, then $M$ contains a $\log7$-semithick point.

The proof of \cite[Theorem 1.4]{fourfree} begins by writing $M$ as $\HH^3/\Gamma$, where $\Gamma$ is a cocompact $4$-free group of orientation-preserving isometries of $\HH^3$. For each maximal cyclic subgroup $C$ of $\Gamma$ and each $\lambda>0$, we denote by $Z_\lambda(C)$ the set of all points $P\in\HH^3$ such that $\dist(P,x\cdot P)<\lambda$ for some $x\in C-\{1\}$. If $M$ contains no $\log7$-semithick point, then in particular it contains no $\log7$-thick point; this is shown to imply that the non-empty sets of the form $Z_{\log7}(C)$ form a covering of $\HH^3$, indexed by certain maximal cyclic subgroups of $\Gamma$. The nerve of this covering (in the sense of algebraic topology) is a simplicial complex $K$ (which we think of as being geometricallhy realized). The vertices are also indexed by certain maximal cyclic subgroups of $\Gamma$. The assumption that $M$ contains no $\log7$-semithick point is shown to imply that each vertex of $M$ has a contractible link, which is in turn used to show that the space $|K|-|K^{(0)}|$ obtained from the underlying space of $K$ by removing its vertices, is  contractible.

For each (open) simplex $\sigma$ of $K$, the maximal cyclic subgroups that index the vertices of $\sigma$ generate a subgroup $\Theta(\sigma)$ of $K$. Using $4$-freeness, together with the case $k=4$ of the so-called $\log(2k-1)$ Theorem, 
one shows that $|K|-|K^{(0)}|$ is set-theoretically a disjoint union $X_2 \cup X_3$, where $X_r$ denotes the union of all 
(open)
simplices $\sigma$ for which $\Theta(\sigma)$ has rank $r$.

(The $\log(2k-1)$ Theorem, which was proved under restrictive hypotheses in \cite{accs} and later generalized using the results of \cite{agol} and \cite{cg} among others, asserts that if orientation-preserving isometries $x_1,\dots,x_k$ of $\HH^3$ freely generate a free, discrete group of isometries, then for every point $P$ of $\HH^3$ we have $\dist(P,x_i\cdot P)\ge\log(2k-1)$ for at least one index $i\in\{1,\ldots,k\}$. The proof of the theorem yields a stronger result which appears in the present paper as Theorem \ref{our log(2k-1)}; it asserts that under the same hypotheses we have $\sum_{i=1}^k1/(1+\exp(\dist(P,x_i\cdot P)))\le1/2$. In this introduction we shall refer to the latter result as the Strong $\log(2k-1)$ Theorem.)

The
connected
components of $X_2$ and of $X_3$ are in a natural bijective correspondence with the vertices of a bipartite graph $T$, and the contractibility of $|K|-|K^{(0)}|$ implies that $T$ is a tree. There is a natural action of $\Gamma$ on $T$, and somewhat involved arguments in combinatorial group theory show that the stabilizer of each vertex under this action is a locally free group (in the sense that its finitely generated subgroups are free). This yields a contradiction, because the fundamental group of a closed hyperbolic $3$-manifold cannot act on a tree in such a way that all vertex stabilizers are locally free.

It would be natural to try to extract stronger information from the argument sketched above under the assumption that $\pi_1(M)$ is $k$-free for a given $k>4$, by using $\log(2k-1)$ in place of $\log7$. The construction that we have described would exhibit $|K|-|K^{(0)}|$ as  a disjoint union $X_2 \cup\cdots\cup X_{k-1}$, which would give rise to a $(k-3)$-complex in place of the tree $T$. This would make it difficult to extract geometric information. (In \cite{gs} we gave a variant of the argument that does involve a tree and allows $k$ to be greater than $4$, but in place of a $\log(2k-1)$-semithick point it gave a point with a weaker property, which seems to be difficult to use in obtaining volume estimates.)

Theorem \ref{before key} of the present paper is a generalization of \cite[Theorem 1.4]{fourfree}, and plays the role in this paper that \cite[Theorem 1.4]{fourfree} plays in \cite{fourfree}. 
Theorem \ref{before key} is stated in terms of a function, denoted $\f_3$ in the body of the present paper, which may be defined directly by $\f_3(x)=\log((e^x+3)/(e^x-5))$ for $x>\log5$. (A more instructive definition is given in the body of the paper.) Then Theorem \ref{before key} asserts that if $M$ is a 
closed, 
orientable hyperbolic $3$-manifold such that
$\pi_1(M)$ is $4$-free, then for every real
number  $\lambda\ge\log7$, the manifold $M$ contains either a $\lambda$-semithick point or an
$\f_3(\lambda)$-thick point. For $\lambda=\log7$ this 
implies
\cite[Theorem 1.4]{fourfree} (see Corollary \ref{corollary before key} below). 
Applying 
Theorem \ref{before key} for suitable values of $\lambda>\log7$ gives stronger information.

The proof of Theorem \ref{before key}, like that of \cite[Theorem 1.4]{fourfree} begins by writing $M$ as $\HH^3/\Gamma$, where $\Gamma$ is a cocompact $4$-free group of orientation-preserving isometries of $\HH^3$. If we assume that $M$ contains no $\lambda$-semithick point, we can again deduce that the non-empty sets of the form $Z_{\lambda}(C)$ form a covering of $\HH^3$, whose nerve $K$ has the property that $|K|-|K^{(0)}|$ is  contractible. One can define $\Theta(\sigma)$ as above for any simplex $\sigma$ of $K$.

If we assume that $M$ has no $\f_3(\lambda)$-thick point, a surprising property of the complex $K$ can be deduced. We define (see Definition \ref{restricting def} below) a 
{\it $3$-\restraining\ vertex} 
of a simplex $\phi$  of $K_{{\mathcal Z}_{\lambda}(\Gamma)}$ to be a vertex $v$  of $\phi$ with the property that for every face $\tau$ of $\phi$ such that $v$ is a vertex of $\tau$, we have $\rank(\Theta(\tau))\le3$. According to
the case $r=3$ of
Proposition \ref{restricting prop} below, the assumption that $M$ has no $\f_3(\lambda)$-thick point implies that, for every simplex $\sigma$ of $K$, there is a simplex $\phi$ of $K$ such that (i) $\sigma$ is a face of $\phi$, and (ii) $\phi$ has at least one
$3$-restraining 
vertex. 
We wish to emphasize that the proof of Proposition \ref{restricting prop} is an application of the Strong $\log(2k-1)$ Theorem, whereas the fact used in the proof of \cite[Theorem 1.4]{fourfree} used only the $\log(2k-1)$ Theorem.
At the time that \cite{fourfree} was written, it was not clear how to use the Strong $\log(2k-1)$ Theorem in this context.

The consequences of the surprising property of the complex $K$ given by Proposition \ref{restricting prop} are expressed in terms of a closure operation that is defined for certain subgroups of 
a  group.
The theory of this closure operation, which we think may be of independent interest, is developed in a purely group-theoretical context in Section \ref{structure section} of the present paper. 
If $\Gamma$ is any group, we define a class of subgroups of $\Gamma$ called the {\it closable} subgroups. If $\Gamma$ is $k$-free for a given positive integer $k$, the class of closable subgroups consists entirely of subgroups of finite local rank, and contains the class of all subgroups of $\Gamma$ having local rank at most $k-1$. (To say that a group $A$ has {\it local rank} at most $r$ means that every finitely generated subgroup of $A$ is contained in a subgroup of rank at most $r$; for a more formal definition of local rank, see Subsection \ref{more refer to me}.) The closure operation assigns a closable subgroup $\clo(A)\ge A$ of $\Gamma$ to each closable subgroup $A$ of $\Gamma$; the operation has a number of formal properties in common with (for example) the operation that assigns to a subfield of a field its relative algebraic closure (see Proposition \ref{closure properties}). One crucial property of $\clo(A)$ is that it contains every element $x$ of $\Gamma$ with the property that the group $\langle A,x\rangle$ generated by $A$ and $x$ has at most the same local rank as $A$.

When the complex $K$ has the property expressed in the conclusion of Proposition \ref{restricting prop}, one can combine this property with the machinery of Section \ref{structure section} 
to deduce strong information about the local ranks of the subgroups $\clo(\Theta(\sigma))$, for simplices $\sigma$ of $K$. This information is expressed in Corollaries \ref{first restricting cor} and \ref{second restricting cor}, and are exploited in the proof of Theorem \ref{before key}. Corollary \ref{first restricting cor}  allows one to write $|K|-|K^{(0)}|$ as a set-theoretical disjoint union $Y_2\cup Y_3$, where $Y_r$ denotes the union of all simplices $\sigma$ for which $\clo(\Theta(\sigma))$ has local rank $r$. The same construction that was used in the proof of \cite[Proposition 1.4]{fourfree}, using the $Y_r$ in place of the $X_r$,
then gives a tree $T$ with a natural action of $\Gamma$; and by combining Corollaries \ref{first restricting cor} and \ref{second restricting cor} with more results from Section \ref{structure section}, one can show that under this action all vertex stabilizers are locally free, which again gives a contradiction.

We note that in the body of the paper we define $r$-restraining vertices for any positive integer $r$, and that the actual statement of Proposition \ref{restricting prop} is a very direct generalization of the special case quoted above for $r=3$, assuming that $\Gamma\cong\pi_1(M)$ is $(r+1)$-free. However, the generalization seems difficult to apply when $r>3$, for the same reasons, pointed out above, for which it is difficult to extract stronger information from the arguments of \cite{fourfree} under the assumption that  $\pi_1(M)$ is $k$-free for a given $k>4$.

One interesting feature of the proof of Theorem \ref{before key} is that it depends (via Proposition \ref{general klm} of this paper) on a result due to Kent and Louder-McReynolds (\cite[Theorem 2]{kent},
\cite[Theorem 1.1]{LMcR}) on ranks of joins and intersections of subgroups of a free group. This result was used in the original proof of \cite[Theorem 1.4]{fourfree}, but it was shown in \cite{gs} that an alternative proof of \cite[Theorem 1.4]{fourfree} could be given without it. For the proof of the stronger Theorem \ref{before key} of this paper, the use of the Kent-Louder-McReynolds result appears to be indispensable.

After the proof of Theorem \ref{before key}, which with needed background and relevant foundational material occupies most of Sections \ref{structure section}--\ref{central section}, the next step is to use the geometric information given by the conclusion of Theorem \ref{before key} to obtain lower bounds for hyperbolic volume. A substantial part of the paper \cite{fourfree} was devoted to the analogous task of deducing \cite[Theorem 1.6]{fourfree}, which gives a lower bound for the volume of a closed, orientable hyperbolic $3$-manifold with $4$-free fundamental group, from \cite[Theorem 1.4]{fourfree}, which asserts the existence of a $\log7$-semithick point in such a manifold. The very geometric arguments needed to carry out this task involved the Strong $\log(2k-1)$ Theorem and a refinement of B\"{o}r\"{o}czky's work \cite{boroczky} on sphere-packing in $\HH^3$.

In this paper, the transition from Theorem \ref{before key}, via its somewhat more technical corollary \ref{keyer corollary}, to the lower bounds for volume given by Propositions  \ref{four-free case} and  \ref{five-free case} depends on involved geometric arguments, which are qualitatively similar to those used to pass from \cite[Theorem 1.4]{fourfree} to \cite[Theorem 1.6]{fourfree}. These arguments occupy Sections \ref{triangle section}---\ref{new very short section} of the present paper.

The results of Section \ref{elem vol section} give lower bounds for the volumes of certain subsets of a hyperbolic $3$-manifold in terms of the geometry of the manifold; some of these are more systematic versions of results from Sections 6, 7, 8, 9 and 11 of \cite{fourfree}, but a number of them are  stronger than the corresponding results from \cite{fourfree}.
Section \ref{triangle section} is an  observation about hyperbolic triangles that is used in Section \ref{elem vol section}.
The results of Section  \ref{Margulis section} concerning volumes, diameters, and Margulis numbers are strict improvements over the corresponding results in Section 10 of \cite{fourfree}. The results of Section  \ref{new very short section}, which concern the case of a hyperbolic $3$-manifold containing a very short closed geodesic, are stronger and more general than those established by the corresponding arguments in \cite{fourfree} (which are contained in the last step of the proof of \cite[Lemma 13.4]{fourfree}).

In Section \ref{my favorite monster}, the results of Sections
% \ref{general volume section}, 
\ref{triangle section}---\ref{new very short section} are combined with Corollary \ref{keyer corollary} to prove Proposition \ref{latest monster} and its corollary \ref{monster corollary}, which for each integer $k\ge4$ provide  sufficient conditions, stated entirely in terms of ranges of analytically defined functions, for a given number $V_0$ to be a lower bound for the volumes of all closed, orientable hyperbolic $3$-manifolds with $k$-free fundamental group. In Section \ref{4 and 5} we prove a result, Lemma \ref{coffee},  which is a fairly direct consequence of Corollary \ref{monster corollary} but is better adapted to numerical calculations. In Lemma \ref{coffee} we resort to a brute-force partitioning method (which was also used in \cite{fourfree}) because we do not have analytical techniques for handling the functions involved. 
Propositions \ref{four-free case} and \ref{five-free case}, of which the statements were given above, are direct applications of Lemma \ref{coffee}.

The transition 
made in \cite{fourfree} from \cite[Theorem 1.6]{fourfree}, which relates $4$-freeness to volume, to \cite[Theorem 1.7]{fourfree}, which relates mod $2$ homology to volume, depended on the 
%involved deep geometric and 
topological results of
%considerations, which are embodied in the results of 
\cite{JS}, \cite{sw}, \cite{bs}, \cite{accs}, \cite{acs-singular} and \cite{singular-two}. These results show that if $M$ is a closed, orientable hyperbolic $3$-manifold for which $\dim H_1(M;\FF_2)$ is sufficiently large in comparison with a given integer $k$, then either $\pi_1(M)$ is $k$-free, or $M$ contains an (embedded) incompressible surface  of genus less than $k$. Using the results of \cite{AST}, one can show that if a closed, orientable hyperbolic $3$-manifold $M$ contains an incompressible surface whose genus is subject to a suitable upper bound, and if $\dim H_1(M;\FF_2)$ is subject to a suitable lower bound, then $\vol M$ is bounded below by $\voct=3.66\ldots$. In the presence of the hypothesis of \cite[Theorem 1.7]{fourfree} this is sufficient to give a contradiction.

The analogous
transition made in the present paper, from
estimates for $\vol M$ under the assumption that $\pi_1(M)$ is $k$-free for a given $k\ge4$
to the proofs of
Propositions \ref{good homology} and \ref{3.77 homology}, which relate mod $2$ homology to volume, uses the same basic approach based on \cite{AST}, but requires more difficult arguments than the ones used in \cite{fourfree}. This is because the information that we mentioned above---an upper bound for the genus of a suitably chosen incompressible surface, and a lower bound for $\dim H_1(M;\FF_2)$---gives only a lower bound of $\voct$ for $\vol M$, and this is not strong enough to  contradict the upper bounds of $3.69$ and $3.77$ that appear in Propositions \ref{good homology} or \ref{3.77 homology} respectively. It is therefore necessary to produce an incompressible surface with stronger properties than an upper bound for its genus. This is done by using the more delicate combinatorial  techniques that were developed in \cite{CDS}. The information that is needed involves the Euler characteristic of the ``guts'' or ``kishkes'' of the manifold obtained by splitting $M$ along a suitably chosen incompressible surface, and is summarized in Proposition \ref{really from CDS} of the present paper. It is interesting to note that the latter result is deduced almost formally from the results in \cite{CDS} but was overlooked when that paper was written.

The method used in Section \ref{homology section} to prove Proposition \ref{3.77 homology}, if it were applied with Proposition \ref{four-free case} in place of Proposition \ref{four-free case}, would be insufficient to prove Proposition \ref{good homology}; it would show only that if $\vol M\le 3.57$, we have $\dim H_1(M;\FF_2)\le7$. In order to prove Proposition \ref{good homology}, which involves the constant $3.69$ rather than $3.57$, we must supplement the ingredients in the proof of Proposition \ref{3.77 homology} with a Dehn drilling argument. It turns out that the case in which the methods described above do not give the information required for (the contrapositive of) Proposition \ref{good homology} is the one in which $\pi_1(M)$ is $4$-free, $\dim H_1(M;\FF_2)\ge8$, and $M$ contains a simple closed geodesic $c$ of length at most $0.5912$. In this case, if we regard $c$ as a simple closed curve in $M$, 
its complement is homeomorphic to a hyperbolic manifold $M_c$. Using that
$\dim H_1(M;\FF_2)\ge8>6$, we use
\cite[Theorem 6.2]{CDS} to deduce that $\vol M_c>5.06$. A result due to Agol and Dunfield, which appears as Theorem 10.1 of \cite{AST}, together with an application of the case $k=2$ of the Strong $\log(2k-1)$ Theorem, 
allows one to relate the volume of $M$ to that of $M_c$; this provides the needed lower bound for $\vol M$. The details of this argument occupy Section \ref{drilling section}.

We thank Marc Culler, Ilya Kapovich, and Jason DeBlois for enlightening discussions of material related to this paper. We are also grateful to Nathan Dunfield for providing the information on known volumes that was given earlier in the introduction, and to Steve Kerckhoff for providing us with the references \cite{kojima-geodesic} and \cite{sakai}.

\section{The structure of $k$-free groups}\label{structure section}

\DefinitionsRemarks\label{refer to me}
If $\Gamma$ is a group, we will write $A\le\Gamma$ to mean that $A$ is a subgroup of $\Gamma$.

We shall say that elements $x_1,\ldots,x_m$ of a group $\Gamma$ are {\it independent} if the subgroup $\langle x_1,\ldots,x_m\rangle$ of $\Gamma$ is free on its generators $x_1,\ldots,x_m$.

A group is said to be {\it locally free} if each of its finitely generated subgroups is free.
(This definition, and some of the ones below, were given in a less formal setting in the Introduction.)

The {\it rank} of a  group $\Gamma$, denoted $\rank \Gamma$,  is defined to  be the minimal cardinality of a generating set for $\Gamma$.
\EndDefinitionsRemarks

A group $\Gamma$ is said to be {\it $k$-free}, where $k$ is a  given
non-negative 
integer, if every subgroup of $\Gamma$ having rank at most $k$ is a free group.

Every group is $0$-free. If $k$ and $k'$ are integers with $0\le k'\le k$, every  $k$-free group is $k'$-free. We define the {\it index of freedom} of a group $\Gamma$, denoted $\iof(\Gamma)$, to be the supremum of all integers $k$ such that $\Gamma$ is $k$-free. Thus we have $\iof(\Gamma)=\infty$ if and only if $\Gamma$ is locally free, i.e. every finitely generated subgroup of $\Gamma$ is free. If $\Gamma$ is not locally free then $\iof(\Gamma)$ is the largest integer $k$ such that $\Gamma$ is $k$-free.

\DefinitionsRemarks\label{more refer to me}
Let $A$ be a group. If there is a non-negative integer $r$ with the property that every finitely generated subgroup of $A$ is contained in a subgroup of $A$ having rank $r$, then the smallest integer with this property will be called the {\it local rank} of $A$. If no such integer $r$ exists, we say that the local rank of $A$ is infinite.

Note that the local rank of a finitely generated group is equal to its rank. Note also that if $k$ is a positive integer and $\Gamma$ is a $k$-free group, then every subgroup of $\Gamma$ having local rank at most $k$ is a locally free group.
\end{definitionsremarks}

\Proposition\label{free prod local rank}
Let $A$ be a group of local rank $r<\infty$, and let $\langle t\rangle$ be an infinite cyclic group. Then the free product $A\star\langle t\rangle$ has local rank $r+1$.
\EndProposition

\Proof
If $A$ is finitely generated, so that $r=\rank A$, then $A\star\ltr$ is finitely generated, and $\localrank(A\star\ltr)=\rank(A\star\ltr)=r+1$ by Grushko's Theorem.

For the proof in the general case, first note that if $B$ is any finitely generated subgroup of $A\star\langle t\rangle$, we have $B\le C_0\star\langle t\rangle$ for some finitely generated subgroup $C_0$ of $A$. Since $\localrank( A)=r$, 
there is a finitely generated subgroup $C_1$ of $A$ with $\rank C_1=r$ and $C_0\le C_1$. Then $B\le C_1\star\ltr$, and $\rank(C_1\star\ltr)=r+1$. Hence $\localrank (A\star\ltr)\le r+1$. 

Now assume that the latter inequality is strict, and consider an arbitrary finitely generated subgroup $D$ of $A$.  Since we have assumed that $\localrank( A\star\ltr)\le r$, there is a finitely generated subgroup  $E$ of $A\star \ltr$ such that $\rank E\le r$ and $D\star\ltr\le E$. According to the Kurosh Subgroup Theorem, we may identify $E$ with a free product $G_1\star\cdots\star G_m\star F$, where $F$ is a free group of finite rank and $G_1,\ldots,G_m$ are precisely the non-trivial subgroups that arise as intersections of $E$ with conjugates of $A$. If we regard $D$ as a subgroup of $D\star\ltr$, we have  $D\le E\cap A$; hence $D\le G_i$ for some $i$, and after re-indexing we may assume  that $D\le G_1$. On the other hand, since $t\in E$, the subgroup $E$ of $A\star\ltr$ cannot be generated by its intersections with conjugates of $A$; hence $F$ has strictly positive rank. By Grushko's Theorem we have $r\ge\rank E=\rank G_1+\cdots+\rank G_m+\rank F\ge\rank G_1+\rank F>\rank G_1$, so that $\rank G_1\le r-1$. As $D$ was an arbitrary finitely generated subgroup of $A$, and $D\le G_1$, this shows that $\localrank( A)\le r-1$, a contradiction.
\EndProof

\NotationRemark\label{j-def}
Let $A$ be a subgroup of a group $\Gamma$, and let $x$ be an element of $\Gamma$. If $\langle t\rangle$ is the standard infinite cyclic multiplicative group, we will denote by $\iota _{A,x}$ the homomorphism from the free product $A\star\langle t\rangle$ to $\Gamma$ defined by taking $\iota _{A,x}|A$ to be the inclusion, and seting  $\iota _{A,x}(t)=x$.

Note that if $B$ is any subgroup of $A$ then $B\star\langle t\rangle$ is canonically identified with a subgroup of $A\star\langle t\rangle$, and that under this identification we have $\iota _{B,x}=\iota _{A,x}|(B\star\langle t\rangle)$.
\EndNotationRemark

\DefinitionsRemarks\label{closable def}
A  subgroup $A$   of a group $\Gamma$ 
will be said to be {\it closable} if (a) $\localrank(A)<\infty$ and (b) there is an element $z$ of $\Gamma$ such that $\localrank(\langle A,z\rangle)<\iof(\Gamma)$.

In the case where $\Gamma$ is locally free, so that 
$\iof(\Gamma)=\infty$, Condition (b) is redundant; in fact, in this case, Condition (a) implies that the inequality in Condition (b) holds for every $z\in\Gamma$.

Note that if $A$ is a subgroup of a group $\Gamma$ such that $\localrank(A)<\iof(\Gamma)$, then $A$ is closable because we may take $z=1$ in the definition given above. 
A subgroup $A$  of $\Gamma$ will be termed {\it strongly closable} if $\localrank(A)<\iof(\Gamma)$.
\EndDefinitionsRemarks

\Proposition\label{when the rank goes up}
Let $A$ be a  subgroup of a 
group $\Gamma$, and let $x$ be an element of $\Gamma$. If $\localrank(A)<\infty$, and if the homomorphism $\iota _{A,x}:A\star\langle t\rangle \to \Gamma$ is injective, then \linebreak $\localrank(\langle A,x\rangle)>\localrank( A)$. Conversely, if $A$ is strongly closable and \linebreak $\localrank(\langle A,x\rangle)>\localrank( A)$, then $\iota _{A,x}$ is injective.
\EndProposition

\Proof
Set $r=\localrank( A)$. 
Under the hypothesis of either of the two assertions we have $r<\infty$.

If $\iota _{A,x}$ is injective then $\langle A,x\rangle$ is isomorphic to $A\star\langle t\rangle$, where $\langle t\rangle$ is infinite cyclic, so that Proposition \ref{free prod local rank} gives $\localrank(\langle A,x\rangle)=1+r>r$. This proves the first assertion.

To prove the 
second assertion, assume that $A$ is strongly closable, so that $r<\iof(\Gamma)$; and assume that $\localrank(\langle A,x\rangle)>r$.
Let us first consider the case in which $A$ is finitely generated. In this case, we have $\rank A=r$ and  $\rank\langle A,x\rangle>r$.
Since $r<\iof(\Gamma)$, the rank-$r$ group $A$ is free. Fix a basis $(u_1,\ldots,u_r)$ for $A$. Then $\{u_1,\ldots,u_r,x\}$ is a generating set for $\langle A,x\rangle$; in particular the rank of $\langle A,x\rangle$ is at most $r+1$. Since this rank has been 
seen 
to be strictly greater than $r$, it must be exactly $r+1$. Now since $r+1\le \iof(\Gamma)$, the rank-$(r+1)$ group $\langle A,x\rangle$ is free. 
But a generating set for a finite-rank free group, whose cardinality is equal to the rank of the group, must be a basis (see  \cite[vol. 2, p. 59]{kurosh}). It follows that the
generating set $\{u_1,\ldots,u_r,x\}$  for $\langle A,x\rangle$ must be a basis for $\langle A,x\rangle$.
 Now since $(u_1,\ldots,u_r)$ and $(u_1,\ldots,u_r,x)$ are bases for $A$ and $\langle A,x\rangle$ respectively, $\iota _{A,x}$ is injective.

For the proof of the second assertion in the general case, 
we 
must show that for an arbitrary element $z$ of $A\star\ltr$ we have $\iota _{A,x}(z)\ne1$. The hypothesis $\localrank(\langle A,x)\rangle
>r$ gives a finitely generated subgroup $B$ of $\langle A,x\rangle$ such that every finitely generated subgroup of $\langle A,x\rangle$ containing $B$ has rank at least $r+1$. Since $\langle B,z\rangle\le\langle A,x\rangle$ is finitely generated, we may choose a finitely generated subgroup $C$ of $A$ so that $\langle B,z\rangle\le\langle C,x\rangle$. Since $\localrank( A)=r$, there is a finitely generated subgroup $D$ of $A$ with $C\le D$ and $\rank D\le r$. We have $B\le\langle C,x\rangle\le\langle D,x\rangle$, and our choice of $B$ therefore guarantees that $\rank\langle D,x\rangle\ge r+1$. 

Hence $\rank\langle D,x\rangle >\rank D$, and $\rank D\le r<\iof(\Gamma)$. Since we have already proved the assertion in the case of a finitely generated subgroup, we may apply this special case with $D$ playing the role of $A$ to deduce that $\iota _{D,x}:D\star\langle t\rangle \to \Gamma$ is injective. But according to \ref{j-def} we have $\iota _{D,x}=\iota _{A,x}|(D\star\langle t\rangle)$. As $z\in C\le D$, we have 
$\iota _{A,x}(z)=\iota _{D,x}(z)\ne1$, as required.
\EndProof

\Proposition\label{new if the right hand don't get ya}
Let $k$ be a strictly positive integer, let  $\Gamma$ be a $k$-free group,  let $S$ be a finite subset of $\Gamma$ with 
$\rank \langle S\rangle\ge k$, 
and let $x_0$ be a non-trivial element of $S$. Then there are elements $x_1,\ldots,x_{k-1}$ of $S$ such that $x_0,\ldots,x_{k-1}$ are independent.
\EndProposition

\Proof
Write $S=\{z_1,\ldots,z_N\}$, where $N$ is the cardinality of $S$ and $z_1=x_0$. For $i=1,\ldots,N$ set $r_i=\rank\langle z_1,\ldots, z_i\rangle$. We have $r_1=1$ since $x_0\ne1$, and $r_N=\rank \langle S\rangle\ge k$ by hypothesis. Furthermore, for $1<i\le N$ we have $r_i\le r_{i-1}+1$. (There may well be values of $i$ for which $r_i<r_{i-1}$.) 
Hence if $J$ denotes the set of all indices $i\in\{2,\ldots,N\}$ for which $r_i= r_{i-1}+1$ (or equivalently for which $r_i> r_{i-1}$), then the cardinality of $J$ is at least $k-1$. 
Let $i_1,\ldots,i_{k-1}$ denote the $k-1$ smallest indices in $J$, labeled in such a way that 
$1<i_1<\cdots<i_{k-1}$. Then we have $r_{i_j}=r_{i_j-1}+1\le k$ for $j=1,\ldots,k-1$. We set $x_j=z_{i_j}$ for $j=1,\ldots,k-1$.

With $x_0$ given by the hypothesis, and with $x_1,\ldots,x_{k-1}\in S$ defined as above, we will show by induction for $t=0,\ldots,k-1$ 
that the elements $x_0,\ldots,x_t$ are independent. For $t=k-1$ this gives the conclusion of the proposition.  For $t=0$ the assertion amounts to saying that $x_0$ has infinite order; this is true because $x_0\ne1$ by hypothesis, and $\Gamma$ is torsion-free because it is $k$-free and $k>0$.

Now suppose that we are given an integer $t$ with $0<t<k$, and that $x_0,\ldots,x_{t-1}$ are independent. We will apply the second assertion of Proposition \ref{when the rank goes up}, taking $A=\langle z_1,\ldots,z_{i_t-1}\rangle$ and $x=x_t=z_{i_t}$. In view of the definitions given above, we have $\rank A=r_{i_t-1}$ and $\rank(\langle A,x\rangle)=r_{i_t}$. Since $r_{i_t}=r_{i_t-1}+1\le k$, we have $\rank(\langle A,x\rangle)>\rank A$, and $\rank A<k$; since $\Gamma$ is $k$-free the latter inequality  implies that $\rank A<\iof(\Gamma)$, i.e. $A$ is strongly closable in $\Gamma$.  It therefore follows from the second assertion of Proposition \ref{when the rank goes up} that $\iota_{A,x}: A\star\langle t\rangle \to \Gamma$ is injective.

Now set $B=\langle x_0,x_1,x_2,\ldots,x_{t-1}\rangle=\langle z_1,z_{i_1},z_{i_2},\ldots,z_{i_{t-1}}\rangle\le\langle z_1,z_2,\ldots,z_{i_{t-1}}\rangle=A$. According to \ref{j-def}, $\iota _{B,x}$ is the restriction of $\iota _{A,x}$ to the subgroup $B\star\langle t\rangle$ of $A\star\langle t\rangle$. Hence $\iota _{B,x}=\iota _{B,x_t}$ is injective. The independence of $x_0,x_1,x_2,\ldots,x_{t-1}$ means that $B$ is free on the generators $x_0,x_1,x_2,\ldots,x_{t-1}$. The injectivity of $\iota _{B,x_t}$ then implies that $\langle B,x_t\rangle$ is free on the generators $x_0,x_1,x_2,\ldots,x_{t}$; that is, $ x_0,x_1,x_2,\ldots,x_{t}$ are independent. This completes the induction.
\EndProof
%P

\Definition\label{admissible elt}
Let $A$ be a   subgroup  of a group $\Gamma$.
An element $x$ of $\Gamma$ will be said to be {\it $A$-\admissible} if 
(a) $\langle A,x\rangle $ is strongly closable and
(b) $\localrank(\langle A,x\rangle)\le\localrank(A))$.
% and
% $\localrank(\langle A,x\rangle)<\iof(\Gamma)$. z

(Recall that the condition that $\langle A,x\rangle $ is strongly closable means that $\localrank(\langle A,x\rangle)<\iof(\Gamma)$.  Hence, in the case where $\Gamma$ is not locally free, so that $\iof(\Gamma)<\infty$, Conditions (a) and (b) above may be expressed 
more simply
by the single inequality
$\localrank(\langle A,x\rangle)\le\min(\iof(\Gamma)-1,\localrank(A))$.

Although the ultimate applications in this paper involve the case where $\Gamma$ is not locally free, we have found it more natural to deal with a definition that applies to all groups.)
\EndDefinition

\Proposition\label{red river}
Let $A$ be a strongly closable subgroup of a group $\Gamma$. For any element $x$ of $\Gamma$, the following conditions are equivalent:
\begin{enumerate}
\item $x$ is $A$-\admissible;
\item $\localrank(\langle A,x\rangle)\le\localrank(A)$;
\item the homomorphism $\iota_{A,x}: A\star\langle t\rangle \to \Gamma$ is non-injective. 
\end{enumerate}
In particular, every element of $A$ is $A$-\admissible.
\EndProposition

\Proof
The equivalence of (2) and (3) follows from Proposition \ref{when the rank goes up}. The implication $(1)\Rightarrow(2)$ is trivial, because Condition (2) is Condition (b) of the  definition of $A$-\admissibility. To prove $(2)\Rightarrow(1)$, note that the strong closability of $A$ means that $\localrank(A)<\iof(\Gamma)$. Hence if $\localrank(\langle A,x\rangle)\le\localrank(A)$, then $\localrank(\langle A,x\rangle) <\iof(\Gamma)$, so that $x$ satisfies Condition (a) of the  definition of $A$-\admissibility\ as well as Condition (b). This proves the first assertion. The second assertion follows, because if $x$ is an element of $A$, then $\langle A,x\rangle=A$, and hence Condition (2) holds.
\EndProof

\Proposition\label{little observation}
Let 
$A$ be a   subgroup  of a group $\Gamma$, 
with $\localrank(A)<\infty$. Then $\Gamma$ contains an $A$-\admissible\ element if and only if $A$ is closable. 
\EndProposition

\Proof
If $x$ is an $A$-admissible element of $\Gamma$, 
then according to Condition (a) of Definition \ref{admissible elt}, $\langle A,x\rangle $ is strongly closable. By definition (see \ref{closable def}) this means that $\localrank(\langle A,x\rangle)
<\iof(\Gamma)$.
This gives Condition (b)  of Definition \ref{closable def}, and Condition (a) of Definition \ref{closable def} holds by hypothesis. Hence $A$ is closable.

Conversely, suppose that $A$ is closable. In the case  where $A$ is strongly closable, it follows from the 
final assertion of Proposition \ref{red river} that $1\in A$ is an $A$-\admissible\ element. 
Now consider the case in which $A$ is not strongly closable, i.e. $\localrank(A)\ge\iof(\Gamma)$. In particular,  $\iof(\Gamma)$ is then finite. Since $A$ is closable, Condition (b) of Definition \ref{closable def}  gives an element $z$ of $\Gamma$ such that $\localrank(\langle A,z\rangle)<\iof(\Gamma)$. This means that $\langle A,z\rangle$ is strongly closable, which is Condition (a) of Definition \ref{admissible elt}. Since $\iof(\Gamma)\le\localrank(A)$, we have $\localrank(\langle A,z\rangle)<\localrank(A)$, which implies  Condition (b) of Definition \ref{admissible elt}. Hence $z$  is an $A$-\admissible\ element in this case.
\EndProof

\Proposition\label{when the rank of a subgroup goes up...}
Let $A$ and $B$ be   subgroups of a 
group 
$\Gamma$. Suppose that $B\le A$, that $\localrank(B)<\infty$, and that $A$ is strongly closable. 
Then every $B$-\admissible\ element of $\Gamma$ is $A$-admissible.
\EndProposition

\Proof
Let $x$ be a $B$-admissible element of $\Gamma$. According  to Proposition \ref{red river}, $x$ is $A$-admissible if and only if
the homomorphism
$\iota _{A,x}:A\star\langle t\rangle \to \Gamma$ is non-injective. Assume to the contrary that
$\iota _{A,x}$ is
injective. According to \ref{j-def} we have $\iota _{B,x}=\iota _{A,x}|(B\star\langle t\rangle)$. Hence $\iota _{B,x}$ is injective. According to 
the first assertion of
Proposition \ref{when the rank goes up}, this implies that $\localrank(\langle B,x\rangle)>\localrank( B)$.
But since $x$ is $B$-admissible, Condition (b) of Definition \ref{admissible elt} gives $\localrank(\langle B,x\rangle)\le\localrank( B)$. This is
 a contradiction.
\EndProof

\DefinitionRemark\label{what's admissible mean}
Let $A$ be a  subgroup of a 
group 
$\Gamma$,
with $\localrank(A)<\infty$. 
Let $x_1,\ldots,x_m$
be a finite sequence of elements of $\Gamma$,
with $m>0$. 
Let us set $A_0=A$, and $A_i=\langle A,x_1,\ldots,x_i\rangle$ for $i=1,\ldots,m$
(so that all the $A_i$ have finite local rank). 
We will say that $x_1,\ldots,x_m$ is an {\it $A$-\admissible\ sequence} in $\Gamma$ if for $i=1,\ldots,m$ the element $x_i$ is $A_{i-1}$-\admissible.

Note that if $x_1,\ldots,x_m$ is an $A$-admissible sequence in $\Gamma$, then for any $s$ with $1\le s\le m$, the sequence $x_1,\ldots,x_s$ is also $A$-admissible.

Note also that if $x_1,\ldots,x_m$ is an $A$-admissible sequence, and if the subgroups $A_0,\ldots,A_m$ are defined as above, we have $A_i=\langle A_{i-1},x_i\rangle$ for $i=1,\ldots,m$. Since $x_i$ is $A_{i-1}$-\admissible, it follows from 
Condition (a) of the definition of \admissibility\ (\ref{admissible elt})
that $A_1,\ldots,A_m$ are strongly closable subgroups of $\Gamma$.
\EndDefinitionRemark

\Lemma\label{still admissible}
Let $A$ and $B$ be subgroups of a 
group 
$\Gamma$.
Suppose that $\localrank(B)<\infty$, 
that $B\le A$, and that $A$ is strongly closable.
Then every $B$-admissible sequence in $\Gamma$ is also $A$-admissible.
\EndLemma

\Proof
Let  $x_1,\ldots,x_m$ be a $B$-admissible sequence. Let us set $B_0=B$ and $A_0=A$; and for $i=1,\ldots,m$ let us set $B_i=\langle B,x_1\ldots,x_i\rangle$ and  $A_i=\langle A,x_1\ldots,x_i\rangle$. Then 
for $i=0,\ldots,m$ we have $B_i\le A_i$, and the local ranks of $A_i$ and $B_i$ are finite. 
Furthermore, the $B$-admissibility of the sequence $x_1,\ldots,x_m$ means, by definition, that $x_i$ is $B_{i-1}$-admissible for $i=1,\ldots,m$.

We will prove  by induction for $i=1,\ldots,m$ that
$x_i$ is $A_{i-1}$-admissible. This will show that $x_1,\ldots,x_m$ is $A$-admissible.

For the base case $i=1$, first note that $A_0$ is strongly closable
by hypothesis. 
Next note that since $x_1$ is $B_0$-admissible, and since by hypothesis  $A_0$ is strongly closable and $B_0\le A_0$, it follows from Proposition \ref{when the rank of a subgroup goes up...}
that $x_1$ is $A_0$-admissible.

Now suppose that  for a given $i$ with $0\le i< m$, the element $x_{i}$ is $A_{i-1}$-admissible. Then according to 
Condition (a) of the definition of \admissibility\  (\ref{admissible elt}), 
the subgroup $A_i=\langle A_{i-1},x_{i}\rangle$ of $\Gamma$ is strongly closable. Next note that since $x_{i+1}$ is $B_i$-admissible, $B_i\le A_i$, and $A_i$ is strongly closable, it follows from Proposition \ref{when the rank of a subgroup goes up...}
that $x_{i+1}$ is $A_i$-admissible. This completes the induction.
\EndProof

\Lemma\label{concatenate}
Let $A$  be a closable subgroup of a 
group 
$\Gamma$,
and let $p$ be a positive integer. For $j=1,\ldots,p$ let $x^{( j )}_1,\ldots,x^{( j )}_{m_{ j }}$ be an $A$-admissible sequence in $\Gamma$. Then the sequence $$x^{(1)}_1,\ldots,x^{(1)}_{m_1},
x^{(2)}_1,\ldots,x^{(2)}_{m_2},\ldots,
x^{(p)}_1,\ldots,x^{(p)}_{m_p}$$ of elements of $\Gamma$ is $A$-admissible. 
\EndLemma

\Proof
It is enough to give the proof in the case $p=2$, as the general case then follows by induction. Thus it suffices to prove that if 
$m$ and $n$ are integers with $0<m<n$, and if $x_1,\ldots,x_n$ are elements of $\Gamma$ such that the sequences $x_1,\ldots,x_m$ and and $x_{m+1},\ldots,x_n$ are  $A$-admissible, then the sequence $x_1,\ldots,x_n$ is $A$-admissible. Let us  define subgroups $A_0,\ldots,A_n$ of $\Gamma$ by setting $A_0=A$, and $A_i=\langle A,x_1,\ldots,x_i\rangle$ for $i=1,\ldots,n$. Let us also define subgroups $B_m,\ldots,B_n$ of $\Gamma$ by setting $B_m=A$ and $B_i=\langle A,x_{m+1},\ldots,x_i\rangle$ for $m<i\le n$; thus we have $B_i\le A_i$ for $i=m,\ldots,n$. The $A_i$ and $B_i$ have finite local rank.
We must show that $x_i$ is $A_{i-1}$-admissible for $i=1,\ldots,n$. If $0\le i\le m$ this follows from the $A$-admissibility of $x_1,\ldots,x_m$. 

Now let $i$ be given with $m<i\le n$. In particular we have $i>1$. Arguing inductively, we may assume that $x_{i-1}$ is $A_{i-2}$-admissible. It follows from 
Condition (a) of the definition of \admissibility\ (\ref{admissible elt}), 
that $A_{i-1}$ is a strongly closable subgroup of $\Gamma$. We have $B_{i-1}\le A_{i-1}$, and the $A$-admissibility of $x_{m+1},\ldots,x_n$ implies that $x_i$ is $B_{i-1}$-admissible. We may therefore apply Proposition \ref{when the rank of a subgroup goes up...},
with $A_{i-1}$, $B_{i-1}$ and $x_i$ playing the respective roles of $A$, $B$ and $x$, to deduce that 
$x_i$ is $A_{i-1}$-admissible, as required.
\EndProof

\NotationRemarks\label{wuzza closure}
If $A$ is a closable subgroup of a 
group 
$\Gamma$, we will denote by $\clo(A)$ the set of all elements $x\in\Gamma$ with the property that there exists an $A$-admissible sequence $x_1,\ldots,x_m$ with $x_m=x$.

We observed in \ref{what's admissible mean} that if $x_1,\ldots,x_m$ is an $A$-admissible sequence then the sequence $x_1,\ldots,x_s$ is also $A$-admissible for any $s$ with $1\le s\le m$. It follows  that if $x_1,\ldots,x_m$ is an $A$-admissible sequence, then we have $x_1,\ldots,x_m\in\clo(A)$.

Note that, a priori, $\clo(A)$ is only a subset of $\Gamma$. 
If we were to use the definition given above in the case of a subgroup $A$ of $\Gamma$ which is not closable, it would follow from Proposition \ref{little observation} that $\clo(A)$ was the empty set. In contrast, 
Proposition \ref{closure properties} below will assert, among other things, that it is a subgroup of $\Gamma$; it will be seen from the proof of Proposition \ref{closure properties} that all the asserted properties of $\clo(A)$ depend on the closability of $ A$. It is for this reason that we have defined $\clo(A)$ only when $A$ is closable. 
\EndNotationRemarks

\Lemma\label{before closure properties I}
Let    $A$ be a closable subgroup of a 
group 
$\Gamma$, and let $S$ be a finite subset of $\clo(A)$. Then there is an $A$-admissible sequence $y_1,\ldots,y_q$ such that 
$\{y_1,\ldots,y_q\}\supset S$. 
\EndLemma

\Proof
Write $S=\{x^{(1)},\ldots,x^{(p)}\}$. For each $j\in\{1,\ldots,p\}$, since $x^{(j)}\in\clo(A)$, there is an $A$-admissible sequence $x^{( j )}_1,\ldots,x^{( j )}_{m_{ j }}$ with $x^{( j )}=x^{( j )}_{m_{ j }}$. 
According to Lemma \ref{concatenate}, the sequence 
\Equation\label{that sequence}
x^{(1)}_1,\ldots,x^{(1)}_{m_1},
x^{(2)}_1,\ldots,x^{(2)}_{m_2},\ldots,
x^{(p)}_1,\ldots,x^{(p)}_{m_p}
\EndEquation
is $A$-admissible. 
The conclusion therefore follows if we set $q=m_1+\cdots+m_p$ and define $y_1,\ldots,y_q$ to be the sequence (\ref{that sequence}). 
\EndProof

\Lemma\label{before closure properties II}
Let    $A$ be a closable subgroup of a 
group 
$\Gamma$, and let $y_1,\ldots,y_q$ be an $A$-admissible sequence. Then
$\localrank(\langle A,y_1,\ldots,y_q\rangle)\le\localrank( A)$, and
$\langle A,y_1,\ldots,y_q\rangle$ is a strongly closable subgroup of $\Gamma$.
\EndLemma

\Proof
Set $A_0=A$, and for $i=1,\ldots,q$ set $A_i=\langle A,y_1,\ldots,y_i\rangle$. By the definition of an $A$-admissible sequence, $y_i$ is $A_{i-1}$-\admissible\ for $i=1,\ldots,q$. Condition (b) of Definition \ref{admissible elt} gives $\localrank(A_{i})\le\localrank(A_{i-1})$ for $i=1,\ldots,q$; hence $\localrank(A_{q})\le\localrank(A_{0})$, which is the first assertion of the lemma. For $i=q$, Condition (a) of Definition \ref{admissible elt} gives that $A_q=\langle A_{q-1},x_q\rangle$ is strongly closable, which is the second assertion of the lemma.
\EndProof

\Proposition\label{closure properties}
Let $A$ be a closable subgroup of a 
group 
$\Gamma$.
Then:
\begin{enumerate}
\item  $\clo(A)$ is a strongly closable subgroup of $\Gamma$, and 
$\localrank(\clo(A)) \le\localrank(A)$; 
\item $A\le\clo(A)$;
\item for any 
closable subgroup  
$B$ of $A$
we have $\clo(B)\le\clo(A)$; and
\item 
$\clo(\clo(A))$ (which is defined since $\clo(A)$ is 
a closable subgroup of $\Gamma$ 
by Assertion (1)) is equal to $\clo(A)$.
\end{enumerate}
\EndProposition

\Proof
We first  show
that the set $\clo(A)$ contains $A$. 
Let an element  $a$ of $A$ be given. According to 
Proposition \ref{little observation}, 
$\Gamma$ contains an $A$-\admissible\ element $z$; and according to 
Condition (a) of the definition of \admissibility\ (\ref{admissible elt}), 
the group $A':=\langle A,z\rangle$ is strongly closable. Since $a\in A\le A'$, we may apply the 
final assertion of Proposition \ref{red river}, 
with the strongly closable subgroup $A'$ playing the role of $A$, to deduce that $a$ is $A'$-admissible. It now follows from the definition that the two-term sequence $z,a$ is an $A$-admissible sequence, which shows that $a\in\clo(A)$. Thus $A\subset\clo(A)$.

 Next we prove that $\clo(A)$ is a subgroup of $\Gamma$. 
Since we have shown $A\subset\clo(A)$, we have $1\in\clo(A)$. It therefore suffices to prove that for any $x,y\in\clo(A)$ we have $xy^{-1}\in\clo(A)$. 

If $x,y\in\clo(A)$ are given, then by definition there are $A$-admissible sequences $x_1,\ldots,x_m$ and $y_1,\ldots,y_n$ such that $x_m=x$ and $y_n=y$. According to Lemma \ref{concatenate}, the sequence $x_1,\ldots,x_m,y_1,\ldots,y_n$ is admissible. 

Now set $A'=\langle A,x_1,\ldots,x_m,y_1,\ldots,y_n\rangle$. According to (the final sentence of) \ref{what's admissible mean}, 
the $A$-\admissibility\ of $x_1,\ldots,x_m,y_1,\ldots,y_n$ implies that
$A'$ is a strongly closable subgroup of $\Gamma$. But 
we have $xy^{-1}=x_my_n^{-1}\in A'$; hence by  
Condition (a) of the definition of \admissibility\ (\ref{admissible elt}), 
$xy^{-1}$ is $A'$-\admissible. Now the 
$A'$-\admissibility\ of $xy^{-1}$, together with
the $A$-\admissibility\ of $x_1,\ldots,x_m,y_1,\ldots,y_n$ implies by definition that
$x_1,\ldots,x_m,y_1,\ldots,y_n,xy^{-1}$ is $A$-\admissible. Hence $xy^{-1}\in\clo(A)$. This 
completes the proof 
that $\clo(A)$ is a subgroup of $\Gamma$.

Now that $\clo(A)$ has been shown to be a subgroup, the inclusion $A\subset\clo(A)$, which was proved above, may be rewritten in the form $A\le\clo(A)$. This proves Assertion (2).

To 
complete the proof of 
Assertion (1), it remains to prove that 
the subgroup $\clo(A)$ is strongly closable and that
$\localrank(\clo(A)) \le\localrank(A)$. 
To say that $\clo(A)$ is strongly closable means by definition that $\localrank(A)<\iof(\Gamma)$. Thus if we set $r=\localrank(A)$, and define an integer $h$ by setting $h=\min(r, \iof(\Gamma)-1)$ if $\iof(\Gamma)<\infty$, and $h=r$ if $\iof(\Gamma)=\infty$, we must show that $\localrank(\clo(A)) \le h$.

Suppose that $B$ is a finitely generated subgroup of $\clo(A)$, and fix a finite generating set $S$ for $B$. 
According to Lemma \ref{before closure properties I}, there is an $A$-admissible sequence $y_1,\ldots,y_q$ such that $\{y_1,\ldots,y_q\}\supset S$. According to 
\ref{wuzza closure}, we have $\{y_1,\ldots,y_q\}\subset\clo(A)$; and according to
Lemma \ref{before closure properties II}, the group $C:=\langle A,y_1,\ldots,y_q\rangle$ 
is strongly closable and has local rank at most $ r$. 
Since strong closability means that $\localrank(C)<\iof(\Gamma)$, we have $\localrank(C)\le h$.
Since $S\subset\{y_1,\ldots,y_q\}\subset C$, we have
$B\le C$. Hence $B$ is contained in some finitely generated subgroup $B'$ of $C$ with $\rank B'\le 
h
$.

Since $\{y_1,\ldots,y_q\}\subset 
%C$,the elements $x^{(j)}_{i}$ lie in $
\clo(A)$, and since $A\le\clo(A)$ by Assertion (2), we have $C\le \clo(A)$. Hence $B'$ is in particular a finitely generated subgroup of $\clo(A)$ which contains $B$ and has rank at most 
$h$. 
As  $B$ was an arbitrary finitely generated subgroup of $\clo(A)$, this shows that $\localrank(\clo(A)) \le 
h
$, and the proof of Assertion (1) is complete.

To prove Assertion (3), consider an arbitrary 
closable subgroup  $B$ of $A$, and an arbitrary element $x$ of $\clo(B)$. 
Let us fix a $B$-\admissible\ sequence $x_1,\ldots,x_m$ with $x_m=x$.  Since $A$ is closable, it follows from Proposition \ref{little observation} that $\Gamma$ contains an $A$-\admissible\ element $z$. According to 
Condition (a) of the definition of \admissibility\ (\ref{admissible elt}), 
the subgroup $A':=\langle A,z\rangle$ of $\Gamma$ is strongly closable. Since $B\le A\le A'$, we may apply Lemma \ref{still admissible} to deduce that $x_1,\ldots,x_m$ is $A'$-\admissible. Now since the element $z$ is $A$-\admissible, and the sequence $x_1,\ldots,x_m$ is $A'$-admissible, where $A'=\langle A,z\rangle$, it follows from the definitions that the sequence $z,x_1,\ldots,x_m$ is $A$-admissible. This shows that $x=x_m\in\clo(A)$, and completes the proof of Assertion (3).

To prove Assertion (4), first note that since 
$\clo(A)$ is a closable subgroup of $\Gamma$
by Assertion (1), we may apply Assertion (2), with $\clo(A)$ playing the role of $A$, to deduce that $\clo(A)\le\clo(\clo(A))$. The main step in the proof of the reverse inclusion will be the proof of the following fact:

\Claim\label{da fact}
Every $\clo(A)$-\admissible\ element of $\Gamma$ lies in $\clo(A)$.
\EndClaim

To prove \ref{da fact}, suppose that $x$ is a $\clo(A)$-\admissible\ element of $\Gamma$. According to Assertion (1), which has already been proved, the subgroup $\clo(A)$ of $\Gamma$ is strongly closable. Hence we may invoke Proposition \ref{red river}, with $\clo(A)$ playing the role of $A$ in that proposition, to deduce that the homomorphism $\iota _{\clo(A),x}:\clo(A)\star\langle t\rangle \to \Gamma$ is non-injective, where $\ltr$ denotes an infinite cyclic group. Choose a non-trivial element $w$ of the kernel of $\iota _{\clo(A),x}$. Let us write $w=u_1t^{e_1}\cdots u_pt^{e_p}$, where $u_1,\ldots,u_p$ are elements of $\clo(A)$ and $e_1,\ldots,e_p$ are integers. (It is unimportant for this argument whether some of the $u_i$ are trivial or some of the $e_i$ are zero.) 
Set $U=\{u_1,\ldots,u_p\}$. 
According to Lemma \ref{before closure properties I}, applied with $U$ playing the role of $S$, there is an $A$-admissible sequence $y_1,\ldots,y_q$ such that $\{y_1,\ldots,y_q\}\supset U$. 
According to 
\ref{wuzza closure}, we have $\{y_1,\ldots,y_q\}\subset\clo(A)$; and according 
to Lemma \ref{before closure properties II},  the group $C:=\langle A,y_1,\ldots,y_q\rangle$ 
is strongly closable. 
Since $U\subset\{y_1,\ldots,y_q\}\subset C$, we have $w\in C\star\ltr$. Since $\{y_1,\ldots,y_q\}\subset\clo(A)$, and since Assertions (1) and (2) give that $\clo(A)$ is a subgroup of $\Gamma$ containing $A$, we have $C\le\clo(A)$. According to  \ref{j-def}, we have $\iota_{C,x}=\iota _{\clo(A),x}|(C\star\langle t\rangle)$. Since $C$ contains the non-trivial element 
$w$ of the kernel of $\iota_{\clo(A),x}$, the homomorphism $\iota_{C,x}$ is non-injective. Applying Proposition \ref{red river}, with the strongly closable subgroup
%follows from 
%Lemma \ref{when the rank goes up} 
%(applied with 
$C$ playing the role of  $A$ in that lemma, we deduce that $x$ is $C$-\admissible.
Now since the sequence $y_1,\ldots,y_q$  is $A$-\admissible, and since $x$ is $C$-\admissible, where $C=\langle A,y_1,\ldots,y_q\rangle$, it follows from the definitions that $y_1,\ldots,y_q,x$  is $A$-\admissible, and hence that $x\in\clo(A)$. This proves \ref{da fact}.

To complete the proof of (4), suppose that $x\in\clo(\clo(A))$ is given. Then there is a $\clo(A)$-\admissible\ sequence $x_1,\ldots,
x_m$ with $x=x_m$. Set $D_0=\clo(A)$, and for $i=1,\ldots,m$ set $D_i=\langle\clo(A),x_1,\ldots,x_i\rangle$. By induction on $i=1,\ldots,m$ we will show that $D_i=\clo(A)$. For $i=0$ the assertion is trivial. Now assume that $i\in\{1,\ldots,m\}$ is given, and that $D_{i-1}=\clo(A)$. The 
$\clo(A)$-\admissibility\ of $x_1,\ldots,
x_m$ implies that $x_i$ is $D_{i-1}$-admissible, i.e. that it is $\clo(A)$-admissible. According to \ref{da fact} we then have $x_i\in\clo(A)$; hence $D_i=\langle D_{i-1},x_i\rangle=\langle \clo(A),x_i\rangle=\clo(A)$. This completes the induction. Taking $i=m$ we find that $x=x_m\in A_m=\clo(A)$. This shows that $\clo(\clo(A))\le\clo(A)$, and completes the proof of Assertion (4).
\EndProof

Note that  Proposition \ref{closure properties} may be regarded as meaning that the operation $A\mapsto\clo(A)$ is a ``closure operation'' on the class of all closable subgroups of a given 
group, 
analogous, for example, to the operation that assigns to each subfield of a field $K$ its relative algebraic closure in $K$. (One difference is that whereas the class of subfields of a field $K$ includes $K$ itself, the class of closable subgroups of a 
group $\Gamma$ does not in general include $\Gamma$.)

\Corollary\label{properties cor} 
Let $A$ and $B$ be closable subgroups of a group
$\Gamma$, and suppose that $A\le B\le\clo(A)$. Then $\clo(B)=\clo(A)$.
\EndCorollary

\Proof
Since $A\le B$, it follows from Assertion (3) of Proposition \ref{closure properties} (with the roles of $A$ and $B$ reversed) that $\clo(A)\le\clo(B)$. By Assertion (1) of Proposition \ref{closure properties}, $\clo(A)$ is a closable subgroup of $\Gamma$. Since $B\le \clo(A)$, it follows from Assertion (3) of Proposition \ref{closure properties} (with $\clo(A)$ playing the role of $A$) that $\clo(B)\le\clo(\clo(A))$. But by Assertion (4) of Proposition \ref{closure properties} we have $\clo(\clo(A))=\clo(A)$, and hence 
$\clo(B)\le\clo(A)$.
\EndProof

\Proposition\label{and the normalizer}
Let $A$ be a closable subgroup of a 
group 
$\Gamma$.
Let $B$ be a non-trivial subgroup of $A$. Then the normalizer of $B$ in $\Gamma$ is contained in $\clo(A)$.
\EndProposition

\Proof
According to  
Proposition \ref{little observation}, 
$\Gamma$ contains an $A$-\admissible\ element $z$; and according to 
Condition (a) of the definition of \admissibility (\ref{admissible elt}), 
the group $A':=\langle A,z\rangle$ is strongly closable. 
Let us fix a non-trivial element $b_0\in B$.

Given any element  $x$  of the normalizer of $B$, we consider the homomorphism $\iota_{A',x}: A'\star\langle t\rangle \to \Gamma$ (see \ref{j-def}). Since $x$ normalizes $B$, we have $b_1:=xb_0x^{-1}\in B$. In particular $b_0$ and $b_1$ are non-trivial elements of $A'$. Hence $b_1$ and $tb_0t^{-1}$ are elements of $A'\star\langle t\rangle$, which are respectively given by words of length $1$ and $3$ in the free product, and are therefore distinct. But we have $\iota_{A',x}(tb_0t^{-1})=xb_0x^{-1}=b_1=\iota_{A',x}(b_1)$. 
This shows that $\iota_{A',x}$ is non-injective. Since $A'$ is strongly closable, it now follows from Proposition \ref{red river} that $x$ is $A'$-\admissible. Since $z$ is $A$-\admissible\ and $A'=\langle A,z\rangle$, it now follows from the definitions that $z,x$ is a two-term $A$-\admissible\ sequence, and hence that $x\in\clo(A)$.
\EndProof

\Proposition\label{conjugation and closure}
Let  $A$ be a closable subgroup of a 
group $\Gamma$. 
Then for every element $u$ of $\Gamma$, the subgroup $uAu^{-1}$ of $\Gamma$ 
is closable, and $\clo(uAu^{-1})=u\clo(A)u^{-1}$.
\EndProposition

\Proof
Since $A$ is closable, Definition \ref{closable def} gives $\localrank(A)<\infty$, and $\localrank(\langle A,z\rangle)<\iof(\Gamma)$ for some $z\in\Gamma$.
We have $\localrank(uAu^{-1})=\localrank(A)<\infty$, and \linebreak
$\localrank(\langle uAu^{-1},z\rangle)=\localrank(\langle A,z\rangle)<\iof(\Gamma)$, so that $uAu^{-1}$ is closable.
Now suppose 
that $x\in\clo(A)$ is given, and fix an $A$-\admissible\ sequence $x_1,\ldots,x_m$ with $x_m=x$. Set $A_0=A$, and for $1\le i\le m$, set $A_i=\langle A,x_1,\ldots,x_i\rangle$. Then  we have $uA_0u^{-1}=uAu^{-1}$, and
$uA_iu^{-1}=\langle uAu^{-1},ux_1u^{-1},\ldots,ux_iu^{-1}\rangle$. Now for $i=1,\ldots,m$ the element $x_i$ is $A_{i-1}$-\admissible; according to Definition \ref{admissible elt} this means that $\localrank(A_i)<\iof(\Gamma)$ and that
$\localrank(A_i)\le\localrank(A_{i-1})$. Hence
$\localrank(uA_iu^{-1})=\localrank(A_i)<\iof(\Gamma)$ and 
$\localrank(uA_iu^{-1})=\localrank(A_i)\le\localrank(A_{i-1})=\localrank(uA_{i-1}u^{-1})$. This shows 
that
$ux_iu^{-1}$ is  $uA_{i-1}u^{-1}$-\admissible\ for $i=1,\ldots,m$; this says that $ux_1u^{-1},\ldots,ux_mu^{-1}$ is a $uAu^{-1}$-\admissible\ sequence, which implies that $uxu^{-1}=ux_mu^{-1}\in\clo(uAu^{-1})$. Thus the inclusion $u\clo(A)u^{-1}\le\clo(uAu^{-1})$ is established. As this holds for every $u\in\Gamma$ and every 
closable subgroup $A$ of $\Gamma$, 
we may replace $u$ by $u^{-1}$ and $A$ by $uAu^{-1}$ to deduce that $u^{-1}\clo(uAu^{-1})u\le\clo(A)$, i.e. $\clo(uAu^{-1})\le u\clo(A)u^{-1}$. 
\EndProof

The following mild generalization of a deep result due to Kent and Louder-McReynolds is different in flavor from the other results in this section. It will be used in a crucial way in the proof of Lemma \ref{commies}, which is in turn a key step in the proof of our central result, Theorem  \ref{key theorem}.

\Proposition\label{general klm}
Suppose that $A$ and $B$ are subgroups of a $4$-free group $\Gamma$. Suppose that $\localrank(A)=\localrank(B)=2$ and that the subgroup $A\cap B$ is non-abelian. Then $\localrank(\langle A,B\rangle)=2$.
\EndProposition

\Proof
Let $C$ be an arbitrary finitely generated subgroup of $\langle A,B\rangle$. Then there exist finitely generated subgroups $A_0$ and $B_0$ of $A$ and $B$ respectively such that $C\le\langle A_0,B_0\rangle$. Now since $A\cap B$ is non-abelian, we may fix non-commuting elements $u,v$ of $A\cap B$. The groups $\langle A_0,u,v\rangle$ and $\langle B_0,u,v\rangle$ are finitely generated subgroups of $A$ and $B$ respectively. Since $A$ and $B$ are of local rank $2$, there are subgroups $X$ and $Y$ of $A$ and $B$ respectively, having rank at most $2$,  such that $\langle A_0,u,v\rangle\le X$ and $\langle B_0,u,v\rangle\le Y$. Since each of the groups $X$ and $Y$ contains the non-commuting elements $u$ and $v$, they must be of rank exactly $2$. Set $F=\langle X,Y\rangle$. We have $\rank F\le\rank X+\rank Y\le4$. Since $\Gamma$ is $4$-free it follows that $F$ is free. The rank-$2$ groups $X$ and $Y$ are subgroups of $F$; and since $X\cap Y$ contains the non-commuting elements $u$ and $v$, we have $\rank(X\cap Y)\ge2$. 
We now apply a result due to Kent and Louder-McReynolds, included in \cite[Theorem 2]{kent} and in 
\cite[Theorem 1.1]{LMcR}, 
which asserts that if $X$ and $Y$ are rank-$2$ subgroups of a free group, and if $\rank(X\cap Y)\ge2$, then $\rank\langle X,Y\rangle=2$. By construction we have $C\le\langle X,Y\rangle$, and since $C$ was an arbitrarily finitely generated subgroup of $\langle A,B\rangle$, it follows that $\localrank (\langle A,B\rangle)\le2$. But $\langle A,B\rangle$ contains the non-abelian subgroup $A\cap B$, and hence $\localrank (\langle A,B\rangle)=2$.
\EndProof

\section{Nerves}\label{nerve section}

This section and the next one include a review of some general conventions from \cite{gs}. 
In the context in which they are applied in this paper, some of these conventions are equivalent to conventions used in \cite{fourfree}, but we have adopted a number of conventions from \cite{gs} because they provide a viewpoint that seems more natural and more flexible.

\Number\label{complex stuff}
Except when we specify otherwise, the term {\it simplicial complex} will be understood in the geometric sense; that is, a simplicial complex is a set $L$  of pairwise-disjoint finite-dimensional open simplices in a (possibly infinite-dimensional) real vector space, with the property that any face of any simplex in $L$ is itself in $L$. The geometric realization of an abstract simplicial complex is a simplicial complex in this sense, and every simplicial complex is simplicially isomorphic to the geometric realization of an abstract simplicial complex.  The simplicial complexes referred to in this paper are not assumed to be locally finite. If $L$ is a simplicial complex, the union of its simplices will be denoted by $|L|$. The set $|L|$ will always be understood to be endowed with the weakest topology which induces the standard topology on the closure of each simplex.

We will not distinguish between vertices and $0$-simplices of a simplicial complex $K$; that is, if $v$ is a vertex of $K$, we will denote by $v$ the $0$-simplex which in more formal notation would be denoted by $\{v\}$.

As is customary, if $d$ is a non-negative integer, we define the $d$-skeleton of a simplicial complex $K$ to be the subcomplex consisting of all simplices of $K$ having dimension at most $d$, and denote it by $K^{(d)}$. In particular, $K^{(0)}$ is 
identified with
the set of vertices of $K$.

Let  $L$ be a simplicial complex. We say that a subset $W$  of   $|L|$  is {\it saturated}  if $W$ is a union of (open)
simplices of $L$.

An indexed family $\calf=(U_i)_{i\in I}$ of nonempty (open) subsets of a topological space $X$ is said to {\it cover} $X$ (or to be a(n open) {\it covering} of $X$) if $X=\bigcup_{i\in I}U_i$. (Here the index set $I$ can be any set whatsoever.) We define the {\it abstract nerve} of a covering $\calf=(U_i)_{i\in I}$ of $X$ to be the abstract simplicial complex that is well defined up to canonical simplicial isomorphism as follows. The vertex set $V$ of the complex is a bijective copy of the index set $I$, equipped with a specific bijection $i\mapsto v_i$ from $I$ to $V$. A simplex $\sigma$  is a set  $\{v_{i_0},\ldots,v_{i_d}\}$, with $d\ge0$ and $i_0,\cdots,i_d\in I$, such that
$U_{i_0}\cap\cdots\cap U_{i_d}\ne\emptyset$. 
The  {\it nerve} of a covering is defined to be the geometric realization of its abstract nerve.
The nerve of a covering $\calf$ will be denoted $K_\calf$.

Note that as in \cite{gs}, the definition of a covering and its nerve given here, unlike the most classical definition,  allows the possibility that the covering $\calf$ is ``non-faithfully indexed'' in the sense that  there exist distinct $i,j\in I$ for which $U_i=U_j$. This affects the definition of the nerve of $\calf$, and will be needed for the proof of Proposition \ref{linky borsuk} below, in which the covering $\calg$ may be ``non-faithfully indexed.'' 
The proof of Proposition \ref{linky borsuk} depends on a version of the Borsuk Nerve Theorem, which is proved in \cite{gs} and is paraphrased below as  Prop. \ref{still our borsuk}, and applies to coverings that are not necessarily ``faithfully indexed.''

\EndNumber

\NotationRemarks\label{nervy notation} 
Let $\calf=(U_i)_{i\in I}$ be  an open
covering of a topological space $X$. For every simplex $\sigma$ of $K_\calf$, we will denote by $\cals _\sigma^{\calf } $ (or by $\cals _\sigma$ when the covering $\calf$ is understood) the set of all indices $i\in I$ such that $v_i$ is a vertex of $\sigma$. (In \cite[Definitions 2.9]{gs}, the object denoted here by $\cals_\sigma$ was defined in a special case, and was denoted by $I_\sigma$.) We will denote by
$\calu_\sigma^{\calf } $ (or by $\calu_\sigma$ when the covering $\calf$ is understood) the set $\bigcap_{i\in \cals _\sigma }U_i$. 

Thus the definition of the nerve $K_\calf$ implies that $\calu_\sigma^{\calf }\ne\emptyset$ for every simplex $\sigma$ of $K_\calf$. 

Note that if $\tau$ is a face of a simplex $\sigma$ of $K_\calf$, we have $\calU_\sigma\subset \calU_\tau$.

If an index $i$ belongs to $\cals _\sigma$,
then by definition we have 
$\calu_\sigma\subset U_i$. We denote by
${\mathscr T}_\sigma^{\calf}$ (or by ${\mathscr T}_\sigma$ when the covering $\calf$ is understood)
the set of all indices $i\in I$ such that $i\notin \cals _\sigma$ but $\calu_\sigma\cap U_i\ne\emptyset$. 
\EndNotationRemarks

Here is our version of the Borsuk Nerve Theorem:

\Proposition\label{still our borsuk}
Let $\calf=(U_i)_{i\in I}$ be  an open
covering of a paracompact space $X$. Suppose that for every simplex $\sigma$ of $K_\calf$, the set $\calu_\sigma$ is contractible. Then $|K_\calf|$ is homotopy-equivalent to $X$. 
\EndProposition

\Proof
This is a paraphrase of Proposition 2.7 of \cite{gs}.
\EndProof

Various special cases of the following result were implicit in the proofs of \cite[Lemma 5.7]{fourfree}, \cite[Lemma 3.8]{rosemary}, and \cite[Lemma 3.3]{gs}.

\Proposition\label{linky borsuk}
Let $\calf=(U_i)_{i\in I}$ be  an open
covering of a paracompact space $X$. Suppose that for every 
simplex 
$\phi$
of $K_\calf$, the set 
$\calu_\phi$ 
is contractible.
Let $\sigma$ be a simplex of $|K_\calf|$. Then the set $\calu_\sigma\cap\bigcup_{i\in {\mathscr T}_\sigma} U_i$ is homotopy-equivalent to $\link_{K_\calf}\sigma$. 
\EndProposition

\Proof
Set $\calu=\calu_\sigma^{\calf } $.
Set $\calw=\calu\cap\bigcup_{i\in {\mathscr T}_\sigma} U_i$.

Set { $V_{{i}} = U_{{i}} \cap {\mathcal U} $
 for each $ {i} \in {\mathscr T}_{\sigma}$,}  and $\calg=(V_{{i}})_{{i} \in {\mathscr T}_{\sigma}}$.
We have $\bigcup_{{i} \in {\mathscr T}_{\sigma}}V_{{i}}=\bigcup_{{i} \in {\mathscr T}_{\sigma}} (U_{{i}} \cap {\mathcal U} )=
(\bigcup_{{i} \in {\mathscr T}_{\sigma}} U_{{i}} )\cap {\mathcal U} =\calw$. Hence
$\calg$ is a cover for $\calw$.

According to \ref{complex stuff}, the nerve {$K_{\calg}$} comes equipped with a bijection from the index set ${\mathscr T}_\sigma$ of the covering $\calg$ to the vertex set 
(or $0$-skeleton)
${K_{\calg}}^{(0)}$  of ${K_{\calg}}$.  We will denote this bijection by $i\mapsto w_i$.

Since ${\mathscr T}_\sigma\subset I$, and since $i\mapsto v_i$ is a bijection from the index set  $I$ of $\calf$ to the vertex set $K_\calf^{(0)}$ of $K_\calf$,
 there is a well-defined injection ${f}^{(0)}:{K_{\calg}}^{(0)}\to K_\calf^{(0)}$ given by ${f}^{(0)}(w_i)=v_i$. If $\tau$ is any simplex of $K_\calg$, we have $\bigcap_{i\in\cals_\tau}V_i\ne\emptyset$ by the definition of the nerve $K_\calg$. But $\bigcap_{i\in\cals_\tau} V_i=\bigcap_{i\in\cals_\tau}(U_i\cap\calu)\subset \bigcap_{i\in\cals_\tau}U_i$, and hence $\bigcap_{i\in\cals_\tau}U_i\ne\emptyset$; in view of the definition of 
the nerve $K_\calf$, it follows that $\{v_i:i\in\cals_\tau\}$, the image  under $f^{(0)}$ of the vertex set $\{w_i:i\in\cals_\tau\}$ of $\tau$, is the vertex set of a simplex of $K_\calf$. Since this is the case for every simplex $\tau$  of $K_\calg$, the map $f^{(0)}$ extends to a simplicial map $f:K_\calg\to K_\calf$. Since $f^{(0)}$ is injective, $f$ is also injective.

We now claim:

\Claim\label{not so fast}
For every simplex $\tau$ of $K_\calg$, we have $f(\tau)\subset {\rm link}_{K_\calf}(\sigma)$, and  the set $\calu_\tau^{\calg}$ is contractible.
\EndClaim

To prove \ref{not so fast}, first note that by the definition of the nerve $K_\calg$
we have $\bigcap_{i\in\cals_\tau}V_i\ne\emptyset$.  But $\bigcap_{i\in\cals_\tau} V_i=\bigcap_{i\in\cals_\tau}(U_i\cap\calu)
=(\bigcap_{i\in\cals_\tau}U_i)\cap\calu
=(\bigcap_{i\in\cals_\tau} U_i)\cap(\bigcap_{i\in\cals_\sigma}U_i)=
\bigcap_{i\in\cals_\tau\cup\cals_\sigma} U_i$. Hence $\bigcap_{i\in\cals_\tau\cup\cals_\sigma} U_i\ne\emptyset$, which by the definition of  $K_\calg$ means that $\cals_\tau\cup\cals_\sigma=\cals_\phi$ for some simplex $\phi$ of $K_\calf$.
But according to \ref{nervy notation},
the index set ${\mathscr T}_\sigma\subset I$ of $\calg$ is disjoint from $\cals_\sigma$; in particular, $\cals_\tau\cap\cals_\sigma=\emptyset$. Thus $\cals_\phi$ is the disjoint union of $\cals_\tau$ and $\cals_\sigma$, so that the vertex set $\{v_i:i\in\cals_\phi\}$ of $\phi$
is the disjoint union of the vertex sets $\{v_i:i\in\cals_\tau\}$ and $\{v_i:i\in\cals_\sigma\}$ of $f(\tau)$ and $\sigma$. This shows that $f(\tau)\subset\link_{K_\calf}\sigma$.

On the other hand, the equality
$\bigcap_{i\in\cals_\tau} V_i=
\bigcap_{i\in\cals_\tau\cup\cals_\sigma} U_i$ may be written as $\calu_\tau^{\calg}=
\bigcap_{i\in\cals_\phi} U_i=\calu_\phi^{\calf } $; since the hypothesis of the proposition implies that 
$\calu_\phi ^\calf$
is contractible, we have shown that $\calu_\tau^{\calg}$ is contractible. Thus (\ref{not so fast}) is proved.

Next, we claim:
\Equation\label{it's the link}
{f}({K_{\calg}})={\rm link}_{K_\calf}(\sigma).
\EndEquation

To prove (\ref{it's the link}), first note that the inclusion ${f}({K_{\calg}})\subset{\rm link}_{K_\calf}(\sigma)$ is immediate from the first assertion of \ref{not so fast}. To establish the reverse inclusion, consider an arbitrary  simplex $\psi$ of ${\rm link}_{K_\calf}(\sigma)$. There is a simplex $\phi$ of $K_\calf$ whose vertex set is the disjoint union of the vertex sets of $\sigma$ and $\tau$. Hence  $\cals_\sigma\cap\cals_\psi=\emptyset$ and $\cals_\sigma\cup\cals_\psi=\cals_\phi$. By the definition of $K_\calf$ we have $\emptyset\ne\bigcap_{i\in\cals_\phi}U_i=
(\bigcap_{i\in\cals_\sigma}U_i)\cap(\bigcap_{i\in\cals_\psi}U_i)=\calu\cap(\bigcap_{i\in\cals_\psi}U_i)$, i.e.
\Equation\label{estate sale}
\bigcap_{i\in\cals_\psi}(U_i\cap\calu)\ne\emptyset.
\EndEquation

In particular, it follows from (\ref{estate sale}) that if $i$ is any index in $\cals_\psi$, we have $U_i\cap\calu\ne\emptyset$; on the other hand, since $\cals_\sigma\cap\cals_\psi=\emptyset$, we have $i\notin\cals_\sigma$. By definition we therefore have $i\in {\mathscr T}_\sigma$. This shows that $\cals_\psi\subset {\mathscr T}_\sigma$. We may now rewrite (\ref{estate sale}) in the form $\bigcap_{i\in\cals_\psi}V_i\ne\emptyset$.
By the definition of $K_\calg$ it follows that $\cals_\psi=\cals_\tau$ for some simplex $\tau$ of $K_\calg$. Hence $f$ maps the vertex set $\{w_i:i\in\cals_\tau\}$ of $\tau$ onto
the vertex set $\{v_i:i\in\cals_\psi\}$ of $\psi$, so that $f(\tau)=\psi$. This establishes the inclusion ${\rm link}_{K_\calf}(\sigma)\subset{f}({K_{\calg}})$, and completes the proof of (\ref{it's the link}).

It follows from (\ref{it's the link}) that the injective simplicial map $f$ may be regarded as a simplicial isomorphism between $K_{\calg}$ and ${\rm link}_{K_\calf}(\sigma)$. But by  the second assertion of \ref{not so fast}, the covering $\calg$ of $\calw$ has the property that for every simplex $\tau$ of $K_\calg$ the set $\calu_\tau^{\calg}$ is contractible. It therefore follows from Proposition \ref{still our borsuk} that $K_\calg$ is homotopy-equivalent to $\calw$. Hence ${\rm link}_{K_\calf}(\sigma)$ is also homotopy-equivalent to $\calw$, as asserted by the present proposition.
\EndProof
%\sigma

\Proposition\label{linky corollary}
Let $\calf=(U_i)_{i\in I}$ be  an open
covering of a paracompact space $X$. Set $K=K_\calf$. Suppose that 
\begin{enumerate}
\item for every point $P$ of $X$ there exist two distinct indices $i,j\in I$ such that $P\in U_i\cap U_j$, and that
\item for
 every 
simplex $\phi$ of $K$, the set $\calu_\phi$ is contractible. 
\end{enumerate}
Then every vertex of $K$ has a contractible link in $K$. 

If in addition we assume that
\begin{itemize}
\item[{\it(3)}]
the complex $K$ is finite-dimensional, 
\end{itemize}
then $|K|-|K^{(0)}|$ is homotopy-equivalent to $X$
(where, as in \ref{complex stuff}, $K^{(0)}$ denotes the $0$-skeleton of $K$).
\EndProposition

\Proof
Let $w$ be an arbitrary vertex of $K$. We may write $w=v_{i_0}$ for some index $i_0\in I$. If we regard $w$ as a $0$-simplex 
(cf. \ref{complex stuff}),
then according to the  definitions (see \ref{nervy notation})
we have $\calu_w=U_{i_0}$ and ${\mathscr T}_w=I-\{i_0\}$.
According to Condition (1) of the hypothesis, for every point $p\in U_{i_0}$ there is an index $i\in I$ such that $i\ne i_0$ and $p\in U_i$. We therefore have $U_{i_0}\subset\bigcup_{i\in I-\{i_0\}} U_i$,
so that
$\calu_w\cap\bigcup_{i\in {\mathscr T}_w} U_i=U_{i_0}\cap\bigcup_{i\in I-\{i_0\}} U_i=U_{i_0}$. It therefore follows from Proposition \ref{linky borsuk}, and Condition (2) of the hypothesis of the present proposition, that 
$\link_{K}w$ is homotopy-equivalent to $U_{i_0}$. But a second application of Condition (2) of the hypothesis gives that $U_{i_0}$ is contractible. This establishes the first assertion of the present proposition. 

According to \cite[Proposition 3.2]{gs}, if $K'$ is a subcomplex of a finite-dimensional simplicial complex $K$ such that $\link_K\sigma$ is contractible for every simplex $\sigma\in K'$, then the inclusion map of the saturated set $|K|-|K'|$ into $|K|$ is a homotopy equivalence. In the present situation, we may apply this by setting $K'=K^{(0)}$. The first assertion of the present proposition gives that  every vertex of $K^{(0)}$ has a contractible link in $K$. It follows that $|K|-|K^{(0)}|$ is homotopy equivalent to $|K|$. But Proposition \ref{still our borsuk}, together with Condition (2) of the hypothesis of the present proposition, implies that $|K|$ is in turn homotopy equivalent to $X$.
\EndProof

\section{Quantitative geometry of hyperbolic manifolds}\label{quant section}

\NotationRemarks\label{nbhd}
Throughout the paper we will use $\dist$ as a default notation for the distance function in a metric space, when it is clear to which distance function we are referring. 
If $A$ and $B$ are non-empty subsets  of a metric space, $\dist(A,B)$ will denote the quantity $\inf_{x\in A,y\in B}\dist(x,y)$. We write $\dist(p,A)=\dist(\{p\},A)$ when $p$ is a point of a metric space and $A$ is a non-empty subset of the space.

If $p$ is a point of a metric space, and $r$ is a real number, we will denote by $\nbhd_r(p)$ the set of all $x\in X$ such that $\dist(x,p)<r$. Thus $\nbhd_r(p)$  is a neighborhood of $p$ if $r>0$, and is empty if $r\le0$.

We will often use the following consequence of the triangle inequality: if $p$ and $p'$  are points of a metric space, and $r$ and $r'$  are real numbers with $r+r'>\dist(p,p')$, then $\nbhd_r(p)\cap\nbhd_{r'}(p')=\emptyset$.

If $X$ is a subset of a hyperbolic manifold $M$, one can make $X$ into a metric space by defining the distance between two points of $X$ to be their  distance in $M$; this is the {\it extrinsic distance function}. 
If $X$ is compact and non-empty, the extrinsic distance function gives rise to the  {\it extrinsic diameter} of $U$, which is the maximum of the set of all 
distances in $M$ between points of $X$. On the other hand, if $X$ is connected and open, the {\it intrinsic distance} between points $p$ and $p'$ of $X$ is the infimum of all lengths of paths in $U$ between $p$ and $p'$. This gives rise to  the notion of an {\it instrinsic isometry}
between open connected subsets $X$ and $X'$ of hyperbolic manifolds $M$ and $M'$: it is a diffeomorphism 
%re {\it intrinsically isometric} if there is a diffeomorphism 
between $X$ and $X'$ that preserves lengths of paths.

A closed geodesic in a hyperbolic $3$-manifold $M$ may be regarded as  a map $c:S^1\to M$. The set $c(S^1)\subset M$ will be called the {\it support} of the geodesic $c$ and will denoted $|c|$.

If $M$ is a closed, orientable hyperbolic $3$-manifold, we will denote
by $\ell_M$ the length of the shortest closed geodesic in $M$.

A {\it hyperbolic ball} in a hyperbolic $3$-manifold $M$ is a subset of $M$ which is intrinsically isometric to a ball  $N\subset\HH^3$; its {\it radius} is the radius of $N$, and its {\it center} is the unique point that is mapped to the center of $N$ by
an intrinsic
isometry.

By a {\it hyperbolic cylinder} in $\HH^3$ we mean a set of the form $Z=\{P\in\HH^3:\dist(P,l)<R\}$, where $l$ is a line in $\HH^3$ and $R$ a positive number; the line $l$ and the number $R$, which are uniquely determined by the set $Z$, will be respectively called the {\it axis} and {\it radius} of the cylinder.
Note that any hyperbolic cylinder is a convex subset of $\HH^3$.

A {\it tube} in a hyperbolic $3$-manifold $M$ is a set $T\subset M$ which is intrinsically isometric to a quotient $Z/\langle\tau\rangle$, where 
  $Z\subset\HH^3$ is a hyperbolic cylinder  and $\tau$ is a loxodromic  transformation
whose axis is the axis
$l$ of $Z$. The {\it radius} of $T$ 
is the radius of $Z$, and its {\it core} is the unique simple closed geodesic  
$c$ such that $|c|$
is mapped to $l/\langle\tau\rangle$ by 
an intrinsic
isometry from $T$ to $Z/\langle\tau\rangle$.

We shall denote by $\isomplus(\HH^3)$ the group of all orientation-preserving isometries of $\HH^3$. For each $x\in\isomplus(\HH^3)$ and for each point  $P\in\HH^3$,
we shall write $d(x,P)=\dist(P,x\cdot P)$.

If $c$ is a closed geodesic in an orientable hyperbolic $3$-manifold $M$,  the {\it tube radius}  of $c$ is defined to be $0$ if $c$ is not simple; and if $c$ is simple it is defined to be the radius of the largest tube having core  $c$. 

If we write a given orientable hyperbolic $3$-manifold $M$ as $\HH^3/\Gamma$, where $\Gamma$ is a discrete, torsion-free subgroup of $\isomplus(\HH^3)$, then for any closed geodesic $c$ in $M$, the pre-image of $|c|$ under the quotient projection $\HH^3\to M$ is a line $A$ in $\HH^3$ whose stabilizer is a cyclic subgroup $C$ of $\Gamma$, and the tube radius $r$ of $c$  then satisfies 
$$2r=\max_{\alpha\in\Gamma-C}\dist(A,\alpha\cdot A).$$

\EndNotationRemarks

\Number\label{i see cg}
As in \cite{gs}, we will say that a group  has 
the {\it infinite cyclic centralizer property}, or  is an {\it \iccg}, if the centralizer of every non-trivial element of $\Gamma$ is infinite cyclic. If $\Gamma$ is an \iccg\ then $\Gamma$ is torsion-free, and
%be a group with the property that the centralizer of every non-trivial element is infinite cyclic. Thus
 every non-trivial element $x$ of $\Gamma$ belongs to a unique maximal cyclic subgroup of
$\Gamma$, namely the centralizer of $x$ in $\Gamma$. Hence two non-trivial elements of $\Gamma$ commute if and only if they lie in the same maximal cyclic subgroup.

Let $\Gamma$ be a discrete, cocompact subgroup of ${\rm Isom}_+(\bf H^3)$. 
Then $\Gamma$ is an \iccg.
Given a real number $\lambda >0$, we define ${\mathcal C}_{\lambda}(\Gamma) 
%\subseteq {\mathcal C}(\Gamma)
$ to be the set of maximal cyclic subgroups $ C$ of $\Gamma$ 
such that a (loxodromic) generator of $ C$ has translation length less than $\lambda$.

As in \cite{gs}, for every non-trivial element $x$ of  $\Gamma$, we set $Z_{\lambda}(x):=\{ P \in {\bf H}^3 :
d(x,P) < \lambda \}$, and for each  $C\in {\mathcal C}_\lambda(\Gamma)$, 
we set $Z_{\lambda}(C)= \bigcup_{1 \not =  x \in {\mathcal C}} Z_{\lambda}( x)$. 
For each ${ C} \in {\mathcal C}_{\lambda}(\Gamma)$, we have $Z_{\lambda} ({ C})=
Z_{\lambda} ( x)$ for some
 element $ x\ne1$ of $C$; furthermore,  $Z_{\lambda} ({ C})$ is a 
hyperbolic cylinder 
(see \ref{nbhd})
whose axis is the common translation axis of all non-trivial elements of $C$.
 %(If $C \in C(\Gamma)-C_{\lambda}(\Gamma)$, then $Z_{\lambda}(C)=\emptyset$.) 
For any $\lambda >0$ we set  $\mathcal{Z}_{\lambda}(\Gamma)= (Z_{\lambda}(C))_{C\in \calc_\lambda(\Gamma)}$.

In this paper we will denote by $\frakX_\lambda(\Gamma)$ the set of all points $P\in\HH^3$ such that $d(x,P)<\lambda$ for some non-trivial element $x$ of $\Gamma$.
It follows from the definitions that $\frakX_\lambda(\Gamma)=\bigcup_{C\in \calc_\lambda(\Gamma)} Z_{\lambda}(C)$.
%so that $ \mathcal{Z}_{\lambda}(\Gamma)$. is a covering  of $\frakX_\lambda(\Gamma)$, in the sense discussed in \ref{complex stuff}. It should be borne in mind that the index set for this covering is $\calc_\lambda(\Gamma)$, a collection of maximal cyclic groups of $\Gamma$.
%The nerve $K_{{\mathcal Z}_{\lambda}(\Gamma)}$ is well defined by \ref{complex stuff}.''

%``This would mean replacing all occurrences of 
%$\frakX_\lambda(\Gamma)$ by $\HH^3$, and always accompanying them by the assumption that 
%there  is a point $p\in M:=\HH^3/\Gamma$ with $\shortone(p)\ge\lambda$. There won't be all that many such occurrences.''

\EndNumber

\DefinitionNotationRemarks\label{short and next}
If $ p $ is 
a point of a closed, orientable hyperbolic $3$-manifold $M$, we denote by $\shortone( p )>0$ the minimum length of a homotopically non-trivial loop based at $ p $. 
In the notation of \ref{nbhd}, we have 
$\ell_M=\min_{p\in M}\shortone(p)$.

We shall say that a point $p\in M$ 
is {\it $\alpha$-thin} for a given $\alpha>0$ if $\shortone(p)<\alpha$, and that $p$ is {\it $\alpha$-thick}  if $\shortone(p)\ge\alpha$. The point $p$ is $\alpha$-thick if and only if it is the center of a hyperbolic ball of radius $\alpha/2$ in $M$.

We denote by $\Mthin(\alpha)$ the set of all $\alpha$-thin points of $M$; thus  $\Mthick(\alpha):=M-\Mthin(\alpha)$ is the set of all $\alpha$-thick points of $M$. Note that if $M$ is written as $\HH^3/\Gamma$ where $\Gamma\le\Isom_+(\HH^3)$ is discrete, torsion-free and cocompact, 
and if $q:\HH^3\to M$ denotes the quotient map, then in the notation of 
\ref{i see cg} 
we have $q^{-1}(\Mthin(\alpha))=\frakX_\alpha(\Gamma)$ for every $\alpha>0$. In particular we have $\frakX_\alpha(\Gamma)=\HH^3$ if and only if $\Mthin(\alpha)=M$, i.e. if and only if $\Mthick(\alpha)
=
\emptyset$.

Given a point $p$ of the closed, orientable hyperbolic $3$-manifold $M$, 
we 
will say that a maximal cyclic subgroup $C$ of   $\pi_1(M,p)$ is {\it short} if some non-trivial element of $C$ is represented
% (or some such term), meaning one that contains a non-trivial element represented 
by a loop of length $\shortone( p )$ based at $p$. It follows from  \ref{i see cg} that $\pi_1(M)$ is an \iccg. Thus every element of $\pi_1(M,p)$ lies in a unique maximal cyclic subgroup, and hence  it follows from the definitions that $\pi_1(M,p)$ has at least one short maximal cyclic subgroup. 

We also define a positive real number $\nextone( p )$ as follows: if $\pi_1(M,p)$ has only one short maximal cyclic subgroup, denoted
%all elements of the \iccg\ $\pi_1(M, p )$ represented by loops of length $\shortone( p )$ lie in the same maximal cyclic subgroup 
$C_0$, we 
define $\nextone( p )$ to be the minimum length of a loop defining an element of $\pi_1(M,p)-C_0$. If there are two or more short maximal cyclic subgroups, we set $\nextone( p )=\shortone( p )$.

Thus in all cases, if $C$ is any short maximal cyclic subgroup,
% containing an element that is  represented by a loop of length $\shortone( p )$, 
then $\nextone( p )$ is the minimum length of a loop defining an element of $\pi_1(M,p)-C$.

If $M$ is given as a quotient $\HH^3/\Gamma$, where $\Gamma\le\isomplus(\HH^3)$ is discrete, cocompact and torsion-free, and if $P\in\HH^3$ is a point that projects to $p\in M$ under the quotient map,  we have $\shortone( p )=\min_{1\ne x\in\Gamma}d(x,P)$. Furthermore, if $z$ is a non-trivial element of $\Gamma$ with $d(z,P)=\shortone( p )$, and if $C$ denotes the maximal cyclic group of the \iccg\ $\Gamma$ containing $z$, then $\nextone( p ) =\min_{x\in\Gamma-C}d(x,P)$.

Notice also that we have $\nextone( p )\ge\shortone( p )$ for every $ p \in M$.
\EndDefinitionNotationRemarks

\Number\label{what about the workers}
As in \cite{fourfree} (see Definition 1.3 of that paper), we define a point $p$ of a closed, orientable hyperbolic $3$-manifold $M$ to be {\it $\lambda$-semithick}, 
where $\lambda$ is a given positive number, if any two loops of length less than $\lambda$ based at $p$ define elements of $\pi_1(M,p)$ which commute. Since $\pi_1(M,p)$ is an \iccg, this is equivalent to saying that any two homotopically non-trivial loops of length less than $\lambda$ based at $p$ define elements in the same maximal cyclic subgroup of $\pi_1(M,p)$. 

It is an immediate consequence of this definition, together with the definition of $\nextone(p)$ given in \ref{short and next}, that $p$ is $\lambda$-semithick if and only if $\nextone(p)\ge\lambda$.

In \cite[3.8]{fourfree}, a set $\frakG_M$ is defined to be the set of all points $p\in M$ for which the elements of $\pi_1(M,p)$ generated by all elements represented by loops of length $\shortone(p)$ is a cyclic group. (The quantity that we denote by $\shortone(p)$ in this paper is denoted by $\ell_p$ in \cite{fourfree}.) In terms of the definition given in \ref{short and next}, this means that $p\in\frakG_M$ if and only if  $\pi_1(M,p)$ has only one short maximal cyclic subgroup. Furthermore, in \cite[3.8]{fourfree}, when $p\in\frakG_M$, the unique maximal cyclic subgroup of $\pi_1(M,p)$ containing at least one non-trivial element represented by   a loop of length less than $\shortone(p)$ is denoted by $C_p$, and $\fraks_M(p)$ denotes the smallest length of a loop based at $p$ that does not represent an element of $C_p$; thus from the point of view of the present paper, $C_p$ is the unique short maximal cyclic subgroup of $\pi_1(M,p)$ when $p\in\frakG_M$, and $\fraks_M$ is simply the restriction of $\nextone$ to $\frakG_M$. Note also that, according to our definitions, we have $\nextone(p)=\shortone(p)$ for any $p\in M-\frakG_M$.

These comparisons of the conventions of the present paper to those of \cite{fourfree} will be useful in the proofs of Proposition \ref{soft caps} and Lemmas \ref{surgeon general} and \ref{surgeon special} below.

\EndNumber

\Lemma\label{continuous}
Let $M$ be a closed, orientable $3$-manifold. 
\begin{enumerate}[(a)] 
\item For every point $p\in M$ we have $\shortone(p)+\nextone(p)=\min_{(\alpha,\beta)}(\length(\alpha)+\length(\beta))$, where $(\alpha,\beta)$ ranges over all pairs of loops based at $p$ such that $[\alpha]$ and $[\beta]$ do not commute in $\pi_1(M,p)$.
\item For any two points $p,p'\in M$ we have $|\shortone(p)-\shortone(p')|\le2\dist(p,p')$ and $|(\shortone(p)+\nextone(p))-(\shortone(p')+\nextone(p'))|\le4\dist(p,p')$.
\item The functions $\shortone$ and $\nextone$ are continuous on $M$.
\end{enumerate}
\EndLemma

\Proof
First note that by the definition of $\shortone(p)$, we may fix a loop $\alpha_0$ based at $p$ such that $1\ne[\alpha_0]\in\pi_1(M,p)$ and $\length(\alpha_0)=\shortone(p)$; and according to the discussion in \ref{short and next}, if $C$ denotes the maximal cyclic subgroup of  the \iccg\ $\pi_1(M,p)$   containing $[\alpha_0]$, we may fix a  loop $\beta_0$ based at $p$ such that $[\beta_0]\notin C$ and $\length(\beta_0)=\nextone(p)$. Since $[\beta_0]\notin C$, the elements $[\alpha_0]$ and $[\beta_0]$ of $\pi_1(M,p)$ do not commute. We have $\length(\alpha_0)+\length(\beta_0)=\shortone(p)+\nextone(p)$. To prove (a), it remains only to prove that if $\alpha$ and $\beta$ are arbitrary loops based at $p$ such that $[\alpha]$ and $[\beta]$ do not commute  in $\pi_1(M,p)$, we have $\length(\alpha)+\length(\beta)\ge\shortone(p)+\nextone(p)$. Since $[\alpha]$ and $[\beta]$ do not commute, they cannot both lie in $C$. Hence by symmetry we may assume that $[\beta]\notin C$. The discussion in \ref{short and next} then shows that $\length(\beta)\ge\nextone(p)$. Furthermore, since $[\alpha]$ is in particular non-trivial, the definition of $\shortone(p)$ gives $\length(\alpha)\ge\shortone(p)$. Hence  $\length(\alpha)+\length(\beta)\ge\shortone(p)+\nextone(p)$ as required, and (a) is proved.

To prove (b), set $d=\dist(p,p')$, and let $\gamma$ denote a geodesic path from $p$ to $p'$ with $\length(\gamma)=d$. If $\alpha$ is a homotopically non-trivial loop based at $p$, then $\alpha':=\overline{\gamma}\star\alpha\star\gamma$ is a homotopically non-trivial loop based at $p'$, and $\length(\alpha')=2d+\length(\alpha)$. Since $\shortone(p)$ and $\shortone(p')$ are by definition the minimum lengths of homotopically non-trivial loops based at $p$ and $p'$ respectively, it follows that
\Equation\label{tinker}
\shortone(p')\le\shortone(p)+2d.
\EndEquation
Likewise, if $\alpha$ and $\beta$ are loops based at $p$ such that $[\alpha]$ and $[\beta]$ do not commute, then $\alpha':=\overline{\gamma}\star\alpha\star\gamma$ and $\beta':=\overline{\gamma}\star\beta\star\gamma$ are loops based at $p'$ such that $[\alpha']$ and $[\beta']$ do not commute, and $\length(\alpha')+\length(\beta')=4d+\length(\alpha)+\length(\beta)$. Since, by (a), we have
$\shortone(p)+\nextone(p)=\min_{(\alpha,\beta)}(\length(\alpha)+\length(\beta))$, where $(\alpha,\beta)$ ranges over all pairs of loops based at $p$ such that $[\alpha]$ and $[\beta]$ do not commute, and
$\shortone(p')+\nextone(p')=\min_{(\alpha',\beta')}(\length(\alpha')+\length(\beta'))$, where $(\alpha',\beta')$ ranges over all pairs of loops based at $p'$ such that $[\alpha']$ and $[\beta']$ do not commute, it follows that
\Equation\label{evers}
\shortone(p')+\nextone(p')\le\shortone(p)+\nextone(p)+4d.
\EndEquation
Interchanging the roles of $p$ and $p'$ in (\ref{tinker}) and in (\ref{evers}), we obtain
\Equation\label{stinker}
\shortone(p)\le\shortone(p')+2d
\EndEquation
and
\Equation\label{severs}
\shortone(p)+\nextone(p)\le\shortone(p')+\nextone(p')+4d.
\EndEquation
Now (\ref{tinker}) and (\ref{stinker}) imply the inequality  $|\shortone(p)-\shortone(p')|\le2d$, while (\ref{evers}) and (\ref{severs}) imply $|(\shortone(p)+\nextone(p))-(\shortone(p')+\nextone(p'))|\le4d$. Thus (b) is proved.

To prove (c), note that it follows from (b) that the functions $\shortone$ and $\shortone+\nextone$ are continuous on $M$; hence $\nextone$ is also continuous.
\EndProof

\RemarkNotation\label{big and little next one}
Let $M$ be a closed, orientable hyperbolic $3$-manifold. Since by Lemma \ref{continuous}  the function $\nextone$ is continuous on the compact space $M$, this function takes a greatest and a least value on $M$, which we will denote by $\lambda_M$ and $\mu_M$ respectively.
\EndRemarkNotation

\Definition\label{margdef}
Recall that a {\it Margulis number} for an orientable hyperbolic $3$-manifold $M$ is defined to be a positive number $\mu$ such that, for every point $p\in M$ and for any two loops $\alpha$ and $\beta$ of based at $p$ and having length less than $\mu$, the elements $[\alpha]$ and $[\beta]$ of $\pi_1(M,p)$ commute.
\EndDefinition

\Proposition\label{Margulis and s-one and s-two}
Let $M$ be a closed, orientable hyperbolic $3$-manifold. Then the interval $[\mu_M,\lambda_M]$ is the range of the function $\nextone$ on $M$, and the interval $(0,\mu_M]$ is the set of all Margulis numbers for $M$.
\EndProposition

\Proof
Since, by Lemma \ref{continuous}, $\nextone$ is a continuous function on the compact connected space $M$, its range is a closed interval. It follows from the definitions of $\mu_M$ and $\lambda_M$ that they are respectively the left-hand and right-hand endpoints of this interval. This proves the first assertion.

To prove the second assertion, first consider an arbitrary number $\mu$ with $0<\mu\le\mu_M$. Let $p$ be an arbitrary point of $M$, and let $\alpha$ and $\beta$ be loops based at $p$ and having length less than $\mu$. Then their lengths are less than $\mu_M$, which is in turn at most $\nextone(p)$ by definition. According to \ref{short and next}, we may choose a  short maximal cyclic subgroup $C$ of   $\pi_1(M,p)$, and $\nextone(p)$ is the minimal length of any loop representing an element of $\pi_1(M,p)-C$. Since $\alpha$ and $\beta$ have length less than $\nextone(p)$, they represent elements of $C$, and therefore commute. This shows that $\mu$ is a Margulis number for $M$.

Now consider an arbitrary number $\nu>\mu_M$. By the definition of $\mu_M$ we may choose a point $p\in M$ with $\nextone(p)=\mu_M$. According to \ref{short and next}, we may choose a  short maximal cyclic subgroup $C$ of   $\pi_1(M,p)$, and there are loops $\alpha$ and $\beta$ based at $p$, having respective lengths $\shortone(p)$ and $\nextone(p)$, and respectively representing elements of $C$ and 
%is the minimal length of any loop representing an element of 
$\pi_1(M,p)-C$. Since $C$ is a maximal cyclic subgroup of the \iccg\ $\pi_1(M,p)$ (see \ref{i see cg}), the elements $\alpha$ and $\beta$ do not commute. But their lengths are both bounded above by $\nextone(p)=\mu_M<\nu$. This shows that $\nu$ is not a Margulis number for $M$, and the proof of the second assertion is complete.
\EndProof

\Proposition\label{my summer vocation}
Let $M$ be a closed, orientable hyperbolic $3$-manifold, and suppose that $\mu$ is a Margulis number for $M$. Then every component of $\Mthin(\mu)$ is an open solid torus. Hence $\Mthick(\mu)$ is connected (and in particular non-empty).
\EndProposition

\Proof
We use the notation of \ref{i see cg}.
Write $M=\HH^3/\Gamma$, where $\Gamma\le\isomplus(\HH^3)$ is discrete, torsion-free and cocompact.
If $C_1$ and $C_2$ are distinct elements of
$\calc_\mu(\Gamma)$ and $P$ is a point of $Z_{\mu}(C_1)\cap Z_{\mu}(C_2)$, then for $i=1,2$ the definition of $Z_{\mu}(C_i)$ gives a non-trivial element $x_i$ of $C_i$ such that $d(x_i,P)<\mu$. Since the $C_i$ are distinct maximal cyclic subgroups of the \iccg\ $\Gamma$, the elements $x_1$ and $x_2$ do not commute. It follows that if  $p\in M$ denotes the image of $P$ under the quotient map, there are loops of length less than $\mu$ based at $p$ representing non-commuting elements of $\pi_1(M,p)$; this contradicts the definition of a Margulis number. Hence the family $(Z_{\mu}(C))_{C\in \calc_\mu(\Gamma)} $  is pairwise disjoint, and the sets in this family are therefore the components of $\frakX_\mu(\Gamma)=\bigcup_{C\in \calc_\mu(\Gamma)} Z_{\mu}(C)$. It follows that each component of $\Mthin(\mu)$ is the quotient of $Z_{\mu}(C)$, for some $C\in \calc_\mu(\Gamma)$, by its stabilizer, which is $C$. The assertion follows.
\EndProof

\section{
Displacement nerves and \restraining\ vertices
}\label{replaces height section}

\NotationRemarks\label{Q-def}

We define a function $Q$ on the non-negative real numbers by $Q(u)=1/(1+e^u)$. 

Note that $Q$ is strictly monotone decreasing and has range $(0,1/2]$. The inverse function $Q^{-1}$, whose domain is $(0,1/2]$, is given by $Q^{-1}(x)=\log(1/x-1)$.

Let $n\ge1$ be an integer. For any $x>\log(2n-1)$ we have $0<1/2- Q(x)<1/2$. We may therefore define a function $\f_n$ on $(\log(2n-1),\infty)$ by $\f_n(x)=Q^{-1}(1/2-nQ(x))$.
Since $Q$ is strictly monotone decreasing on its domain, the function $\f_n$ is also strictly monotone decreasing. Note also that $\log(2n+1)$ is the unique fixed point of $\f_n$.

Note that the function $\f_1$ has domain $(0,\infty)$, and that for any $x,y>0$ we have
 $y>\fthree(x)$ if and only if $Q(y)+Q(x)<1/2$. 
Similarly, we have
 $y\ge\fthree(x)$ if and only if $Q(y)+Q(x)\le1/2$.

If $P$ is a point of $\HH^3$ and $x$ is an isometry of $\HH^3$, we set $\calq(x,P)=Q(d(x,P))$
where $d(x,P)$ is defined as in \ref{nbhd}).

\EndNotationRemarks

The following result is fundamental for the arguments in this paper. 

\Theorem\label{our log(2k-1)} 
Suppose that $p$ is a positive  integer and that $x_1,\ldots,x_p$ are independent elements of $\isomplus(\HH^3)$ which generate a discrete group. Then for every point $P\in\HH^3$ we have 
$$\calq(x_1,P)+\cdots+\calq(x_p,P)\le\frac12.$$
\EndTheorem

\Proof
For $p=1$ the assertion is trivial. For $p\ge2$ it is a paraphrase of \cite[Theorem 4.1]{acs-surgery}. (The latter result is an extension of \cite[Theorem 6.1]{accs} obtained by combining the results of \cite{accs} with those of \cite{agol}, \cite{cg}, \cite{nz}, and \cite{ohshika}.) 
\EndProof

\Corollary\label{nebulizer}
Let $M$ be a closed, orientable hyperbolic $3$-manifold such that $\pi_1(M)$ is $2$-free. Then for every point $p\in M$ we have $\nextone(p)\ge\f_1(\shortone(p))$. 
\EndCorollary

\Proof
Let us write $M$  as a quotient $\HH^3/\Gamma$, where $\Gamma\le\isomplus(\HH^3)$ is discrete, cocompact and torsion-free, and choose a point $P\in\HH^3$ which projects to $p\in M$ under the quotient map.  According to \ref{short and next}, there is a non-trivial element $z$ of $\Gamma$ with $d(z,P)=\shortone( p )$, and if $C$ denotes the maximal cyclic group of the \iccg\ $\Gamma$ containing $z$, then $\nextone( p ) =\min_{x\in\Gamma-C}d(x,P)$. Let us fix an element $y\in\Gamma-C$ such that $\nextone( p )=d(y,P)$. Since $\Gamma$ is an \iccg, $C$ is the centralizer of $z$ in $\Gamma$ (see \ref{i see cg}), and hence $y$ and $z$ do not commute. But $F:=\langle y,z\rangle$ is free since $\Gamma\cong\pi_1(M)$ is $2$-free. Since $y$ and $z$ do not commute we must have $\rank F=2$, and since the two elements $y$ and $z$ generate $F$, they must form a basis for $F$ (see  \cite[vol. 2, p. 59]{kurosh}). Thus  $ y$ and $z$  are independent. It then follows from Theorem \ref{our log(2k-1)} that $Q(\shortone(p))+Q(\nextone(p))=Q(d(z,P))+Q(d(y,P))=\calq(z,P)+\calq(y,P)\le1/2$, which by \ref{Q-def} implies that
$\nextone(p)\ge\f_1(\shortone(p))$. 
\EndProof

We record here a result, Proposition \ref{big radius corollary from cusp}, which was proved in \cite{cusp} and was deduced from Theorem \ref{our log(2k-1)}. It will be needed at a couple of points in the present paper.

\Proposition\label{big radius corollary from cusp}
Let $k$ and $m$ be integers with $k\ge2$ and $0\le m\le k$. Suppose that
$M$ is a closed, orientable hyperbolic $3$-manifold such
that $\pi_1(M)$ is $k$-free, and let $\mu$ be a Margulis
number for $M$.  Let $p$ be a point of $M$, let $u_1,\ldots,u_m$
be independent elements of $\pi_1(M,p)$ represented by loops
$\alpha_1,\ldots,\alpha_m$ based at $p$, and let $d_j$ denote the
length of $\alpha_j$.  Then there is a $\mu$-thick point $p'\in M$ such that
 $\rho:=\dist (p,p')$ satisfies
$$(k-m)Q(2\rho)+
\sum_{j=1}^mQ(d_j)
\le\frac12.$$
\EndProposition

\Proof
This is a paraphrase of \cite[Corollary 6.2]{cusp}.
\EndProof
%P

\Reformulation\label{reformulated big radius corollary from cusp}
For applications of 
Proposition
\ref{big radius corollary from cusp}, it will be convenient to define, for every  positive integer $n$, a function $\xi_n$ on the interval $(0,1/2)$, by setting $$\xi_n(u)=\frac12Q^{-1}\bigg(\frac 1n\bigg(\frac12-u\bigg)\bigg)$$
whenever $0<u<1/2$. Since $Q$ is strictly monotone decreasing and has range $(0,1/2)$ (see \ref{Q-def}), the function $\xi_n$ is well defined,
and is strictly monotone increasing.
In terms of these definitions, we may reformulate Proposition \ref{big radius corollary from cusp} 
(at least in the case $k>m$)
as follows: if $k$ and $m$ are integers with $k\ge2$ and 
$0\le m< k$,
if
$M$ is a closed, orientable hyperbolic $3$-manifold such
that $\pi_1(M)$ is $k$-free, if $\mu$ is a Margulis
number for $M$, if $p$ is a point of $M$, if $u_1,\ldots,u_m$ are independent elements of $\pi_1(M,p)$ represented by loops
$\alpha_1,\ldots,\alpha_m$ based at $p$, and if $d_j$ denotes the
length of $\alpha_j$,  then 
$\sum_{j=1}^mQ(d_j)<1/2$ (so that
$\xi_{k-m}(\sum_{j=1}^mQ(d_j))$ is defined), and
there is a $\mu$-thick point $p'\in M$ such that
 $\dist (p,p')\ge\xi_{k-m}(\sum_{j=1}^mQ(d_j))$.
\EndReformulation

\Lemma\label{when m is two}
Let $k>2$ be an integer, let $M$ be a closed, orientable hyperbolic $3$-manifold such
that $\pi_1(M)$ is $k$-free, and let $\mu$ be a Margulis
number for $M$. Then for every point $p\in M$ we have $Q(\shortone(p)+Q(\nextone(p))<1/2$ (so that
$\xi_{k-2}(Q(\shortone(p)+Q(\nextone(p)))$ is defined), and there is a $\mu$-thick point $p'\in M$ such that 
 $\dist (p,p')\ge\xi_{k-2}(Q(\shortone(p)+Q(\nextone(p)))$.
\EndLemma

\Proof
By \ref{short and next} we may choose a short maximal cyclic subgroup $C$ of $\pi_1(M,p)$, and there are elements $u_1\in C-\{1\}$ and $u_2\in\pi_1(M,p)-C$ that are represented by loops of respective lengths $\shortone(p)$ and $\nextone(p)$. Since $\pi_1(M)$ is an \iccg\ (see \ref{i see cg}), the elements $u_1$ and $u_2$ do not commute; and since the $k$-free group $\pi_1(M)$ is in particular $2$-free, $u_1$ and $u_2$ are independent. The conclusion now follows from the case $m=2$ of \ref{reformulated big radius corollary from cusp}.
\EndProof

\Number\label{fitz}
This subsection is in large part devoted to review of material from \cite{fourfree}, \cite{rosemary} and \cite{gs}. 
One minor difference between the conventions that we give here and those used in the cited papers is that we use the notation $d(x,P)$ that was introduced in 
\ref{nbhd}; another minor difference will be mentioned below.

Let $\Gamma$ be a 
cocompact,
discrete subgroup of $\isomplus(\HH^3)$. We observed in \ref{i see cg} that
$\frakX_\lambda(\Gamma)=\bigcup_{C\in \calc_\lambda(\Gamma)} Z_{\lambda}(C)$. Hence we may regard  $ \mathcal{Z}_{\lambda}(\Gamma)$ as a covering  of $\frakX_\lambda(\Gamma)$, in the sense discussed in \ref{complex stuff}. It should be borne in mind that the index set for this covering is $\calc_\lambda(\Gamma)$, a collection of maximal cyclic groups of $\Gamma$.
 The nerve $K_{{\mathcal Z}_{\lambda}(\Gamma)}$ is well defined by \ref{complex stuff}.

(In the corresponding discussion in \cite[Section 2]{gs} it was assumed that $\frakX_\lambda(\Gamma)=\HH^3$. 
Although, for the sake of a smooth exposition, we have not built this condition into the general machinery described in this subsection, the condition will hold in the context of some of the results which build on this machinery and are quoted in the final application, Theorem \ref{before key}. More specifically, the condition is shown to hold in the context of the proof of Proposition \ref{hum wow wow} below, where it is used to justify quoting Proposition 2.13 of \cite{gs}. We have avoided quoting \cite{gs} in contexts where this condition does not necessarily hold; 
in a few places---all of them in the present subsection---this has involved pointing out proofs of statements
similar to those proved in \cite{gs}, but fortunately the statements in question are almost trivial.)

Recall from \ref{complex stuff} that the nerve $K_{{\mathcal Z}_{\lambda}(\Gamma)}$ comes equipped with a bijection $C\mapsto v_C$ from $\calc_\lambda(\Gamma)$, the index set of the covering ${\mathcal Z}_{\lambda}(\Gamma)$ of 
$\frakX_\lambda(\Gamma)$, to the vertex set 
(or $0$-skeleton)
$K_{{\mathcal Z}_\lambda(\Gamma)}^{(0)}$ of $K_{{\mathcal Z}_\lambda(\Gamma)}$.   We will denote the inverse bijection by $v\mapsto C_v$. 

According to the conventions of \ref{nervy notation}, for any simplex 
$\sigma$ of $K_{{\mathcal Z}_\lambda(\Gamma)}$, the set $\cals _\sigma=\cals _\sigma^{{\mathcal Z}_{\lambda}(\Gamma)}$
%{\calf } $ (or by $\cals _\sigma$ when the covering %$\calf$ is understood) the set 
consists of all maximal cyclic subgroups $C\in 
\calc_\lambda(\Gamma)$ such that $v_C$ is a vertex of $\sigma$.
Equivalently,
%Thus for every simplex $\sigma$ of $K_{{\mathcal Z}_{\lambda}(\Gamma)}$, the definition of 
$\cals_\sigma$
%=\cals_\sigma^
%{{\mathcal Z}_{\lambda}(\Gamma)}$ given in \ref{nervy notation} 
may be described as the set of all indices $C_v$, where $v$ ranges over the vertices of $\sigma$. 

For each simplex $\sigma$ we will set $\Theta(\sigma)=\langle\cals_\sigma\rangle$, so that $\Theta(\sigma)$ is a finitely generated subgroup of $\Gamma$.

If $W$ is a saturated subset of $K_{{\mathcal Z}_{\lambda}(\Gamma)}$, we define $\Theta(W)$ to be the subgroup of $\Gamma$ generated by all the subgroups $\Theta(\sigma)$ where $\sigma$ ranges over the simplices contained in $W$. (The definitions of $\Theta(\sigma)$ and $\Theta(W)$ given here are equivalent to the definitions given in \cite{gs}, although the latter are framed in terms of the notion of a ``labeled complex,'' which will not be used in this paper.)

Note that for any simplex $\sigma$ of $K_{{\mathcal Z}_{\lambda}(\Gamma)}$, the group $\Theta(\sigma)$ is by definition generated by a non-empty collection of maximal cyclic subgroups of $\Gamma$, and is in particular non-trivial. It follows that $\Theta(W)$ is a non-trivial subgroup of $\Gamma$ for every non-empty saturated subset  $W$ of $K_{{\mathcal Z}_{\lambda}(\Gamma)}$.

Note that if $\tau$ is a face of a simplex $\sigma$ of $K_{{\mathcal Z}_{\lambda}(\Gamma)}$, we have $\cals_\tau\subset\cals_\sigma$ and therefore $\Theta(\tau)\le\Theta(\sigma)$.

The conventions of \ref{nervy notation}  give
\Equation\label{that's what it means}
\calu_\sigma=
\calu_\sigma^{
K_{{\mathcal Z}_{\lambda}(\Gamma)}}
=
\bigcap_{C\in\cals_\sigma}Z_\lambda(C)
\EndEquation
for each simplex $\sigma$ of $K_{{\mathcal Z}_{\lambda}(\Gamma)}$.

We have observed that 
the hyperbolic cylinder
$Z_{\lambda} ({ C})$ is convex for each $C\in \calc_\lambda(\Gamma)$; hence for each simplex $\sigma$ of $K_{{\mathcal Z}_{\lambda}(\Gamma)}$, the equality (\ref{that's what it means}) exhibits $\calu_\sigma$ as an intersection of convex sets. 
Thus:
\Claim\label{convex}
The non-empty set
 $\calu_\sigma$ is convex, and in particular contractible,
for each simplex $\sigma$ of $K_{{\mathcal Z}_{\lambda}(\Gamma)}$.
\EndClaim
(The observation \ref{convex} is  related to the proof of \cite[Lemma 3.3]{gs}.)

Since $C\mapsto v_C$ is a bijection, we may define an action of $\Gamma$ on $K_{{\mathcal Z}_\lambda(\Gamma)}^{(0)}$ by setting $ x\cdot v_C=v_{ x C x^{-1}}$ for each $ x\in\Gamma$ and each $C\in\calc_\lambda(\Gamma)$. Under this  action we have 
\Equation\label{volcano}
C_{x\cdot v}=xC_vx^{-1} \text{ for every } x\in\Gamma \text{ and every }v\in K_{{\mathcal Z}_\lambda(\Gamma)}^{(0)}. 
\EndEquation

Note also that for every $C\in\calc_\lambda(\Gamma)$ and every $x\in\Gamma$, we have
\Equation\label{corona corona}
Z_\lambda(xCx^{-1})=x\cdot Z_\lambda(C).
\EndEquation

If $\sigma$ is any simplex of $K_{{\mathcal Z}_\lambda(\Gamma)}$, using (\ref{that's what it means}),
(\ref{corona corona}) and (\ref{volcano}) we find that
$$\emptyset\ne
x\cdot\calu_\sigma=
\bigcap_{C\in\cals_\sigma}x\cdot Z_\lambda(C)
=\bigcap_{C\in\cals_\sigma} Z_\lambda(xCx^{-1})=
\bigcap_{v\le\sigma}
%{C\in\cals_\sigma} 
Z_\lambda(xC_vx^{-1})
=\bigcap_{v\le\sigma}
%{C\in\cals_\sigma} 
Z_\lambda(C_{x\cdot v}),
$$
where the notation $v\le\sigma$ means that $v$ is a vertex of $\sigma$. This shows that the set of vertices of the form $x\cdot v$, where $v$ is a vertex of $\sigma$, is itself the vertex set of a simplex of $K_{{\mathcal Z}_\lambda(\Gamma)}$. Thus the action of $\Gamma$ on $K_{{\mathcal Z}_\lambda(\Gamma)}^{(0)}$ extends to a simplicial 
action on $K_{{\mathcal Z}_\lambda(\Gamma)}$, which will be referred to as the {\it canonical action} of $\Gamma$ on 
$K_{{\mathcal Z}_{\lambda}(\Gamma)}$ (cf. \cite[Subsection 2.12]{gs}).

It follows from (\ref{volcano}) that under the canonical action we have 
\Equation\label{expo}
\Theta(x\cdot W)=x\Theta(W)x^{-1}
\EndEquation
for every $x\in\Gamma$ and for every saturated set $W\subset K_{{\mathcal Z}_{\lambda}(\Gamma)}$.

\EndNumber

\Proposition\label{hum wow wow}
Let $\Gamma$ be a torsion-free, discrete, 
cocompact 
subgroup of ${\rm Isom}_+(\bf H^3)$, and let $\lambda$ be a positive real number. 
Suppose that 
for every point $p$ in the manifold $\HH^3/\Gamma$ 
we have $\nextone(p)<\lambda$.
Set $K=K_{{\mathcal Z}_{\lambda}(\Gamma)}$. Then $|K|-|K^{(0)}|$ is contractible
(where as in \ref{complex stuff}, $K^{(0)}$ denotes the $0$-skeleton of $K$).
\EndProposition

\Proof
Set $M=\HH^3/\Gamma$. According to \ref{short and next} we have $\shortone( p )\le\nextone( p )$ for every $ p \in M$. The hypothesis therefore implies that
for every point $p$ in  $M$ we have $\shortone(p)<\lambda$. 
According to \ref{short and next}, this means that 
$\Mthin(\lambda)=M$, and hence that
$\frakX_\lambda(\Gamma)=\HH^3$. Thus
${\mathcal Z}_{\lambda}(\Gamma)$ is a covering of $\HH^3$, and $K$ is its nerve.

We will apply Proposition \ref{linky corollary}, with $\HH^3$ and $\mathcal{Z}_{\lambda}(\Gamma)= (Z_{\lambda}(C))_{C\in \calc_\lambda(\Gamma)}$ playing the respective roles of the space $X$ and the (indexed) covering $\calf=(U_i)_{i\in I}$ in that proposition. We will show that Conditions (1)--(3) of the hypothesis of Proposition \ref{linky corollary} hold in this setting.

To verify Condition (1), consider an arbitrary point $P$ of $\HH^3$. If $p$ denotes the image of $P$ under the quotient projection from $\HH^3$ to $M:=\HH^3/\Gamma$, then by the hypothesis of the present proposition we have $\nextone(p)<\lambda$. According to \ref{what about the workers},
this means that $p$ is not a $\lambda$-semithick point of $M$; that is, there exist two loops of length less than $\lambda$ based at $p$ which define non-commuting elements of $\pi_1(M,p)$. Hence there exist non-commuting elements $x_1$ and $x_2$ of $\Gamma$ such that $d(x_m,P)<\lambda$ for $m=1,2$. Since $\Gamma$ is an \iccg\ (see \ref{i see cg}), each $x_m$ belongs to a maximal cyclic subgroup $C_m$ of $\Gamma$, and we have $C_1\ne C_2$ since $x_1$ and $x_2$ do not commute. For $m=1,2$, since $x_m$ is a non-trivial element of $C_m$ and $d(x_m,P)<\lambda$, we have $C_m\in{\mathcal C}_{\lambda}(\Gamma)$, and $P\in Z_\lambda(C_m)$. 
To summarize, $C_1$ and $C_2$ are distinct elements of the index set $\calc_\lambda(\Gamma)$, and the given point $P\in\HH^3$ lies in $ Z_\lambda(C_1)\cap Z_\lambda(C_2)$. This establishes Condition (1) of Proposition \ref{linky corollary} in the present setting.

According to \ref{convex},
the non-empty set
$\calu_\sigma$
is   convex and hence contractible 
for each simplex 
$\sigma$ 
of $K=K_{{\mathcal Z}_{\lambda}(\Gamma)}$.
Thus Condition (2) of Proposition \ref{linky corollary} holds in the present setting.

Since $\Gamma$ is cocompact, discrete and torsion-free,
and since $\frakX_\lambda(\Gamma)=\HH^3$,
it follows from \cite[Proposition 2.13]{gs} that $K$ is finite-dimensional.
Thus Condition (3) of Proposition \ref{linky corollary} holds in the present setting.
 
We may therefore apply Proposition \ref{linky corollary} to deduce that $|K|-|K^{(0)}|$ is homotopy equivalent to $\HH^3$, and is therefore contractible.
\EndProof

The following definition is needed for the statement of Proposition \ref{restricting prop}, which will play a key role in this paper.

\Definition\label{restricting def}
Let $\Gamma$ be a discrete, 
cocompact
subgroup of ${\rm Isom}_+(\bf H^3)$, 
let $\lambda$ be a positive real number, and let $r$ be a positive integer.

Suppose that $\phi$ is a simplex of $K_{{\mathcal Z}_{\lambda}(\Gamma)}$ and that $v$ is a vertex of $\phi$. We will say that $v$ is 
an {\it $r$-\restraining\ vertex}
of $\phi$ if, for every face $\tau$ of $\phi$ such that $v$ is a vertex of $\tau$, we have 
$\rank(\Theta(\tau))\le r$.
\EndDefinition

\Proposition\label{restricting prop}
Let 
$r$ be a positive integer, let
$\Gamma$ be an $(r+1)$-free, 
discrete, 
cocompact 
subgroup of ${\rm Isom}_+(\bf H^3)$, and let 
$\lambda\ge\log(2r+1)$ be a real number. 
Suppose that 
$\shortone(p)<\f_r(\lambda)$
for every point $p$ of the manifold  $\HH^3/\Gamma$ 
(or equivalently that $\Mthin(\fone(\lambda))=M$, cf. \ref{short and next}).
Set $K=K_{{\mathcal Z}_{\lambda}(\Gamma)}$.  Then for every simplex $\sigma$ of $K$ there is a simplex $\phi$ of $K$ such that (i) $\sigma$ is a face of $\phi$, and (ii) $\phi$ has at least one 
$r$-\restraining\ vertex.
\EndProposition

\Proof
Set $\nu=\f_r(\lambda)$. 
According to 
\ref{Q-def},
$\f_r$ 
is strictly monotone decreasing on 
$(\log5,\infty)$,
and 
$\f_r(\log(2r+1))=\log(2r+1)$. Hence the condition $\lambda\ge\log(2r+1)$ implies that $\nu=\f_r(\lambda)\le\log(2r+1)$,
and in particular that $\nu\le\lambda$.

Choose a point $P\in \calU_\sigma$,
and let $p\in\HH^3/\Gamma$ denote the image of $P$ under the quotient projection. 
 Since the hypothesis implies that
$\shortone(p)<\nu$,
there is an element $y_0$ of $\Gamma$ such that $d(y_0,P)<\nu$. Let $C_0$ denote the maximal cyclic subgroup of $\Gamma$ containing $y_0$. Since 
$\nu\le\log(2r+1)\le\lambda$, 
we in particular have 
$d(y_0,P)<\lambda$;
hence $C_0\in{\mathcal C}_{\lambda}(\Gamma)$ and $P\in Z_{\lambda}(C_0)$. Set $v_0=v_{C_0}$. The set $\calU_\sigma\cap Z_{\lambda}(C_0)$ contains $P$ and in particular is non-empty; this implies that there is a simplex $\phi$ of $K$ whose vertex set is the union of $\{v_0\}$
with the vertex set of $\sigma$. Thus $\sigma$ is a face of $\phi$. (It is not necessarily a proper face, as it may happen that $v_0$ is a vertex of $\sigma$.) We have $P\in Z_{\lambda}(C_0)\cap\calU_\sigma=\calU_\phi$.

We will complete the proof by showing that $v_0$ is 
an $r$-\restraining\ vertex 
of $\phi$. Suppose that $\tau$ is a  face of $\phi$ having $v_0$ as a vertex. We must show that 
$\rank(\Theta(\tau))\le r$. 
Assume to the contrary that 
$\rank(\Theta(\tau))\ge r+1$.

For each vertex $v$ of 
$\phi$, 
choose a generator $x_v$ of the infinite cyclic subgroup $C_v$ of $\Gamma$; in particular the $x_v$ are non-trivial elements of $\Gamma$. Set $x_0=x_{v_0}$,
so that $C_{v_0}=C_0$.
Let $S$ denote the set consisting of all elements $x_v$ where $v$ ranges over the vertices of $\tau$. We have $\Theta(\tau)=\langle S\rangle$, so that $\rank \langle S\rangle\ge r+1$. The group $\Gamma$ is 
$(r+1)$-free 
by hypothesis.
Proposition \ref{new if the right hand don't get ya}, applied with $k=r+1$, therefore gives elements $x_1,\ldots,x_r$ of $S$ such that $x_0, x_1,\ldots,x_r$ are independent. For $i=1,\ldots,r$ 
we have $x_i=x_{v_i}$ for some vertex $v_i$ of $\tau$.

Since $\tau$ is a face of $\phi$, we have $\calU_\tau\supset\calU_\phi$ by \ref{nervy notation}. In particular, $P\in\calU_\tau$, and hence $P\in Z_{\lambda}(C_{v_i})$ for 
$i=1,\ldots,r$. We may therefore choose, for $i=1,\ldots,r$, 
a non-trivial element $y_i$ of $C_{v_i}$ such that $d(y_i,P)<\lambda$. 
We have already seen that $y_0\in C_0=C_{v_0}=\langle x_0\rangle$. Thus for 
$i=0,1,\ldots,r$, 
the 
infinite cyclic group $C_{v_i}$ is generated by $x_i$ and has $y_i$ as a non-trivial element. We may write $y_i=x_i^{e_i}$ for 
$i=0,1,\ldots,r$, 
where the $e_i$ are non-zero integers. Since 
$x_0,\ldots,x_r$ are independent and $e_0,\ldots,e_r$ are non-zero, the elements $y_0,\ldots,y_r$ 
of $\Gamma$ are also independent. As $\Gamma\le\isomplus(\HH^3)$ is discrete, it follows from Theorem \ref{our log(2k-1)} that 
$$
\sum_{i=0}^r
\calq(y_i,P)
%+\calq(y_1,P)+\calq(y_2,P)+\calq(y_3,P)
\le\frac12,$$
i.e. 
$$
\sum_{i=0}^r
Q(d(y_i,P))
%+Q(d(y_1,P))+Q(d(y_2,P))+Q(d(y_3,P))
\le\frac12.$$
Since $d(y_0,P)<\nu$, and $d(y_i,P)<\lambda$ for $i=1,2,\ldots,r$, 
and since $Q$ is strictly monotone decreasing, it follows that 
$Q(\nu)+rQ(\lambda)<1/2$. But the definition of $\nu=\f_r(\lambda)$ implies that $Q(\nu)+rQ(\lambda)=1/2$. 
This contradiction completes the proof.
\EndProof

\Corollary\label{first restricting cor}
Let $\Gamma$ be a $4$-free, discrete, 
cocompact
subgroup of ${\rm Isom}_+(\bf H^3)$, and let $\lambda\ge\log7$ be a real number. 
Suppose that 
$\shortone(p)<\fone(\lambda)$
for every point $p$ of the manifold  $\HH^3/\Gamma$ 
(or equivalently that $\Mthin(\fone(\lambda))=M$).
Set $K=K_{{\mathcal Z}_{\lambda}(\Gamma)}$.  Then for every simplex $\sigma$ of $K$,  the subgroup $\Theta(\sigma)$ is closable (see \ref{closable def}), and we have (in the notation of \ref{wuzza closure}) $\localrank(\clo(\Theta(\sigma)))\le3$.
\EndCorollary

\Proof 
The hypotheses of the corollary imply those of
%Since $\Gamma$ is $4$-free, we may apply 
Proposition \ref{restricting prop} with $r=3$. Hence, if $\sigma$ is any simplex of $K$, there is
a simplex $\phi$ of $K$ such that $\sigma$ is a face of $\phi$, and  $\phi$ has a 
$3$-\restraining\ vertex, 
say $v_0$ (which may or may not be a vertex of $\sigma$). Let $\omega$ denote the face of $\phi$ whose vertices are $v_0$ and the vertices of $\sigma$; then $\sigma$ is a (not necessarily proper) face of $\omega$. If we fix a generator $x_0$ of $C_{v_0}\le\Gamma$, we have $\Theta(\omega)=\langle\Theta(\sigma),x_0\rangle$. 

Since $v_0$ is a 
$3$-\restraining\ vertex 
for $\phi$, we have $\rank(\Theta(\omega))\le3$
(see Definition \ref{restricting def}).
On the other hand, the hypothesis that $\Gamma$ is $4$-free says that $\iof(\Gamma)\ge4$. Hence  $
\rank(\langle\Theta(\sigma),x_0\rangle)=
\rank(\Theta(\omega))<\iof(\Gamma)$, so that Condition (b) of Definition \ref{closable def} holds with $\Theta(\sigma)$ and $x_0$ playing the respective roles of $A$ and $z$. Condition (a) of of Definition \ref{closable def} also holds with $A=\Theta(\sigma)$, since $\Theta(\sigma)$ is by definition of finite rank. This shows that $\Theta(\sigma)$ is closable, which is the first assertion of the corollary.

To prove the second assertion, that $\localrank(\clo(\Theta(\sigma)))\le3$, we first note that according to Assertion (1) of Proposition \ref{closure properties}, the local rank of the group $\clo(\Theta(\sigma))$ is at most $\rank (\Theta(\sigma))$. Hence the second assertion of the present corollary holds in the case where $\rank (\Theta(\sigma))\le3$. 

It remains to prove the second assertion of the corollary in the case where $\rank(\Theta(\sigma))>3$. Since $
\rank(\langle\Theta(\sigma),x_0\rangle)=
\rank(\Theta(\omega))\le3$, we have $\rank(\langle\Theta(\sigma),x_0\rangle)<
\rank(\Theta(\sigma))$ in this case; in particular, Condition (b) of Definition  \ref{admissible elt} holds with $\Theta(\sigma)$ and $x_0$ playing the respective roles of $A$ and $x$. Since we have also seen that $
\rank(\langle\Theta(\sigma),x_0\rangle)<\iof(\Gamma)$, the subgroup 
$
\langle\Theta(\sigma),x_0\rangle)$ of $\Gamma$ is by definition strongly closable; thus Condition (a)  of Definition  \ref{admissible elt} also holds with $\Theta(\sigma)$ and $x_0$ playing the roles of $A$ and $x$. This shows that $x_0$ is
$\Theta(\sigma)$-\admissible. Thus $x_0$ is a length-one $\Theta(\sigma)$-\admissible\  sequence, and therefore $x_0\in\clo(\Theta(\sigma))$. Since, by Assertions (1) and (2) of Proposition \ref{closure properties}, $\clo(\Theta(\sigma))$ is a subgroup of $\Gamma$ containing $\Theta(\sigma)$, we have 
$$\Theta(\sigma)\le\langle \Theta(\sigma),x_0\rangle\le\clo(\Theta(\sigma)).$$
We have shown that $\Theta(\sigma)$ is closable, and that $\langle \Theta(\sigma),x_0\rangle$ is strongly closable. Thus the hypotheses of
Corollary \ref{properties cor} hold with 
$\Theta(\sigma)$  and $\langle \Theta(\sigma),x_0\rangle=\Theta(\omega)$ playing the respective roles of $A$ and $B$. It follows that $\clo(\Theta(\sigma))=\clo( \Theta(\omega))$. But according to Assertion (1) of Proposition \ref{closure properties}, we have $\localrank(\clo(\Theta(\omega))) \le\localrank(\Theta(\omega))$, and we have seen that $\rank(\Theta(\omega))\le3$;  hence we have $\localrank(\clo(\Theta(\sigma)))\le3$, and the second assertion is proved in this case as well.
\EndProof

\Corollary\label{second restricting cor}
Let $\Gamma$ be a $4$-free, discrete, 
cocompact
subgroup of ${\rm Isom}_+(\bf H^3)$, and let $\lambda\ge\log7$ be a real number. 
Suppose that 
for every point $p$ in the manifold $\HH^3/\Gamma$ 
we have $\shortone(p)<\fone(\lambda)$
(or equivalently that $\Mthin(\fone(\lambda))=M$).
Set $K=K_{{\mathcal Z}_{\lambda}(\Gamma)}$. Suppose that $\sigma$ is a simplex of $K$ and that $\tau$ is a face of $\sigma$ (so that $\Theta(\sigma)$ and $\Theta(\tau)$ are closable by Corollary \ref{first restricting cor}). Suppose that $\localrank(\clo(\Theta(\tau)))=3$. Then $\clo(\Theta(\sigma))=\clo(\Theta(\tau))$.
\EndCorollary

\Proof
The hypotheses of the corollary imply those of
%Since $\Gamma$ is $4$-free, we may apply 
Proposition \ref{restricting prop} with $r=3$. Hence there is
a simplex $\phi$ of $K$ which has $\sigma$ as a face and has a 
$3$-\restraining\ 
vertex $v_0$. Let  $\omega$ denote the face of $\phi$ whose vertex set is the union of $\{v_0\}$ with the vertex set of $\tau$. (This is not necessarily a disjoint union.) For each vertex $v$ of $\phi$ let us choose a generator $x_v$ of the cyclic group $C_v\le\Gamma$. Set $x_0=x_{v_0}$, and note that $\Theta(\omega)=\langle\Theta(\tau),x_0\rangle$. 

According to Assertion (1) of Proposition \ref{closure properties}, and the hypothesis of the present corollary, we have $\rank(\Theta(\tau))=\localrank(\Theta(\tau))\ge \localrank(\clo(\Theta(\tau)))=3$. On the other hand, since $v_0$ is a 
$3$-\restraining\ 
vertex of $\phi$, and $\omega$ is a face of $\phi$, we have $\rank(\langle\Theta(\tau),x_0\rangle)=
\rank(\Theta(\omega))\le3$
(see Definition \ref{restricting def}).
Thus Condition (b) of Definition \ref{admissible elt} holds if we let $\Theta(\tau)$ and $x_0$ play the respective roles of $A$ and $x$. Furthermore, since $\Gamma$ is $4$-free we have $\iof(\Gamma)\ge4$; since $\rank(
\langle\Theta(\tau),x_0\rangle)\le3$, the subgroup $\langle\Theta(\tau),x_0\rangle$ of $\Gamma$ is strongly closable, and thus Condition (a) of Definition \ref{admissible elt} also holds with these choices of $A$ and $x$. This shows that $x_0$ is $\Theta(\tau)$-\admissible. 

We may therefore regard $x_0$ as a one-term $\Theta(\tau)$-\admissible\ sequence, so that by definition we have $x_0\in\clo(\Theta(\tau))$. By Assertions (1) and (2) of Proposition \ref{closure properties}, $\clo(\Theta(\tau))$ is a subgroup of $\Gamma$ containing $\Theta(\tau)$. Hence 
$\Theta(\tau)\le\langle\Theta(\tau),x_0\rangle\le\clo(\Theta(\tau))$, i.e.
$$\Theta(\tau)\le\Theta(\omega)\le\clo(\Theta(\tau)).$$
Furthermore, the subgroups $\Theta(\tau)$ and $\Theta(\omega)$ of $\Gamma$ are both closable according to the first assertion of Corollary \ref{first restricting cor}. We may therefore apply Corollary \ref{properties cor}, with $\Theta(\tau)$ and $\Theta(\omega)$ playing the respective roles of $A$ and $B$ in that corollary, to deduce that $\clo(\Theta(\tau))=\clo(\Theta(\omega))$. 

According to Assertion (1) of Proposition \ref{closure properties},  we have $\rank(\Theta(\omega))=\localrank(\Theta(\omega))\ge \localrank(\clo(\Theta(\omega)))$. But we have seen that $\clo(\Theta(\omega))=\clo(\Theta(\tau))$, and by hypothesis we have $\localrank(\clo(\Theta(\omega)))=3$. Hence $\rank(\Theta(\omega))\ge3$. As
we have already observed that
$\rank(\Theta(\omega))\le3$, we have $\rank(\Theta(\omega))=3$.

Now consider an arbitrary vertex $v$ of $\phi$. Let $\psi_v$ denote the face of $\phi$ whose vertex set is the  union of $\{v\}$ with the vertex set of $\omega$. (Again, this union is not necessarily disjoint.) We have $\Theta(\psi_v)=\langle\Theta(\omega),x_v\rangle$. The vertex $v_0$ of $\omega$ is in particular a vertex of $\psi_v$, and is a  
$3$-\restraining\ vertex 
of $\phi$; hence $\rank(\Theta(\psi_v))\le3$, i.e. $\rank(\langle\Theta(\omega),x_v\rangle)\le3=\rank(\Theta(\omega))$. Thus Condition (b) of \ref{admissible elt} holds with $\Theta(\omega)$ and $x_v$ playing the respective roles of $A$ and $x$. Since $\iof(\Gamma)\ge4$, Condition (a) of \ref{admissible elt} also holds with these choices of $A$ and $x$. This shows that $x_v$ is $\Theta(\omega)$-\admissible\ for each vertex $v$ of $\phi$.

Thus for each vertex $v$ of $\phi$ we may regard $x_v$ as a one-term $\Theta(\omega)$-\admissible\ sequence, and hence we have $x_v\in\clo(\Theta(\omega))=\clo(\Theta(\tau))$. But $\Theta(\phi)$ is generated by the elements $x_v$ as $v$ ranges over the vertex set of %$\sigma$, which is in turn a subset of the vertex set of 
$\phi$. Since $\Theta(\phi)$ is generated by a family of elements of
$\clo(\Theta(\tau))$, which by Assertion (1) of Proposition \ref{closure properties} is a subgroup of $\Gamma$, we have $\Theta(\phi)\le\clo(\Theta(\tau))$. Now since $\sigma$ is a face of $\phi$, we have $\Theta(\sigma)\le\Theta(\phi)$ and hence 
$\Theta(\sigma)\le\clo(\Theta(\tau))$.

On the other hand,
since $\tau$ is a face of  $\sigma$ we also have $\Theta(\tau)\le\Theta(\sigma)$. As was pointed out in the statement of the present corollary, it follows from Corollary \ref{first restricting cor} that $\Theta(\tau)$ and $\Theta(\sigma)$ are closable subgroups of $\Gamma$. We may now apply 
 Corollary \ref{properties cor}, with $\Theta(\tau)$ and $\Theta(\sigma)$ playing the respective roles of $A$ and $B$, to deduce that
$\clo(\Theta(\tau))=\clo(\Theta(\sigma))$, as required.
\EndProof

\Lemma\label{commies}
Let $\Gamma$ be a $4$-free, discrete, 
cocompact 
subgroup of ${\rm Isom}_+(\bf H^3)$, and let $\lambda\ge\log7$ be a real number. 
Suppose that 
for every point $p$ in the manifold $\HH^3/\Gamma$ 
we have $\shortone(p)<\fone(\lambda)$
(or equivalently that $\Mthin(\fone(\lambda))=M$).
Set $K=K_{{\mathcal Z}_{\lambda}(\Gamma)}$. Suppose that $r$ is an integer equal to either $2$ or $3$, and that $X$ is a connected, saturated subset of $|K|$
such that for every simplex $\sigma\subset X$ we have $\localrank(\clo(\Theta(\sigma)))=r$. (Here the group $\Theta(\sigma)$ is closable by the first assertion of Corollary \ref{first restricting cor}, so that $\clo(\Theta(\sigma))$ is a well- defined subgroup of $\Gamma$ by Assertion (1) of Proposition \ref{closure properties}.) Then there is a subgroup $C$ of $\Gamma$ such that $\Theta(X)\le C$ and $\localrank(C)=r$.
\EndLemma

\Proof
We will say that a simplex $\sigma$ of $K$ is {\it incident} to a simplex $\sigma'$ of $K$ if 
either $\sigma$ is a face of $\sigma'$ or $\sigma'$ is a face of $\sigma$.

Since $\Gamma$ is finitely generated, the set $\calc_\lambda(\Gamma)$ is countable, and hence $K$ has only countably many simplices. Hence we may write the saturated set $X\subset|K|$ in the form $\bigcup_{i\in I}\sigma_i$, where the index set $I$ is either the set of all non-negative integers or the set $\{0,\ldots,N\}$ for some $N\ge0$. Since  $X$ is connected, we may take the $\sigma_i$ to be indexed in such a way that for every $i\in I-\{0\}$ there is an index $l\in I$, with $l<i$, such that  $\sigma_i$ is incident to $\sigma_{l}$.

Consider the case $r=3$. If $i$ and $l$ are indices in $I$ such that $\sigma_i$ is incident to $\sigma_l$, then since $\localrank(\clo(\Theta(\sigma_i)))=\localrank(\clo(\Theta(\sigma_l)))=3$ it follows from Corollary \ref{second restricting cor} that 
$\clo(\Theta(\sigma_i))=\clo(\Theta(\sigma_l))$. Our indexing of the $\sigma_i$ therefore guarantees, by strong induction, that $\clo(\Theta(\sigma_i))=\clo(\Theta(\sigma_0))$ for every $i\in I$. According to Assertion (2) of Proposition \ref{closure properties}, 
$\Theta(\sigma_i)$ is a subgroup of
$\clo(\Theta(\sigma_i))=\clo(\Theta(\sigma_0))$
for each $i\in I$. Since $\Theta(X)$ is generated by the subgroups $\Theta(\sigma_i)$ for $i\in I$, we have $\Theta(X)\le \clo(\Theta(\sigma_0))$. Since $\localrank \clo(\Theta(\sigma_0))=r=3$, the conclusion of the lemma follows in this case upon setting $C=\clo(\Theta(\sigma_0))$.

We now turn to the case $r=2$. In this case, for every $n\in I$ we set \linebreak $C_n=\langle \clo(\Theta(\sigma_0)),\ldots,\clo(\Theta(\sigma_n))\rangle\le\Gamma$.
We claim that
\Equation\label{lomax}
\localrank(C_n)=2\text{ for every }n\in I.
\EndEquation
The proof of (\ref{lomax}) proceeds by induction on $n\in I$. For the base case, we observe that $C_0=\clo(\Theta(\sigma_0))$, so that $\localrank(C_0)=r=2$. Now suppose that $n\in I-\{0\}$ is given, and that $\localrank(C_{n-1})=2$. Choose an index $l\in I$, with $l<n$, such that $\sigma_n$ is incident to $\sigma_l$. 

Consider first the subcase in which $\sigma_n$ is a face of $\sigma_l$. We then have $\Theta(\sigma_n)\le\Theta(\sigma_l)$, and Assertion (3) of Proposition \ref{closure properties} gives $\clo(\Theta(\sigma_n))\le\clo(\Theta(\sigma_l))\le C_{n-1}$, so that $C_{n}=\langle C_{n-1}, \clo(\Theta(\sigma_n))\rangle=C_{n-1}$. In particular, $\localrank( C_n)=\localrank (C_{n-1})=2$, and the induction is complete in this subcase.

Now consider the subcase in which $\sigma_l$ is a face of $\sigma_n$. Then $\Theta(\sigma_l)\le\Theta(\sigma_n)$, and Assertion (3) of Proposition \ref{closure properties} gives $\clo(\Theta(\sigma_l))\le\clo(\Theta(\sigma_n))$. By the definition of $C_{n-1}$ we also have $\clo(\Theta(\sigma_l))\le C_{n-1}$. Thus if we set $A=C_{n-1}$ and $B=\clo(\Theta(\sigma_n))$, we have $A\cap B\ge \clo(\Theta(\sigma_l))$. Since $\localrank(\clo(\Theta(\sigma_l)))=r=2$, the group $\clo(\Theta(\sigma_l))$ is non-abelian, and hence $A\cap B$ is non-abelian. The groups $A$ and $B$ are subgroups of the $4$-free group $\Gamma$. We have $\localrank(A)=\localrank(C_{n-1})=2$ by the induction hypothesis, and $\localrank(B)=\localrank(\clo(\Theta(\sigma_n)))=r=2$. Thus all the hypotheses of Proposition \ref{general klm} hold with these choices of $A$ and $B$, and the latter proposition therefore gives $\localrank(C_n)=\localrank(\langle C_{n-1}, \clo(\Theta(\sigma_n))\rangle)=\localrank(\langle A, B\rangle)=2$. This completes the induction, and (\ref{lomax}) is established.

To complete the proof of the lemma in the case $r=2$, we set $C=C_N$ in the subcase where $I=\{0,\ldots,N\}$ for some $N\ge0$, and define $C$ to be the monotone union $\bigcup_{n\ge0}C_n$ in the subcase where $I$ is the set of all non-negative integers. In the first subcase, (\ref{lomax}) directly implies that $\localrank(C)=2$; in the second subcase, it follows from (\ref{lomax}) that $C$ is a monotone union of subgroups having local rank $2$, and is therefore of local rank $2$. In either case, it follows from the construction of $C$ that $C$ has $\clo(\Theta(\sigma_i))$ as a subgroup for each $i\in I$; hence by Assertion (2) of Proposition \ref{closure properties},
$C$ has $\Theta(\sigma_i)$ as a subgroup for each $i\in I$. But $\Theta(X)$ is generated by the subgroups $\Theta(\sigma_i)$  for each $i\in I$, and hence $\Theta(X)\le C$.
\EndProof

\section{
The central result
}\label{central section}

\Theorem\label{before key}
Let  
$M$ be a 
closed, 
orientable hyperbolic $3$-manifold such that
$\pi_1(M)$ is $4$-free. Let $\lambda\ge\log7$ be a real
number
(so that in particular $\lambda$ is in the domain of $\fone$ (see \ref{Q-def})).
Then
either there  is a point $p\in M$ with $\nextone(p)\ge\lambda$, or
there is a point  $p\in M$ with $\shortone(p) 
\ge
\fone(\lambda)$.
\EndTheorem

\Proof
Let us assume that the conclusion is false; that is, we assume  that for every point $p\in M$ we have $\nextone(p)<\lambda$ and $\shortone(p)
<
\fone(\lambda)$.

Let us write $M=\HH^3/\Gamma$, where $\Gamma\le\isomplus(\HH^3)$ is a discrete group which is 
cocompact and 
$4$-free. 
Since we have  $\nextone(p)<\lambda$
for every point 
$p\in M$,
Proposition \ref{hum wow wow} guarantees that $|K|-|K^{(0)}|$ is contractible.

Since we have $\shortone(p)<\fone(\lambda)$ for every point $p\in M$,
it follows from Corollary \ref{first restricting cor} that
for every simplex $\sigma$ of $K$,  the subgroup $\Theta(\sigma)$ is closable (see \ref{closable def}), and we have (in the notation of \ref{wuzza closure}) $\localrank(\clo(\Theta(\sigma)))\le3$.

If $\sigma$ is any simplex contained in $|K|-|K^{(0)}|$, then by definition we have $\dim\sigma>0$. Thus $\sigma$ has at least two vertices, and hence $\Theta(\sigma)$ is generated by two or more distinct maximal cyclic subgroups. Since $\Gamma$ is an \iccg\ (see \ref{i see cg}), it follows that $\Theta(\sigma)$ is non-abelian. Since $\clo(\Theta(\sigma))\ge\Theta(\sigma)$ by Assertion (2) of Proposition \ref{closure properties}, the group $\clo(\Theta(\sigma))$ is also non-abelian, and hence $\localrank (\clo(\Theta(\sigma)))>1$. As we have shown that
$\localrank(\clo(\Theta(\sigma)))\le3$ for every simplex $\sigma$ of $K$, it now follows that for each simplex $\sigma\subset|K|-|K^{(0)}|$ we have
$\localrank(\clo(\Theta(\sigma)))\in\{2,3\}$.
We may therefore write $|K|-|K^{(0)}|$ as a set-theoretical disjoint union of two saturated sets $X_2$ and $X_3$, where for $r=2,3$ we denote by $X_r$ the union of all simplices $\sigma$ of $K$ such that $\localrank(\clo(\Theta(\sigma)))=r$.

For each $r \in \{ 2, 3 \}$, let $\mathcal W_r$ denote the set of connected components of $X_r$. 
We will say that elements $W_{2}$ of $\calw_{2}$ and $W_{3}$ of $\calw_{3}$ are {\it adjacent} if  there  are simplices $\sigma_{2}\subset{ W}_{2}$ and
%which is a face of a simplex
$\sigma_{3}\subset{ W}_{3}$ such that 
either $\sigma_{2}$ is a face of $\sigma_{3}$ or $\sigma_{3}$ is a face of $\sigma_{2}$.
(It is easy to show using Corollary \ref{second restricting cor} that $\sigma_3$ cannot be a face of $\sigma_2$, but this does not matter for the proof.)

We construct an abstract bipartite graph $\mathcal G$ as follows: The vertex set $Y$ of $\mathcal G$ is a disjoint union $Y_{2}\,\Dot\cup\, Y_{3}$, where for $r=2,3$ the set $Y_r$ is a bijective copy of  $\mathcal W_{r}$ and is equipped with a bijection from $\calw_r$ to $Y_r$, which we denote $W\mapsto s_W$. An edge of $\calg$ is defined to be an unordered pair of the form $\{s_{W_{2}},s_{W_{3}}\}$, where $W_r$ is an element of $\calw_r$ for $r=2,3$, and $W_{2}$ and $W_{3}$ are adjacent. 

Let $T$ denote the geometric realization of $\calg$ (regarded as a simplicial complex of dimension at most $1$).
Because 
$|K|-|K^{(0)}|$
 is the set-theoretical disjoint union of the saturated subsets  $X_{2}$ and $ X_{3}$ of $|K|$, it follows from \cite[Lemma 5.12]{fourfree} that $|{T}|$ is a homotopy-retract of 
$|K|-|K^{(0)}|$.
%\shmel{K}{k-3}|$}.
%$|K|-|K_{(k-3)}'|$.
%the saturated subset $X_{2} \Dot \cup X_{3}$. 
%As noted above, Lemma \ref{ContractibilityLemma} gives 
Since 
$|K|-|K^{(0)}|$
%\shmel{K}{k-3}|$}
%$|K|-|K_{(k-3)}'|$ 
is contractible, it follows that
%that $X_{2} \Dot \cup X_{3}$ is both connected and simply connected, and therefore 
$T$ is a tree. 
%Set $T={\mathcal G} (X_{2}, X_{3})$.

We consider the canonical action of
$\Gamma$  on $K$ that was discussed in \ref{fitz}. 
According to (\ref{expo}),  for every
 $x\in\Gamma$ and for every saturated set $W\subset K_{{\mathcal Z}_{\lambda}(\Gamma)}$, we have $\Theta(x\cdot W)=x\Theta(W)x^{-1}$. In particular for every $x\in\Gamma$ 
and every simplex $\sigma$ of $K$ 
we have
$\Theta(x\cdot \sigma)=x\Theta(\sigma)x^{-1}$. Using Proposition \ref{conjugation and closure} (with $\Theta(\sigma)$ and $x$ playing the respective roles of $A$ and $u$ in that proposition), 
we then find $\clo(\Theta(x\cdot \sigma))=\clo(x\Theta(\sigma)x^{-1})=x\clo(\Theta(\sigma))x^{-1}$. In particular we have
$\localrank(\clo(\Theta(x\cdot\sigma)))=
\localrank(\clo(\Theta(\sigma)))$ for every $\sigma\in K$ and every $x\in\Gamma$. 
Consequently, $X_{2}$ and $X_{3}$ are invariant under the action of $\Gamma$. Hence for each $x\in\Gamma$, for each $r\in\{2,3\}$, and for each $W\in\calw_{r}$, we have $ x\cdot W\in\calw_r$. Thus we may define an action of $\Gamma$ on
$Y_r$  by setting $ x\cdot s_W=s_{ x\cdot W}$ for each $W\in\calw_r$. The actions of $\Gamma$ on $Y_2$ and $Y_3$ defined in this way give rise to an action of $\Gamma$ on $Y$; the subsets $Y_2$ and $Y_3$ are invariant under this action.

Since the canonical action of $\Gamma$ on $K$ is simplicial, it defines a continuous action of $\Gamma$ on $|K|$; hence for $r=2,3$, the restricted action of $\Gamma$ on $X_r$ gives rise to an action on $\calw_r$. The simplicial nature of the action on $K$ also implies that if $W_2\in\calw_2$ and $W_3\in\calw_3$ are adjacent, then $ x\cdot W_2$ and $ x\cdot W_3$ are adjacent for every $ x\in\Gamma$. Hence the action of $\Gamma$ on 
$Y$ that we have
defined 
extends to a simplicial action of $\Gamma$ on $\calg$, which gives rise to an action on $T$.

Since  the sets $Y_{2}, Y_{3}\subset\calg$ are  $\Gamma$-invariant, and since each edge of $\calg$ has one vertex in $Y_{2}$ and one in $Y_{3}$, the action of $\Gamma$ on $T$ has no inversions.

If $s$ is an arbitrary vertex of $\calg$, we may write $s=s_W$ for some $W\in\calw_r$, where $r\in\{2,3\}$. If $ x$ is any element of the vertex stabilizer $\Gamma_s\le\Gamma$, it follows from the definition of the action of $\Gamma$ on $\calg$ that $ x\cdot W=W$. Again using 
(\ref{expo})
we deduce that  
$\Theta( W)= x\Theta(W) x^{-1}$ for every $x\in\Gamma_s$. 
This shows that $\Gamma_s$ is contained in the normalizer of $\Theta(W)$ in $\Gamma$. But since $W\in\calw_r$ is by definition a component of $X_r$, it is a connected saturated subset of $K$, and for each simplex contained in $W$ we have $\localrank(\clo(\Theta(\sigma)))=r$. Hence, again using the assumption that $\shortone(p)<\fone(\lambda)$ for every $p\in M$, we may apply Lemma \ref{commies} to deduce 
that $\Theta(W)$ is contained in a subgroup $C$ of $\Gamma$ such that  $\localrank(C)=r$. 
In particular we have $\localrank(C)<4$, and since $\Gamma$ is $4$-free it follows that $C$ is strongly closable in $\Gamma$. By Assertion (1) of Proposition \ref{closure properties}, $\clo(C)$ is a subgroup of $\Gamma$ having local rank at most $3$. We may now apply  Proposition \ref{and the normalizer}, with $\Theta(W)$
and $C$ playing the respective roles of $B$ and $A$, to deduce that the normalizer of $\Theta(W)$ in $\Gamma$ is a subgroup of $\clo(C)\le\Gamma$.  (The non-triviality of $\Theta(W)$, which is required for the application of Proposition \ref{and the normalizer}, holds because $W\neq\emptyset$; see \ref{fitz}.)
In particular we have $\Gamma_s\le\clo(C)$.

Since $\Gamma$ is $4$-free and $\localrank(\clo(C))\le3$, the group $\clo(C)$ is locally free. Hence $\Gamma_s\le\clo(C)$ is locally free.

We have shown that $\Gamma$ acts simplicially, without inversions, on the tree $T$, 
 and that the stabilizer of every vertex in $\Gamma$ is locally free. But since $\Gamma$ is cocompact, it is isomorphic to the fundamental group of a compact, orientable hyperbolic $3$-manifold.
According to \cite[Lemma 5.13]{fourfree}, the fundamental group of a closed, orientable, aspherical $3$-manifold cannot  act simplicially, without inversions, on a tree in such a way that each vertex stabilizer is locally free. This  contradiction completes the proof.
\EndProof
%fitz

\Corollary\label{corollary before key}
If
$M$ is any closed, orientable hyperbolic $3$-manifold such that
$\pi_1(M)$ is $4$-free, then there  is a point $p\in M$ with $\nextone(p)\ge\log7$.
\EndCorollary

\Proof
We apply 
Theorem 
\ref{before key}, taking $\lambda=\log7$. This
gives a point $p\in M$ such that either $\nextone(p)\ge\log7$ or
$\shortone(p) \ge\fone(\log7)=\log7$. Since $\nextone(p)\ge\shortone(p)$
by \ref{short and next}, the conclusion holds in both cases.
\EndProof

According to the discussion in \ref{what about the workers}, the conclusion of Corollary \ref{corollary before key} may be interpreted as saying that $M$ contains a $\log7$-semithick point. The corollary is therefore equivalent to Theorem 1.4 of \cite{fourfree}. Thus 
Theorem \ref{before key}
is a generalization of the latter result.

\begin{corollary}\label{key theorem}
Let  
$M$ be a closed, orientable hyperbolic $3$-manifold such that
$\pi_1(M)$ is $4$-free. Then  there  is a point $p_0\in M$ such that $\nextone(p_0)\ge\log7$ 
(so that in particular $\nextone(p_0)$ is in the domain of $\fone$)
and $\max_{p\in M}\shortone(p) \ge\fone(\nextone(p_0))$.
\end{corollary}

\Proof
According to Lemma \ref{continuous},  the functions $\shortone$ and $\nextone$ are continuous on the
closed manifold $M$. We set
$\nu_0=\max_{p\in M}\shortone(p)$
and
 $\lambda_0=\max_{p\in M}\nextone(p)$. Since by \ref{short and next}
 we have $\nextone(p)\ge\shortone(p)$ for every $p\in M$, we have
 $\lambda_0\ge\nu_0$. By Corollary \ref{corollary before key}, we have $\lambda_0\ge\log7$.

 We fix a point $p_0\in M$ with
$\nextone(p_0)=\lambda_0$. We will establish the conclusion of the
theorem with this choice of $p_0$. Equivalently, we need to  show that
 $\nu_0\ge\fone(\lambda_0)$.

If $\lambda$ is any number strictly greater than $\lambda_0$, then by
Theorem 
\ref{before key}, either there  is a point $p\in M$ with $\nextone(p)\ge\lambda$, or
there is a point  $p\in M$ with $\shortone(p) \ge\fone(\lambda)$. The
first of these alternatives contradicts the maximality of
$\lambda_0$. Hence there is a point  $p\in M$ with $\shortone(p)
\ge\fone(\lambda)$. In particular we have
$\nu_0\ge\fone(\lambda)$. Since this holds for every
$\lambda>\lambda_0$, and since $\fone$ is continuous by its
definition, we have $\nu_0\ge\fone(\lambda_0)$, as required.
\EndProof

\Corollary\label{key corollary}
Let  
$M$ be a closed, orientable hyperbolic $3$-manifold such that $\pi_1(M)$ is $4$-free. 
Let 
$\lambda^+$ and $\lambda^-$
be real numbers with $\log7
\le
\lambda^-\le\lambda^+
$. 
Then  at least one of the following alternatives holds:
\Alternatives
\item  there is a point $p\in M$ with $\nextone(p)=\lambda^+$;
\item  there is a point $p\in M$ with $\shortone(p)\ge\fone(\lambda^-)$;
\item  there  is a point $p_0\in M$ with $\lambda^-\le\nextone(p_0)\le\lambda^+$ and $\max_{p\in M}\shortone(p) \ge\fone(\nextone(p_0))$.
\EndAlternatives
\EndCorollary

\Proof
According to  
Corollary 
\ref{key theorem}, we may fix a point $p_0\in M$
such that  $\nextone(p_0)\ge\log7$ and $\max_{p\in M}\shortone(p) \ge\fone(\nextone(p_0))$.
Set $\lambda_0=\nextone(p_0)$. If
$\lambda^-\le\lambda_0\le\lambda^+$, then Alternative (iii) of the
present corollary holds. 

Next consider the case in which $\lambda_0<\lambda^-$. Since $\fone$
is strictly monotone decreasing, we have
$\fone(\lambda_0)>\fone(\lambda^-)$; hence $\max_{p\in
  M}\shortone(p)>\fone(\lambda^-)$, so that Alternative (ii) of the
corollary holds. 

We now turn to  the case in which $\lambda_0>\lambda^+$,
i.e. $\nextone(p_0)>\lambda^+$. 
Since the
function $\nextone$ is continuous on the connected manifold $M$ by Lemma \ref{continuous}, either
$\nextone$ takes the value $\lambda^+$ at some point of $M$, or the
range of $\nextone$ is entirely contained in the interval
$(\lambda^+,\infty)$. The first of these subcases is Alternative (i)
of the conclusion of the corollary. If we are in the second subcase,
then according to 
Proposition \ref{Margulis and s-one and s-two}
we have $\mu_M\ge\lambda^+$, i.e.
$\lambda^+$ is a Margulis
number for $M$. Hence by
Proposition \ref{my summer vocation} we have
$\Mthick(\mu)\ne\emptyset$; that is, there is a point 
$p_1\in
M$ 
such that 
$
\shortone(p_1)\ge\lambda^+
\ge
\log7$. 
But since $\fone$ is
monotone decreasing and $\fone(\log7)=\log7$, and since
$\lambda^-\ge\log7$, 
we have 
$\fone(\lambda^-)\le\log7$. 
Hence
$
\shortone(p_1)
\ge
\fone(\lambda^-)$, so that Alternative (ii)  holds in
this remaining subcase.
\EndProof

\Corollary\label{keyer corollary}
Let  
$M$ be a closed, orientable hyperbolic $3$-manifold such that $\pi_1(M)$ is $4$-free. 
Let 
 $\waslambdaminus$
be a real number with $\log7
\le
\waslambdaminus\le\log8
$. 
Then  at least one of the following alternatives holds:
\Alternatives
\item  there is a point $p\in M$ with $\shortone(p)\ge\fone(\waslambdaminus)$;
\item  there  is a point $p_0\in M$ with $\waslambdaminus\le\nextone(p_0)\le\log8$ and $\max_{p\in M}\shortone(p) \ge\fone(\nextone(p_0))$.
\EndAlternatives
\EndCorollary

\Proof
Note that the hypotheses of the present corollary imply those of Corollary \ref{key corollary} if we set 
$\lambda^-=\waslambdaminus$ and $\lambda^+ =\log8$. 
Hence one of the alternatives (i)--(iii) of the conclusion of Corollary \ref{key corollary} holds with 
these choices of $\lambda^-$ and $\lambda^+$. 
Alternatives (ii) and (iii) of Corollary \ref{key corollary} (with 
$\lambda^-=\waslambdaminus$ and $\lambda^+ =\log8$) 
respectively imply Alternatives (i) and (ii) of the present corollary. Now suppose that Alternative (i) of Corollary \ref{key corollary} holds with
$\lambda^-=\waslambdaminus$ and $\lambda^+ =\log8$, 
and fix a point $p_0\in M$ with $\nextone(p_0)=\log8$. 
We now apply
\cite[Corollary 9.3]{acs-singular} (a generalization of a result proved in \cite{accs}), which asserts that every closed, orientable hyperbolic $3$-manifold with a $3$-free fundamental group contains a hyperbolic ball of radius $(\log5)/2$. In the present situation, since $\pi_1(M)$ is $4$-free, it is in particular $3$-free, and hence
$\max_{p\in M}\shortone(p)\ge \log5$. But we have $\fone(\log8)=\log5$, and hence Alternative (ii) of the present corollary holds in this case.
\EndProof

\section{An observation about hyperbolic triangles}\label{triangle section}
The lemma established in this brief section will be needed in the proof of Lemma \ref{surgeon general} to give an upper bound for an angle of a hyperbolic triangle, given certain constraints on the sides of the triangle.

\NotationRemarks\label{all about Omega}
We define a function $\omega$ on $(0,\infty)^3$ by
$$\omega(x,y,z)=\frac{(\cosh x)(\cosh y)-\cosh z}{(\sinh x)(\sinh y)}.$$
We denote by $\frakN$ the subset of $\RR^3$ consisting of all points $(x,y,z)$ such that  $x$, $y$ and $z$ are all positive, and satisfy $x+y\ge z$, $x+z\ge y$ and $y+z\ge x$. For any $(x,y,z)\in\frakN$ we have $|\omega(x,y,z)|\le1$, and we define a function $\Omega$ on $\frakN$ by setting
$\Omega(x,y,z)=\arccos(\omega(x,y,z))\in[0,\pi]$ whenever $(x,y,z)\in\frakN$.

If a hyperbolic triangle has sides of lengths $x$, $y$ and $z$, then we have $(x,y,z)\in\frakN$, and according to the hyperbolic law of cosines,
$\Omega(x,y,z)$ is the angle between the sides of lengths $x$ and $y$. 

We shall extend $\Omega$ to a function $\barOmega$ on $(0,\infty)^3$ by setting
$$\barOmega(x,y,z)=\arccos(\min(\max(\omega(x,y,z),-1),1))\in[0,\pi]$$
for all $x,y,z>0$.

Note that if $x>0$ and $y>z>0$, the definition of $\omega$ implies that $\omega(x,y,z)>0$. Hence:
\Equation\label{music playing}
\barOmega(x,y,z)<\frac\pi2\text{ whenever } x>0\text{ and }y>z>0.
\EndEquation

We define a function $\theta$ on $(0,\infty)^2$ by setting
\Equation\label{threedle}
\theta(C,x)=
\barOmega\bigg(\arccosh\bigg(\frac{\cosh x}{\cosh C}\bigg), x, C\bigg)
\EndEquation
 if $x>C$, and
$\theta(C,x)=
\pi/2
$ if $x\le C$. 
(It is easily checked that $\theta$ is continuous, although this fact will not be used.)

We observe that if $x>C$, it follows from (\ref{music playing}) that the right-hand side of (\ref{threedle}) is less than $\pi/2$. Hence for all positive numbers $x$ and $C$, we have $0\le\theta(C,x)
\le
\pi/2$.

If $C$ is a positive number, and if $a_{1,0}$, $a_{1,1}$, $a_{2,0}$ and $a_{2,1}$ are numbers such that 
$0< a_{1,0}\le a_{1,1}$ and  $0< a_{2,0}\le a_{2,1}$, we set
$$A_1 (C,a_{1,0},a_{1,1},a_{2,0},a_{2,1})=\max_{(i, j)\in\{0,1\}\times\{0,1\}}\barOmega(a_{1,i},a_{2,j},C),$$
$$A_2 (C,a_{1,0},a_{1,1},a_{2,0},a_{2,1})=\max_{(m,i)\in\{1,2\}\times\{0,1\}}\theta(C,a_{m,i}),$$ 
and 
$$A (C,a_{1,0},a_{1,1},a_{2,0},a_{2,1})=\max(A_1 (C,a_{1,0},a_{1,1},a_{2,0},a_{2,1}),A_2 (C,a_{1,0},a_{1,1},a_{2,0},a_{2,1})).$$

\EndNotationRemarks

\Lemma\label{why o high o}
Suppose that for $m=1,2$, we are given positive numbers $a_{m,0}$, $a_{m,1}$ with $a_{m,0}\le a_{m,1}$. Let a positive number $C$ be given, and suppose that $(x_1^{(0)},x_2^{(0)})$ is a point of $[a_{1,0},a_{1,1}]\times[a_{2,0},a_{2,1}]$ such that
$(x_1^{(0)},x_2^{(0)},C)\in\frakN$, Then we have $\Omega(x_1,x_2,C)\le A (C,a_{1,0},a_{1,1},a_{2,0},a_{2,1})$.
\EndLemma

\Proof
Set $R=[a_{1,0},a_{1,1}]\times[a_{2,0},a_{2,1}]$.
The function $f$ defined on $(0,\infty)^2\supset R$ by
$f(x_1,x_2)=\omega(x_1,x_2,C)$ is continuous and therefore takes a least value on $R$; we choose a  point 
 $(\eta_1,\eta_2)\in R$ where this least value is achieved.

Since $(x_1^{(0)},x_2^{(0)},C)\in\frakN$, we have $f(x_1^{(0)},x_2^{(0)})=\omega(x_1^{(0)},x_2^{(0)},C)\le1$; since  
%$f(\eta_1,\eta_2)\le f(x_1^{(0)},x_2^{(0)})$, it follows that
$(x_1^{(0)},x_2^{(0)})\in R$, it follows that
$f(\eta_1,\eta_2)\le1$. This means that $\omega(\eta_1,\eta_2,C)\le1$, which with the definition of $\barOmega$ gives $\barOmega(\eta_1,\eta_2,C)=\arccos(\max(\omega(\eta_1,\eta_2,C),-1)$, i.\,e.
\Equation\label{freakin' stupid}
\barOmega(\eta_1,\eta_2,C)=\arccos(\max(f(\eta_1,\eta_2),-1)).
\EndEquation
 On the other hand, since 
$(x_1^{(0)},x_2^{(0)},C)\in\frakN$, we have 
$\omega(x_1^{(0)},x_2^{(0)},C)\ge-1$; 
%$f(x_1^{(0)},x_2^{(0)})=\omega(x_1^{(0)},x_2^{(0)},C)\ge-1$; 
%again using that 
and since $(x_1^{(0)},x_2^{(0)})\in R$, we have
$\omega(x_1^{(0)},x_2^{(0)},C) \ge f(\eta_1,\eta_2)$, Hence $\omega(x_1^{(0)},x_2^{(0)},C) \ge \max(f(\eta_1,\eta_2),-1)$. With (\ref{freakin' stupid}) and the definition of $\Omega$, this gives 
\Equation\label{I still say it's stupid}
\Omega(x_1^{(0)},x_2^{(0)},C) \le \barOmega(\eta_1,\eta_2,C).
\EndEquation

%\le f(x_1^{(0)},x_2^{(0)})$, we deduce that
%$f(x_1^{(0)},x_2^{(0)})\ge\max(f(\eta_1,\eta_2),-1)

%we have =   
We have $f(x_1,x_2)=(\coth x_1)(\coth x_2)-(\cosh C)(\cosech x_1)(\cosech x_2)$ for all $x_1,x_2>0$. Hence the partial derivative of $f$ with respect to the first variable is 
$$f_1'(x_1,x_2)=-(\cosech^2 x_1)(\coth x_2)+(\cosh C)(\cosech x_1)(\coth x_1)(\cosech x_2),$$
which vanishes precisely when $\cosh x_2=(\cosh C)(\cosh x_1)$. Likewise, we have $f_2'(x_1,x_2)=0$ precisely when $\cosh x_1=(\cosh C)(\cosh x_2)$. It follows that the two partial derivatives can vanish simultaneously only if $\cosh C=1$, which is not the case since $C>0$. Hence $f$ has no critical points in $R$, and $(\eta_1,\eta_2)$ must be a boundary point of $R$.

If $(\eta_1,\eta_2)$ is a corner point of the rectangle $R$, we have $(\eta_1,\eta_2)= (a_{1,i},a_{2,j})$ for some $(i, j)\in\{0,1\}\times\{0,1\}$. Using (\ref{I still say it's stupid}) and the definitions of $A_1 (C,a_{1,0},a_{1,1},a_{2,0},a_{2,1})$ and $A (C,a_{1,0},a_{1,1},a_{2,0},a_{2,1})$, we find that
$$
\begin{aligned}
\Omega(x_1^{(0)},x_2^{(0)},C) &\le \barOmega(\eta_1,\eta_2,C)=
%.{, and the definitions give
\barOmega(a_{1,i},a_{2,j},C)\le A_1  (C,a_{1,0},a_{1,1},a_{2,0},a_{2,1})\\
&\le A (C,a_{1,0},a_{1,1},a_{2,0},a_{2,1}).
\end{aligned}
$$

%the  definition of $A_1$ gives is the largest value that $\Omega(x_1,x_2,C)$ achieves at a corner of $R$. 

There remains the case in which $(\eta_1,\eta_2)$ is an interior point of a side of $R$. In this case we will show that 
%We shall complete the proof by showing that if (b) holds then
%$f$ achieves its minimum value at an interior point $(x_1,x_2)$ of a side of $R$ then 
$\barOmega(\eta_1,\eta_2,C)\le A_2\le A$, which with (\ref{I still say it's stupid}) implies the conclusion. By
%Since the definitions of $\Omega$, $f$ and $A_2$ are 
the symmetry in the definitions of $f$, $\Omega$ and $A_2$,
%in the variables  $x_1$ and $x_2$, 
we may 
assume that the side containing $(\eta_1,\eta_2)$ has the form $[a_{1,0},a_{1,1}]\times\{a_{2,i}\}$ for some $i\in\{0,1\}$. Thus we have $\eta_2=a_{2,i}$ and $f_1'(\eta_1,a_{2,i})=0$, so that $\cosh a_{2,i}=\cosh \eta_2=(\cosh C)(\cosh \eta_1)$. Thus we have $\cosh a_{2,i}>\cosh C$, and $\eta_1=\arccosh((\cosh a_{2,i})/(\cosh C))$. According to the definition of the function $\theta$, it follows that in this case we have
$$\barOmega(\eta_1,\eta_2,C)=\theta(C,a_{2,i})\le A_2$$
as required.
\EndProof

\section{Elementary facts about volumes}\label{elem vol section}

This section and the next two concern
ways of using  quantitative geometric information about a closed, orientable hyperbolic $3$-manifold $M$ 
% $M$ given by that proposition 
to obtain lower bounds for the volumes of certain subsets of $M$. 
In Section \ref{my favorite monster},
these methods will be used to make the transition from 
Corollary \ref{keyer corollary}
to a sufficient condition for a number to be a lower bound for $\vol M$  
when $\pi_1(M)$ is $k$-free for a given $k\ge4$; this will in  turn be used in later sections
to give the explicit lower volume bounds for $M$, described in the introduction, that hold when $M$ has sufficient topological complexity.

In Subsections \ref{bee of eggs}--\ref{reformulated surgeon general}, we obtain lower bounds for the volume of $\nbhd_{\nextone(p)}(p)$, where  $p$ is a point of $M$; these lower bounds depend on $\nextone(p)$ and on the minimal length of a loop based at $p$ representing   the generator of a short maximal cyclic subgroup of $\pi_1(M,p)$ (see \ref{short and next}). 
Lemma \ref{surgeon special} 
is similar in nature, but deal more specifically with the case in which $p$ lies on a suitably short closed geodesic. In Subsection \ref{Boroczky number} we obtain lower bounds for the volume of a certain metric neighborhood of a point $p\in M$ in terms of $\shortone(p)$, given that there is a  point of $M$ which is suitably distant from $p$. Subsections \ref{is that you wthree?} and \ref {new distant thunder} are devoted to establishing a lower bound for the {\it complement} of a certain metric neighborhood of a point $p$ of $M$, given that there is a $\mu$-thick point---where $\mu$ is a Margulis number for $M$---which is suitably distant from $p$. In Subsections \ref{stoppeth two}--\ref{stoppeth one}, some of these results are combined in ways that will prove useful in 
Sections \ref{my favorite monster}, \ref{4 and 5} and \ref{drilling section}.

A significant part of this section consists of  review of material from \cite{fourfree}, but  much of this material will be reorganized here,
and a number of the results will be improved upon.

\Notation\label{bee of eggs}
We define  a 
strictly  increasing
function $B(x)=\pi(\sinh(2x)-2x)$ for $x>0$. Geometrically, $B(x)$ gives the volume of a ball  in $\HH^3$ of
radius $x$. 
\EndNotation

\Number\label{what to call it}
We shall review some material
from Subsections 6.1 and 6.5 of \cite{fourfree}.
%Let $R$ be a positive real number,  let 
%If $N$ denote an open ball of
%radius $R$ in $\HH^3$, and let $S$ denote the sphere $\partial\overline N$.
If $N$ is an open ball 
in $\HH^3$, 
%and let $S$ denote the sphere $\partial\overline N$.
%Let z0 be a point of H3. For each positive real number R we shall
%denote by S(R,z0) the sphere of radius R centered at z0. Thus S(R,z0)
%is the boundary of the closed ball N(z0, R)= ̇ N(z0, R) of radius R
%about z0. 
then for each point $\zeta$ in the sphere $\partial\overline N$, we will  denote by   $\eta_{N,\zeta}$ 
the ray
originating at the center of $N$  and passing through $\zeta$. 
For each point $\zeta\in \partial\overline N$ and each number $w\ge 0$ we shall denote
by $\Pi_N(\zeta,w)$ the plane which meets 
$\eta_{N,\zeta}$
perpendicularly at a distance $w$  from the center of $N$, and by
$H_N(\zeta, w)$ the closed half-space which is bounded by $\Pi_N(\zeta,w)$
and has unbounded intersection with 
$\eta_{N,\zeta}$.
For any given $w>0$ and
for an arbitrary point
$\zeta\in 
\partial\overline N
$, we denote by $K_N(\zeta , w)$
%$\kappa(R,w)$ the
%volume of
 the ``cap'': $N \cap H(\zeta, w)$. When the choice of the ball $N$ is understood, we will write $\eta_\zeta$, $\Pi(\zeta,w)$, $H(\zeta,w)$ and $K(\zeta,w)$ in place of  $\eta_{N,\zeta}$, $\Pi_N(\zeta,w)$, $H_N(\zeta,w)$ and $K_N(\zeta,w)$ respectively.

Now let $R$ be a positive real number,  and let 
 $N$ denote an open ball of
radius $R$ in $\HH^3$.
We set $\kappa(R,w)=\vol K(\zeta,w)$ for an arbitrary point $\zeta\in \partial\overline N$. We
%the
%volume of
% the ``cap'': $N \cap H(\zeta, w)$
%shall
 endow $\partial\overline N$ with the spherical metric in which the distance
between two points $\zeta,\zeta'\in S$ is the angle between
$\eta_\zeta$ and $\eta_{\zeta'}$, and
we denote by $\iota(R,w,w',\alpha)$ and
$\sigma(R,w,w',\alpha)$ the respective volumes of $K(\zeta, w) \cap
K(\zeta', w')$ and $K(\zeta, w) \cup
K(\zeta', w')$, where $w,w'>0$ and $\alpha\in[0,\pi]$ are given, and
$\zeta$  and $\zeta'$ are  points in $S$ whose spherical  distance is
$\alpha$. Then $\kappa$ is a well-defined function of two positive
real variables, 
monotone increasing in its first argument and monotone decreasing in its second;
while $\iota$ and $\sigma$ are well-defined functions of three positive
real variables and a fourth variable whose values are restricted to
$[0,\pi]$. 
(If $w\ge R>0$, then $\kappa(R,w)$ vanishes, as does $\iota(R,w,w',\alpha)$ for any $w'>0$ and any $\alpha\in[0,\pi]$.)

 Note that we have 
\Equation\label{six of one}
\sigma(R, w, w', \alpha) = \kappa(R, w)
+ \kappa(R, w') - \iota(R, w, w', \alpha)
\EndEquation
 for any positive numbers $R,w,w'$
and any $\alpha\in[0,\pi]$. Analytic expressions 
for
the functions
$\kappa$ and $\iota$ (and hence 
for
$\sigma$) are given in \cite[Section 14]{fourfree}.

According to \cite[Proposition 6.7]{fourfree}, $\sigma$ is monotone decreasing in its third argument and monotone increasing in its fourth argument. Furthermore, the first paragraph of the proof of \cite[Proposition 6.7]{fourfree}, with union signs replaced by intersection signs, shows that $\iota$ is monotone decreasing in its third argument.  Since $\sigma$ and $\iota$ are symmetric in their second and third arguments (cf. \cite[6.5]{fourfree}), they are also monotone decreasing in their second argument. In view of (\ref{six of one}), the fact that $\sigma$  is monotone increasing in its fourth argument may be  interpreted as meaning that $\iota$ is monotone decreasing in its fourth argument. 

\EndNumber

Here is the first of several results that we will extract from Sections 6 and 7 of \cite{fourfree} and apply later in this paper:

\Proposition\label{soft caps}
Let $M$ be a closed, orientable hyperbolic $3$-manifold, and
write $M=\HH^3/\Gamma$
where $\Gamma\le\Isom_+(\HH^3)$ is discrete and torsion-free.
Let $\alpha$ be a positive 
number, and let $p$ be a point of $M$ with $\nextone(p )\ge\alpha$. 
Let $P$ be a point of $\HH^3$ which projects to $p$ under the quotient map, let $j:\pi_1(M,p)\to\Gamma$ denote the isomorphism determined by the compatible base points $P\in\HH^3$ and $p\in M$, let $x$ denote the image under $j$ of a generator of a short maximal cyclic subgroup (\ref{short and next}) of  $\pi_1(M,p)$. Let $N\subset\HH^3$ denote the open ball of radius $\alpha/2$ centered at $P$. 
%Let $F\subset\ZZ$ denote the finite set consisting of all non-zero integers $n$ such that $d(x^n,P)<\alpha$ (in the notation of \ref{Q-def}). 
For each  integer $n\ne0$, set $d_n=
d(x^n,P)$ (in the notation of \ref{Q-def}). 
Let $\zeta_n$ denote the point of intersection of $\partial\overline N$  with the ray originating at $P$ and passing through $x^n\cdot P$. 
%Suppose that $g$ is the generator of a short maximal cyclic subgroup (see \ref{short and next}) of %$\pi_1(M,p)$. For each $n\in\ZZ$, let $d_n$ denote the length of the shortest loop based at $p$ %representing $g^n$, and  let $F$ denote the finite set $\{n\in\ZZ-\{0\}:d_n\le\alpha\}$. 
Then 
%there is a family $(\zeta_n)_{n\in F}$ of points in $\partial\overline{N}$ such that 
in the notation of \ref{what to call it} we have
\Equation\label{the wood}
\vol\nbhd\nolimits_{\alpha/2}(p)
=B(\alpha/2)-\vol\bigg(\bigcup_{0\ne n\in \ZZ}K\bigg(\zeta_n,\frac{d_n}2\bigg)\bigg).
\EndEquation
\EndProposition

\Proof
First consider the case in which $\pi_1(M,p)$ has only one short maximal cyclic subgroup. As we pointed out in \ref{what about the workers}, this says, in the notation of \cite{fourfree}, that $p\in\frakG_M$, and we then have $\fraks_M(p)=\nextone(p)\ge\alpha$, where $\fraks_M$ is defined as in \cite{fourfree}; furthermore, the short maximal cyclic subgroup of $\pi_1(M,p)$ is denoted $C_p$ in the notation of \cite{fourfree}. In this  case, the statement of the present proposition is precisely that of \cite[Proposition 6.2]{fourfree}, with $p$, $P$ and $\alpha$ playing the respective roles of $P$, $\tP$ and $\lambda$ in the latter result.

In the case where $\pi_1(M,p)$ has more than one short maximal cyclic subgroup, we have by definition (see \ref{short and next}) that $\nextone(p)=\shortone(p)$; hence $\shortone(p)\ge\alpha$, so that $p$ is the center of a hyperbolic ball of radius $\alpha/2$. Hence $
\vol\nbhd\nolimits_{\alpha/2}(p)
=B(\alpha/2)$. But in this case, for each integer $n\ne0$, the definition of $\shortone(p)$ gives $d_n\ge\shortone(p)\ge\alpha$ so that $K(\zeta_n,d_n/2)=\emptyset$. Hence the left-hand side of (\ref{the wood}) is also equal to $B(\alpha/2)$.
\EndProof

\Number\label{audrey junior}
We shall review  the definitions and basic properties of some more functions that are introduced in \cite{fourfree}.

As in \cite[Subsection 7.1]{fourfree},  we define, for each integer $n\ge1$, a function $\Phi_n$ on the domain $\{(\delta, D) : 0 < \delta\le D\}\subset\RR^2$
by
$$\Phi_n(\delta,D)=\arccosh\bigg(\cosh (n\delta)+\frac{(\cosh (n\delta)-1)(\cosh D-\cosh\delta)}{\cosh\delta+1}\bigg).$$
Note that $\Phi_n$ is monotonically increasing in its second argument.

The significance of the function $\Phi_n$ arises from
Lemma 7.3 of \cite{fourfree}, which asserts (in the notation of \ref{Q-def} above) that if  $P$ is a point of $\HH^3$, and $x$ is a loxodromic isometry of $H^3$ whose translation length is bounded below by  a given positive number $\delta$ then for every positive integer $n$ we have
\Equation\label{wi-fi}
d(x^n,P)\ge\Phi_n(\delta, d(x,P)).
\EndEquation
(The statement of Lemma 7.3 given in \cite{fourfree} unfortunately contained a typographical error: the inequality (\ref{wi-fi}) was  reversed. The form given here matches both the proof of Lemma 7.3 given in \cite{fourfree}, and the applications given there.)

We will also use Lemma 7.4 of \cite{fourfree}, which asserts that if $n$ is a positive integer, and if $\delta$ and $D$ are real numbers with $0<\delta\le D$, then 
\Equation\label{appomatox}
n\delta\le\Phi_n(\delta,D)\le nD.
\EndEquation

\EndNumber
%\nu

\Number\label{moore}
As in \cite[Subsection 7.5]{fourfree}, we define a function $\Psi$ on the domain
$\{(x,y) : 0 < y \le 2x\}\subset\RR^2$,  with values in $[0, \pi/2]$, by
$$\Psi(x, y) = \arccos((\coth x)(\coth y - \cosech y)).$$
Thus in the notation of \ref{all about Omega}, we have $\Psi(x,y)=\Omega(y,x,x)$. Hence if an isosceles hyperbolic triangle has base $y$ and has its other two sides equal to $x$, we have $y\le2x$ (i.e. $(y,x,x)\in \frakN$ in the notation of \ref{all about Omega}), and the base angles of the triangle are equal to $\Psi(x, y)$.

It is pointed out in \cite[Subsection 7.5]{fourfree} that $\Psi$ is monotone increasing in its first argument and monotone decreasing in its second.

As in \cite[Subsection 6.3]{fourfree}, we denote by $\Theta$ the
real-valued function with domain $\{(w, R) : 0 < w < R\}\subset \RR^2$
which is
       defined by
$$\Theta(w,R)=\arccos\bigg(\frac{\tanh w}{\tanh R}\bigg)$$
and takes values in $(0, \pi/2)$. 
Note that $\Theta$ is monotonically 
decreasing
in its first argument and monotonically 
increasing
in its second argument.

\EndNumber
%w_

\Notation\label{whats Lambda}
We will need one function that is not defined in \cite{fourfree}. Recall from \ref{all about Omega} that
$A (C,a_{1,0},a_{1,1},a_{2,0},a_{2,1})$ is defined whenever
$0\le a_{1,0}\le a_{1,1}$ and  $0\le a_{2,0}\le a_{2,1}$. If $D$ and 
$\delta$ are  numbers with $0\le\delta\le D$,
then for $n=2,3$ we have 
$n\delta\le\Phi_n(\delta,D)$ by (\ref {appomatox}). Hence  we may define a function $\Lambda$ on 
$\{(\delta,D):0\le\delta\le D\}\subset\RR^2\}$
by
$$\Lambda(\delta,D)=A( D,2\delta,\Phi_2(\delta,D),3\delta, \Phi_3(\delta,D)).$$
\EndNotation

%\lambda
\Lemma\label{surgeon general}
Let $p$ be a
point of  a closed, orientable hyperbolic 3-manifold $M$,
and let $C$ be a short maximal cyclic subgroup (see \ref{short and next}) of $\pi_1(M,p)$.
Let $\delta$ and $\alpha$ be constants with $0 < \alpha < 4\delta$. Assume that
\Bullets
\item $\nextone(p) \ge\alpha$, and that
\item
the conjugacy class of a generator of $C$ is represented by a closed geodesic in $M$ having length
  at least $\delta$. 
\EndBullets
Let $D$ denote the minimal length of a loop based at $p$ that represents a generator of $C\le \pi_1(M,p)$.
Then 
$\delta\le D$, 
so that
$T_n := \Phi_n(\delta, D)$ is defined for every $n \ge 1$. Furthermore, 
%lies
%in the domain of $\Psi$,  
%and 
we have
$$\vol \nbhd\nolimits_{\alpha/2}(p)\ge B\bigg(\frac \alpha2\bigg)- 2\sigma\bigg(\frac\alpha2 , \frac D2 , \frac {T_2}2 , \Psi(D, T_2)\bigg) - 2\kappa\bigg(\frac\alpha2 , \frac{T_3}2\bigg)+2\iota\bigg(\frac\alpha2,\frac{T_2}2,\frac{T_3}2,\Lambda(\delta,D)\bigg),$$
%$$\vol \nbhd\nolimits_{\alpha/2}(p)\ge B\bigg(\frac \alpha2\bigg)- 2\sigma\bigg(\frac\alpha2 , \frac D2 , \frac {T_2}2 , \Psi(D, T_2)\bigg) - 2\kappa\bigg(\frac\alpha2 , \frac{T_3}2\bigg)$$
(where
$\Psi(D, T_2)$ is defined because $T_2\le 2D$ by 
\ref{wi-fi}).
If in addition we have $D < T_3 < \alpha$, so that, in particular,
the quantities $\Theta(D/2, \alpha/2)$ and $\Theta(T_3/2 ,
\alpha/2)$ 
are defined 
(see \ref{audrey junior}), 
and if
$$\cos\bigg(\Theta\bigg(\frac D2,\frac\alpha2\bigg)-\Theta\bigg(\frac{T_3}2,\frac\alpha2\bigg)\bigg)<
\frac{\cosh D\cosh T_3-\cosh2D}{\sinh D\sinh T_3},$$
     then
$$ \vol \nbhd\nolimits_{\alpha/2}( p)\ge B\bigg(\frac\alpha2\bigg) -2\sigma\bigg(\frac\alpha2 , \frac D2 , \frac {T_2}2 , \Psi(D, T_2)\bigg).$$
\EndLemma

\begin{specialremark}\label{ok, dude}
\textnormal{The quantity that is denoted $D$ in the statement of Lemma \ref{surgeon general} is not necessarily the same as  $\shortone(p)$. Whereas $D$ is the minimal length of a loop based at $p$ that represents a generator $x$ of  $C\le \pi_1(M,p)$, the definitions in \ref{short and next} imply that $\shortone(p)$ is the minimal length of a loop based at $p$ that represents any non-trivial power of $x$. It is a standard observation that these need not be equal: for example, if we write $M=\HH^3/\Gamma$, where $\Gamma\le\isomplus(\HH^3)$, if $P$ is a point of $\HH^3$ which projects to $p$ under the quotient map, if $j:\pi_1(M,p)\to\Gamma$ denotes the isomorphism determined by the compatible base points $P\in\HH^3$ and $p\in M$, and if $j(x)$ is a loxodromic isometry having very small translation length and twist angle very close to $\pi$, then the minimal length of a loop representing  $x^2$ will be less than the minimal length of a loop representing $x$.}

\textnormal{Because $D$ and $\shortone(p)$ may be distinct, one cannot rule out the possibility that $D>\nextone(p)$. (This is related to the issue addressed in Lemma \ref{other woney} below.)}
\end{specialremark}

\Proof[Proof of Lemma \ref{surgeon general}]
Let $P$ be a point of $\HH^3$ which projects to $p$ under the quotient map, let $j:\pi_1(M,p)\to\Gamma$ denote the isomorphism determined by the compatible base points $P\in\HH^3$ and $p\in M$, and let $x$ denote the image under $j$ of a generator of $C$. Let $N$ denote the ball of radius $\lambda/2$ centered at $P$, and for each integer $n\ne0$ define the quantity $d_n$ ,  and the point $\zeta_n\in\partial\overline N$, as in Proposition \ref{soft caps}. Note that for any $n\ne0$ we have $d_{-n}=d_n$. The hypothesis of the present lemma implies that $x$ has translation length at least $\delta$, so that in particular $D\ge\delta$, and thus each $T_n$ is defined. According to (\ref{wi-fi}) we have 
\Equation\label{cuz}
d_n=d_{|n|}\ge T_{|n|}
\EndEquation
for each $n\ne0$.
 Let us set $K_n=K(\zeta_n,d_n/2)$ for each integer $n\ne0$. 
For  $0\ne n\in\ZZ$ we have $d_n
\ge
|n|\delta$; since $\lambda<4\delta$, it follows that 
\Equation\label{pontex}
K_n=\emptyset\text{ for }|n|>3.
\EndEquation

Let us set $L_n=K(\zeta_n,T_{|n|}/2)$ for each $n\ne0$. It follows from (\ref{cuz}) that 
\Equation\label{caps too}
K_n\subset L_n\text{ for each }n\ne0.
\EndEquation

Combining  Proposition \ref{soft caps} with (\ref{pontex}) and (\ref{caps too}), we obtain
\Equation\label{unmaligned}
\begin{aligned}
\vol\nbhd\nolimits_{\lambda/2}(p)
&=B\bigg(\frac\lambda2\bigg)-\vol\bigg(\bigcup_{0\ne n\in\ZZ}K_n\bigg)\\
&=B\bigg(\frac\lambda2\bigg)-\vol\bigg(\bigcup_{0<|n|\le3}K_n\bigg)\\
&\ge B\bigg(\frac\lambda2\bigg)-\vol\bigg(\bigcup_{0<|n|\le3}L_n\bigg)\\
%\bigg(\bigcup_{0\ne n\in \ZZ}K_n\bigg)\\
%&=B\bigg(\frac\lambda2\bigg)-\vol(( K_1\cup K_2\cup K_3)\cup ( K_{-1}\cup K_{-2}\cup K_{-3}))\\
&\ge B\bigg(\frac\lambda2\bigg)-(\vol( L_1\cup L_2\cup L_3)+\vol(  L_{-1}\cup L_{-2}\cup L_{-3})).
\end{aligned}
\EndEquation

Let $\epsilon\in\{-1,1\}$ be given.
Consider the triangle with vertices $P$, $x^\epsilon\cdot P$ and $x^{2\epsilon}\cdot P$. The side joining $P$ to $x^{2\epsilon}\cdot P$ has length $d_2$, and each of the other sides has length $D$. It therefore follows from  \ref{moore} that the angle of the triangle at $P$ is $\Psi(D,d_2)$. Hence $\Psi(D,d_2)$ is the spherical distance between $\zeta_\epsilon$ and $\zeta_{2\epsilon}$. The  definition of the function $\sigma$ given in \ref{what to call it} then implies that $\vol(L_\epsilon\cup L_{2\epsilon})=\sigma(\lambda/2,D,T_2/2,
\Psi(D,d_2))$.  But we have $d_2\ge T_2$ by (\ref{cuz}), and since the function $\Psi$ is monotone decreasing in its second argument by \ref{moore}, we have $\Psi(D,d_2)\le\Psi(D,T_2)$. Since, according to 
\cite[Proposition 6.7]{fourfree} (cf.  \ref{what to call it}), the function $\sigma$ is monotone increasing in its  fourth argument, we deduce that
\Equation\label{ole man hudson}
\vol(L_\epsilon\cup L_{2\epsilon})\le\sigma\bigg(\frac\lambda2,\frac D2,\frac{T_2}2, \Psi(D,T_2)\bigg).
\EndEquation

By definition we have $
\vol L_{3\epsilon}
=\kappa(\lambda/2,T_{3}/2)$, so that
\Equation\label{toot-toot}
\vol(L_\epsilon\cup L_{2\epsilon}\cup L_{3\epsilon})=\vol(L_\epsilon\cup L_{2\epsilon})+\kappa\bigg(\frac\lambda2,\frac{T_3}2\bigg)-\vol((L_\epsilon\cup L_{2\epsilon})\cap L_{3\epsilon}).
\EndEquation

Now 
consider
the triangle with vertices 
$P$, $x^{3\epsilon}\cdot P$ and $x^{2\epsilon}\cdot P$;
 the sides opposite these three vertices have respective lengths $D$, 
$d_2$ and $d_3$.
It therefore follows from \ref{all about Omega} that we have 
$ (d_2,d_3,D)\in\frakN$,
and that the angle of the triangle at $P$ is 
$\Omega (d_2,d_3,D)$. Hence $\Omega (d_2,d_3,D)$
is the spherical distance between $\zeta_{2\epsilon}$ and $\zeta_{3\epsilon}$. In view of the definition of $\iota$ (see \ref{what to call it}), it follows that
\Equation\label{epsigh-lon}
%\begin{aligned}
 \vol((L_\epsilon\cup L_{2\epsilon})\cap L_{3\epsilon})\ge \vol(L_{2\epsilon}\cap L_{3\epsilon})=\iota\bigg(\frac\lambda2,\frac{T_2}2,\frac{T_3}2, 
\Omega (d_2,d_3,D)\bigg).
%\\
%&=\iota(\lambda/2,d_{2}/2,d_{3}/2, \Omega (d_3,d_2,D)).
%\end{aligned}
\EndEquation
From (\ref{ole man hudson}), (\ref{toot-toot}) and (\ref{epsigh-lon}) it follows that
\Equation\label{snuff it}
\vol(L_\epsilon\cup L_{2\epsilon}\cup L_{3\epsilon})\le\sigma\bigg(\frac\lambda2,\frac D2,\frac{T_2}2, \Psi(D,T_2)\bigg) +\kappa\bigg(\frac\lambda2,\frac{T_3}2\bigg)-
\iota\bigg(\frac\lambda2,\frac{T_2}2,\frac {T_3}2, \Omega (d_2,d_3,D)\bigg).
\EndEquation

Now by \ref{appomatox} we have $n\delta\le d_n\le T_n$ for $n=2,3$. We may therefore apply Lemma \ref{why o high o}, letting $D$, 
$2\delta$, $T_2$, $3\delta$ and $T_3$
play the respective roles of $C$, $a_{1,0}$, $a_{1,1}$, $a_{2,0}$ and $a_{2,1}$, and taking 
$x_1^{(0)}=d_2$ and $x_2^{(0)}=d_3$,
to deduce that 
$
\Omega (d_2,d_3,D)
\le A (C,2\delta,T_2, 3\delta,T_3)=A (C,2\delta,\Phi_2(\delta,D), 3\delta, \Phi_3(\delta,D))$. In view of the definition given in \ref{whats Lambda}, this means that
$\Omega (d_2,d_3,D)\le\Lambda(\delta,D)$. 
But $\iota$ is monotone decreasing in its fourth argument; indeed, it was pointed out in  \ref{what to call it} that this is included in \cite[Proposition 6.7]{fourfree}. It
now follows that 
$\iota(\lambda/2,T_{2}/2,T_{3}/2, \Omega (d_2,d_3,D))\ge \iota(\lambda/2,T_{2}/2,T_{3}/2, \Lambda(\delta,D))$.
Combining this inequality with (\ref{snuff it}), we deduce that
\Equation\label{cuckoo}
\vol(L_\epsilon\cup L_{2\epsilon}\cup L_{3\epsilon})\le\sigma\bigg(\frac\lambda2,\frac D2,\frac{T_2}2, \Psi(D,T_2)\bigg) +\kappa\bigg(\frac\lambda2,\frac{T_3}2\bigg)-
\iota\bigg(\frac\lambda2,\frac{T_2}2,\frac{T_3}2, \Lambda(\delta,D)\bigg).
\EndEquation

We may now combine (\ref{unmaligned}) with
% (\ref{pontex}) and
 the cases  $\epsilon=1$ 
and
$\epsilon=-1$ of (\ref{cuckoo}) to obtain 
$$
\vol\nbhd\nolimits_{\lambda/2}(p)
%&=B\bigg(\frac\lambda2\bigg)-\vol\bigg(\bigcup_{0\ne n\in \ZZ}K_n\bigg)\\
%&=B\bigg(\frac\lambda2\bigg)-\vol(( K_1\cup K_2\cup K_3)\cup ( K_{-1}\cup K_{-2}\cup K_{-3}))\\
%&\ge B\bigg(\frac\lambda2\bigg)-(\vol( K_1\cup K_2\cup K_3)+\vol(  K_{-1}\cup K_{-2}\cup K_{-3}))\\
%&
\ge B\bigg(\frac\lambda2\bigg)-2\bigg(\sigma\bigg(\frac\lambda2,\frac D2,\frac{T_2}2, \Psi(D,T_2)\bigg) +\kappa\bigg(\frac\lambda2,\frac{T_3}2\bigg)-
\iota\bigg(\frac\lambda2,\frac{T_{2}}2,\frac{T_{3}}2, \Lambda(\delta,D)\bigg)\bigg),
%\end{aligned}
$$
which  proves the first assertion of the lemma.
%\kappa(\lambda/2,d_3) K
%\le\sigma(\lambda/2,D,d_2/2, \Psi(D,T_2)).
%\Psi(D,d_2)

To prove the second assertion, we first consider the case in which $\pi_1(M,p)$ has only one short maximal cyclic subgroup. As we pointed out in \ref{what about the workers}, this says, in the notation of \cite{fourfree}, that $p\in\frakG_M$, and we then have $\fraks_M(p)=\nextone(p)\ge\alpha$, where $\fraks_M$ is defined as in \cite{fourfree}; furthermore, the short maximal cyclic subgroup of $\pi_1(M,p)$ is denoted $C_p$ in the notation of \cite{fourfree}. Furthermore, in this case the quantity denoted by $D$ in the statement of the present proposition is the  minimal length of a loop based at $p$ that represents a generator of the unique maximal cyclic subgroup of $\pi_1(M,p)$; in the notation of \cite{fourfree}, this quantity is denoted by $D_M(p)$. The statement of the second assertion of the present proposition is then seen to be precisely that of the second assertion of \cite[Lemma 7.6]{fourfree}, with $p$ and $\alpha$ playing the respective roles of $P$ and $\lambda$ in the latter result. 

In the case where $\pi_1(M,p)$ has more than one short maximal cyclic subgroup, we have by definition (see \ref{short and next}) that $\nextone(p)=\shortone(p)$; hence $\shortone(p)\ge\alpha$, so that $p$ is the center of a hyperbolic ball of radius $\alpha/2$. Hence $
\vol\nbhd\nolimits_{\alpha/2}(p)
=B(\alpha/2)$. This is enough to imply the second assertion in this case, and the proof of the lemma is thus complete.

(It may be noted that the argument that has been used to deduce the second assertion of the present lemma from the second assertion of \cite[Lemma 7.6]{fourfree} would also allow us to deduce a weaker version of the first assertion of the present lemma from the first assertion of \cite[Lemma 7.6]{fourfree}. The stronger version of the first assertion of the present lemma was made possible by the use of
the material in Section \ref{triangle section} above.)
\EndProof

\Reformulation\label{reformulated surgeon general}
For applications of Lemma \ref{surgeon general}, it will be convenient to define subsets $\frakW $, $\frakW '$ and $\frakW ''$ of $\RR^3$, and  functions $\W$ and $\Wone$ with domain $\frakW $, as follows. 
We set
$\frakW =\{(\alpha,\delta,D)\in\RR^3: 0<\delta\le D \text{ and } \alpha>0\}$. 
Note that if $(\alpha,\delta,D)\in \frakW $, then $\Phi_n(\delta, D)$ is defined for every $n \ge 1$. Furthermore, 
we have
$\Phi_2(\delta, D) \le 2D$ by 
\ref{appomatox},
so that
$\Psi(D, 
\Phi_2(\delta, D))$ is defined. For each 
$(\alpha,\delta,D)\in \frakW $ we set
$$\W(\alpha,\delta,D)=B(\alpha/2)-2\sigma(\alpha/2,D/2,\Phi_2(\delta,D)/2,\Psi(D, \Phi_2(\delta,D))).$$
Next we define $\frakW '=\{(\alpha,\delta,D)\in \frakW : 
D<\Phi_3(\delta,D)<\alpha\}
$, and observe that by 
\ref{moore},
the quantities $\Theta(D/2, \alpha/2)$ and $\Theta(\Phi_3(\delta,D)/2 ,
\alpha/2)$ 
are defined for every $(\alpha,\delta,D)\in \frakW '$. We 
define
$\frakW ''$ to be the set of all $(\alpha,\delta,D)\in \frakW '$ such that
\Equation\label{erev birthday}
\cos\bigg(\Theta\bigg(\frac D2,\frac\alpha2\bigg)-\Theta\bigg(\frac{\Phi_3(\delta,D)}2,\frac\alpha2\bigg)\bigg)<
\frac{\cosh D\cosh \Phi_3(\delta,D)-\cosh2D}{\sinh D\sinh \Phi_3(\delta,D)}.
\EndEquation
We then define the function $\Wone$ on $\frakW $ by setting
$$\Wone(\alpha,\delta,D)=\W(\alpha,\delta,D)$$
if 
$(\alpha,\delta,D)\in\frakW ''$, and
$$
\Wone(\alpha,\delta,D)=
\W(\alpha,\delta,D)-2\kappa\bigg(\frac\alpha2,\frac{\Phi_3(\delta,D)}2 \bigg)+2\iota\bigg(\frac\alpha2,\frac{\Phi_2(\delta,D)}2,\frac{\Phi_3(\delta,D)}2,\Lambda(\delta,D)\bigg)
$$
if 
$(\alpha,\delta,D)\in\frakW -\frakW ''$.

In terms of these definitions, we may reformulate Lemma \ref{surgeon general} as follows. Let $p$ be a
point of a closed, orientable hyperbolic 3-manifold $M$,
and let $C$ be a short maximal cyclic subgroup  of $\pi_1(M,p)$.
Let $\delta$ and $\alpha$ be constants with $0 < \alpha < 4\delta$. Assume that
$\nextone(p) \ge\alpha$, and that
the conjugacy class of a generator of $C$ is represented by a closed geodesic in $M$ having length
  at least $\delta$. 
Let $D$ denote the minimal length of a loop based at $p$ that represents a generator of $C$. Then 
$(\alpha,\delta,D)\in\frakW $, and $\vol \nbhd_{\alpha/2}( p)\ge\Wone(\alpha,\delta,D)$.
(The superscript ``near'' and the subscript ``ST'' are meant to indicate that $\Wone$ gives a lower bound for the volume of a suitable neighborhood of a point which is $\alpha$-semithick (see  \ref{what about the workers}).)
\EndReformulation
%\frakW \lambda

\Remark\label{cases and cases}
For any $(\alpha,\delta,D)\in\frakW $, it follows from the geometric definitions of $\kappa$ and $\iota$ (see \ref{what to call it}) that 
$$\kappa\bigg(\frac\alpha2,\frac{\Phi_3(\delta,D)}2 \bigg)
\ge
\iota\bigg(\frac\alpha2,\frac{\Phi_2(\delta,D)}2,\frac{\Phi_3(\delta,D)}2,\Lambda(\delta,D)\bigg).$$
In view of the definition of $\Wone$, it then follows that for any $(\alpha,\delta,D)\in\frakW $  we have
\Equation\label{inequality version}
\Wone(\alpha,\delta,D)\ge
\W(\alpha,\delta,D)-2\kappa\bigg(\frac\alpha2,\frac{\Phi_3(\delta,D)}2 \bigg)+2\iota\bigg(\frac\alpha2,\frac{\Phi_2(\delta,D)}2,\frac{\Phi_3(\delta,D)}2,\Lambda(\delta,D)\bigg).
\EndEquation
\EndRemark

We shall extract one more result, Lemma \ref{surgeon special} below, from \cite[Section 7]{fourfree}.

\Notation\label{short geo def}
We define a function $\Wtwo$ on the quadrant $(0,\infty)^2\subset\RR^2$ by $\Wtwo(\lambda,l)=B(\lambda/2)-2\kappa(\lambda/2,l/2)$.

(The superscript ``near'' and the subscript ``SG'' are meant to indicate that $\Wtwo$ gives a lower bound for the volume of a suitable neighborhood of a point which lies on a short geodesic in a hyperbolic manifold;
this is the content of Lemma \ref{surgeon special} below.)
\EndNotation

\Lemma\label{surgeon special}
Let $M$ be a closed, orientable hyperbolic $3$-manifold, and let $\mu$ be a Margulis number for $M$. Suppose that $c$ is a closed geodesic in $M$ of length $l < \mu$, and let $p$ be any point of 
$|c|$ (in the notation of \ref{nbhd}). 
Then  $\shortone(p) = l$, and
$$\vol \nbhd \nolimits_{\nextone(p )/2} (p) = \Wtwo(\nextone (p ),l).$$
\EndLemma

\Proof
Since $c$ is a closed geodesic of length $l$, there is an element $g$ of $\pi_1(M,p)$, generating a maximal cyclic subgroup $C$, such that for every non-zero integer $n$, the shortest loop based at $p$ and representing $g^n$ has length $nl$. For any $h\in\pi_1(M)-C$, the elements $g$ and $h$ do not commute; since $l<\mu$, and $\mu$ is a Margulis number, it follows that a loop representing $h$ must have length at least $\mu>l$. This shows that $\shortone(p)=l$, which is the first assertion of the lemma.

According to \cite[Proposition 13.1]{fourfree}, we have $p\in\frakG_M$ and  
\Equation\label{busted knee}
\vol \nbhd \nolimits_{\fraks_M(p )/2} (p) = B(\fraks_M (p )/2) - 2\kappa(\fraks_M (p )/2, l/2),
\EndEquation
where $\frakG_M$ is the set, and $\fraks_M$ the function with domain $\frakG_M$, whose definitions were reviewed in \ref{what about the workers}. We observed in \ref{what about the workers} that $\fraks_M$ is simply the restriction of $\nextone$ to $\frakG_M$. Thus (\ref{busted knee}) becomes $\vol \nbhd \nolimits_{\nextone(p )/2} (p) = B(\fraks_M (p )/2) - 2\kappa(\nextone (p )/2, l/2)$, which by the definition of $\Wtwo$ is equivalent to the second assertion of the present lemma.
\EndProof

The following result will be needed in the proofs of Lemmas \ref{cranky-poo} and \ref{coffee}.

\Lemma [cf. \cite{correction}]\label{on account of}
The function $\Wtwo$ is monotonically increasing (in the weak sense) in each of its  arguments.
\EndLemma

\Proof
Suppose that numbers $\lambda_1$, $\lambda_2$ and $l$ are given, with $l>0$ and $0<\lambda_1\le\lambda_2$.
Fix a line $L\subset\HH^3$ and a point $P\in L$, and let $Y$ denote the closed connected subset of $\HH^3$ bounded by the two planes that are perpendicular to $L$ and have distance $l/2$ from $P$. Then for $i=1,2$, if $N_i$ denotes the ball of radius $\lambda_i/2$ centered at $P$, we have $\vol(Y\cap N_i)=\Wtwo(\lambda_i,l)$. Since $\lambda_1\le\lambda_2$, we have $N_1\subset N_2$ and 
therefore
 $Y\cap N_1\subset Y\cap N_2$; hence $\Wtwo(\lambda_1,l)\le \Wtwo(\lambda_2,l)$. This shows that $\Wtwo$ is monotonically increasing  in  its first  argument. It is  monotonically increasing in its second  argument because $\kappa$  monotonically decreasing in its second argument (see \ref{what to call it}).
\EndProof

\NotationReviewRemarks\label{Boroczky number}
As in \cite{boroczky} and \cite{fourfree}, for any $n\ge2$ and any $R>0$ we shall denote by $h_n(R)$ the distance from the barycenter to a vertex of a regular hyperbolic $n$-simplex $\scrD_{n,R}$ with sides of length $2R$. For any $n\ge2$ the function $h_n(R)$ is strictly monotone increasing. Formulae for $h_2(R)$ and $h_3(R)$ are given in \cite[Subsection 9.1]{fourfree}.

For $R>0$
we define a function $\density(R)$ (denoted 
$d_3(R)$ in \cite{boroczky} and $d(R)$ in \cite{fourfree})
by
$$\density(R)=(3\beta(R)-\pi)(\sinh((2R)-2R)/\tau(r),$$
where the functions 
$$\beta(R)=\arcsec(\sech(2R)+2)\text{ and }\tau(R)=3\int_{\beta(R)}^{\arcsec3}\arcsech((\sec t)-2)\,dt$$
respectively give the dihedral angle and the volume of $\scrD_{n,R}$.

Suppose that $p$ is a point of a hyperbolic $3$-manifold $M$ and that $R$ is a positive number. 
%As was pointed out in  \ref{short and next}, $p$ is the center of a hyperbolic ball of radius $R$.
It is pointed out in \cite[Subsection 9.2]{fourfree}, as a consequence of the results on sphere packing proved in \cite{boroczky}, 
that if there is a radius-$R$ hyperbolic ball in $M$ with center $p$ (see \ref{nbhd}) 
then
$
\vol\nbhd_{ h_3(R)}(p)
\ge B(R)/\density(R)$. 
It is also observed in \cite[Subsection 9.4]{fourfree} that we have $B(h_3(R))\ge B(R)/\density(R)$ for any $R>0$.

A stronger version of the lower bound for 
$\vol\nbhd_{ h_3(R)}(p)$ 
mentioned above is established in \cite{fourfree}. For $R>0$ and $\rho>h_3(R)$, one defines
$$
\phi(R,\rho)=
\arcsin\bigg(\frac{\sqrt{\cosh^2\rho-\cosh^2 R}}{\sinh
\rho\cosh R}\bigg)-
\arcsin\bigg(\frac{\sqrt{\cosh^2 h_3(R)-\cosh^2 R}}{\sinh h_3(R)\cosh R}\bigg).
$$
and
$$
\Vbor(R,\rho)=\bigg(\frac{1-\cos\phi(R,\rho)}2\bigg) B(h_3(R) )+
\bigg(\frac{1+\cos\phi(R,\rho)}2\bigg)
\frac{B(R)}{d(R)}.
$$
Note that since 
$B(h_3(R))\ge B(R)/\density(R)$, we have $\Vbor(R,\rho)\ge B(R)/\density(R)$. 

According to \cite[Remark 9.6]{fourfree}, the function $\Vbor$ is monotone increasing in its second argument. Furthermore, \cite[Proposition 9.7]{fourfree} asserts that if $p$ is a point  of a hyperbolic $3$-manifold $M$,  if $R$ and $\rho$ satisfy  $\rho>h_3(R)$, if $p$ is the center of a hyperbolic ball of radius $R$, and if there is a point of $M$ whose distance from $p$ is  $\rho$, then $
\vol\nbhd_{ h_3(R)}(p)
\ge \Vbor(R,\rho)$.

\EndNotationReviewRemarks

\NotationRemark\label{is that you wthree?}
We define a function $\Wthree$ on $(0,\infty)^3$ by
$$\Wthree(\rho,R,\mu)=\begin{cases}
B(\max(0,\rho-R))\text{ if }\rho-R<\mu/2\\
B(\mu/2)\text{ if }\mu/2\le\rho-R<h_3(\mu/2)\\
\Vbor(\mu/2,\rho)+B(\min(\mu/2,(\rho-R-h_3(\mu/2))/2))\text{ if }h_3(\mu/2)\le\rho-R<
3 
h_3(\mu/2)\\
\Vbor(\mu/2,\rho)+B(\mu/2)/\density(\mu/2)\text{ if } \rho-R\ge
3 
h_3(\mu/2).
\end{cases}
$$
Note that $\Wthree$ is monotone decreasing (in the weak sense) in its second argument, since we have $\Vbor(\mu/2,\rho)\ge B(\mu/2)/\density(\mu/2)\ge B(\mu/2)$ whenever $\rho>\mu/2$ (see \ref{Boroczky number}). Furthermore, the latter fact, together with the fact that $\Vbor$ is monotone increasing in its second argument by \cite[Remark 9.6]{fourfree} (cf. \ref{Boroczky number}), 
implies that
$\Wthree$ is  (weakly) monotone increasing in its first argument.
\EndNotationRemark

The notation $\Wthree$ is meant to suggest that this function provides a lower bound for the volume of the complement of a suitable neighborhood of a point in a hyperbolic manifold under certain hypotheses. More precisely, we have the following result, the proof of which will incorporate ideas from the proofs of Lemmas 8.3 and 11.3 of \cite{fourfree}.

\Lemma\label{new distant thunder}
Let $M$ be a hyperbolic $3$-manifold, let $\mu$ be a Margulis number for $M$, let $\rho$ and $R$ be positive numbers. Let $p_0\in M$  and $p_1\in\Mthick(\mu)$ be points  such that %$\shortone(p_0)\ge\mu$ and 
$\dist(p_0,p_1)\ge\rho$. Suppose that
$R\ge\nextone(p_0)/2$.  Then
$$\vol (M-\nbhd\nolimits_R(p_0))\ge\Wthree(\rho,R,\mu).$$
\EndLemma

\Proof
We set $\rho_1=\dist(p_0,p_1)\ge\rho$.

If $\rho\le R$ we have $\Wthree(\rho,R,\mu)=0$, and the result is therefore trivial in this case. For the rest of the proof we will assume $\rho>R$. 

Set $s=\min(\rho-R,h_3(\mu/2))>0$.
Since $s+R\le\rho\le\rho_1=\dist(p_0,p_1)$, it follows from the triangle inequality (see \ref{nbhd}) that $\nbhd_s(p_1)\cap\nbhd_R(p_0)=\emptyset$,
%a contradiction to the hypothesis. Hence 
i.e. 
\Equation\label{nu garelick label}
M-\nbhd\nolimits_R(p_0)\supset \nbhd\nolimits_s(p_1)
\EndEquation

The 
hypothesis $p_1\in\Mthick(\mu)$ means that $\shortone(p_1)\ge\mu$. Hence
there is a hyperbolic ball in $M$ having radius $\mu/2$ and center $p_1$ (see \ref{short and next}).
In particular, if we set $s'=\min(\rho-R,\mu/2)>0$, so that $0<s'\le s$, then there is a hyperbolic ball in $M$ having radius $s'$ and center $p_1$. Hence $\vol\nbhd_{s}(p_1)\ge\vol\nbhd_{s'}(p_1)=B(s')=B(\min(\rho-R,\mu/2))$. But it follows from (\ref{nu garelick label}) that $\vol(M-\nbhd_R(p_0))\ge
\vol\nbhd_s(p_1)
$, and therefore
\Equation\label{easy on the garelick}
\vol(M-\nbhd\nolimits_R(p_0))\ge B(\min(\rho-R,\mu/2)). 
\EndEquation
In the case where $\rho-R<h_3(\mu/2)$, the left-hand side of (\ref{easy on the garelick}) is equal to 
$\Wthree(\rho,R,\mu)$, and hence the conclusion of the lemma is true in this case.

The remainder of the proof will be devoted to the case in which $\rho-R\ge h_3(\mu/2)$. In this case we have $s=h_3(\mu/2)$. We set $u=(\rho-s-R)/2\ge0$ and $\lambda=\nextone(p_0)$. By hypothesis we have $R\ge\lambda/2$.

According to \ref{short and next}, we may choose a short maximal cyclic subgroup $C$ of $\pi_1(M,p_0)$; furthermore, some non-trivial element of $C$ is represented by a loop $\alpha$ of length $\shortone(p_0)$, while some  element of $\pi_1(M,p_0)-C$ is represented by a loop $\beta$ of length $\lambda$. Since $\pi_1(M)$ is an \iccg, the elements $[\alpha]$ and $[\beta]$ of $\pi_1(M,p_0)$ do not commute. If we set $K=\alpha([0,1])\cup\beta([0,1])$, it follows that the image of the inclusion homomorphism $\pi_1(K)\to\pi_1(M)$ is non-abelian. But 
by Proposition \ref{my summer vocation},
each component of $\Mthin(\mu)=M-\Mthick(\mu)$ has an abelian fundamental group. Hence $K\not\subset\Mthin(\mu)$. We may therefore choose a point $p_2\in K\cap\Mthick(\mu)$. 

Let us define a continuous function $F$ on $M$ by $F(p)=\dist(p_0,p)$. Since the respective lengths of $\alpha$ and $\beta$ are $\shortone(p_0)\le\lambda$ and $\lambda$, we have $F(K)\subset[0,\lambda/2]$. In particular $F(p_2)\le\lambda/2\le R<\rho\le\rho_1$. On the  other hand, the definition of $\rho_1$ gives $F(p_1)=\rho_1$. Since 
%$R<\rho\le\rho_1$, and since 
$p_1$ and $p_2$ both lie in 
the  set $\Mthick(\mu)$, which is connected by Proposition \ref{my summer vocation},
we have $F(\Mthick(\mu))\supset[F(p_2),\rho_1]\supset[R,\rho_1]$. Since the definition of $u$ implies that $R\le R+u<\rho\le\rho_1$,  
there is a point $p^*\in\Mthick(\mu)$  such that $\dist(p_0,p^*) =F(p^*)=R+u$.
it then follows from the triangle inequality (see \ref{nbhd}) that $\nbhd_R(p_0)\cap\nbhd_u(p^*)=\emptyset$, i.e.
\Equation\label{take your limits}
M-\nbhd\nolimits_R(p_0)\supset\nbhd\nolimits_u(p^*).
\EndEquation

On the other hand, we have $\dist(p_1,p^*)\ge \dist(p_0,p_1)
-\dist(p_0,p^*)=\rho_1-(R+u)\ge\rho-(R+u)$, and the definition of $u$ implies that $\rho-(R+u)=u+s$.
Thus $\dist(p_1,p^*)\ge u+s$, and the triangle inequality gives 
\Equation\label{start your engines}
\nbhd\nolimits_s(p_1)\cap\nbhd\nolimits_u(p^*)=\emptyset.
\EndEquation

From (\ref{nu garelick label}), (\ref{take your limits}) and (\ref{start your engines}), it follows that
\Equation\label{enough already}
\vol(M-\nbhd\nolimits_R(p_0))\ge\vol(\nbhd\nolimits_s(p_1))+\vol(\nbhd\nolimits_u(p^*)).
\EndEquation

We have observed that $p_1$ is the center of
%hypothesis $p_1\in\Mthick(\mu)$ means that $\shortone(p_1)\ge\mu$. Hence
%there is 
a hyperbolic ball in $M$ having radius $\mu/2$. Since the point $p_0$ lies at a distance  $\rho_1$ from $p_1$, and since $s=h_3(\mu/2)$ in the present case, it follows from \cite[Proposition 9.7]{fourfree}, which was reviewed in \ref{Boroczky number}, that $
\vol\nbhd\nolimits_s(p_1)
\ge\Vbor(\mu/2,\rho_1)$; as $\Vbor$ is monotone increasing in its second argument by \cite[Remark 9.6]{fourfree}, which was also reviewed in \ref{Boroczky number}, we obtain
\Equation\label{when o lord}
\vol\nbhd\nolimits_s(p_1)
\ge\Vbor(\mu/2,\rho).
\EndEquation

Since $p^*\in\Mthick(\mu)$, we have $\shortone(p^*)\ge\mu$, so that $p^*$ is the center of a hyperbolic ball of radius $\mu/2$;
if we set $u'=\min(u,\mu/2)>0$, it follows that $p^*$ is the center of a hyperbolic ball of radius $u'$, and hence
\Equation\label{why not lordy}
\vol\nbhd\nolimits_u(p^*)\ge \vol\nbhd\nolimits_{u'}(p^*)
=B(u'). 
\EndEquation

From (\ref{enough already}), (\ref{when o lord}) and (\ref{why not lordy}), we obtain 
\Equation\label{inky dinky doo}
%\begin{aligned}
\vol(M-\nbhd\nolimits_R(p_0))\ge\Vbor(\mu/2,\rho)+B(u').
\EndEquation
In the subcase where $h_3(\mu/2)\le\rho-R<
3 
h_3(\mu/2)$, the right-hand side of (\ref{inky dinky doo}) is equal to
$\Vbor(\mu/2,\rho)+B(\min((\rho-h_3(\mu/2)-R)/2,\mu/2))=\Wthree(\rho,R,\mu)$ according to the definitions, and the conclusion of the lemma follows in this subcase.

Finally, consider the subcase in which $\rho-R\ge
3 
h_3(\mu/2)$=3s. We then have $u=(\rho-s-R)/2\ge s$. Hence $\nbhd_u(p^*)\supset\nbhd_s(p^*)$. But since $p^*$ is the center of a hyperbolic ball of radius $\mu/2$, it follows from the discussion in \ref{Boroczky number} that $
\vol\nbhd_s(p^*)=\vol\nbhd_{h_3(\mu/2)}(p^*)
\ge B(\mu/2)/\density(\mu/2)$. Hence
\Equation\label{last one maybe}
\vol\nbhd\nolimits_u(p^*)
\ge B(\mu/2)/\density(\mu/2).
\EndEquation
In this subcase,  (\ref{enough already}), (\ref{when o lord}) and (\ref{last one maybe}) give
$$
%\begin{aligned}
\vol(M-\nbhd\nolimits_R(p_0))\ge\Vbor(\mu/2,\rho)+B(\mu/2)/\density(\mu/2)=\Wthree(\rho,R,\mu),
$$
where the last equality is simply the definition of $\Wthree(\rho,R,\mu)$ in the subcase $\rho-R\ge 3h_3(\mu/2)$.
\EndProof
%p'p_2

%\Number\label{diddly}
%\EndNumber

\Notation\label{chi def}
Suppose that $k>2$ is an integer, 
that $\alpha$, $\delta$, $\lambda$ and $\mu$ are positive numbers, and that $D$ is a number such that $D\ge\delta$ and $D>\f_1(\lambda)$.
Since $\alpha>0$ and $D>\delta>0$, the quantity $\Wone(\alpha,\delta,D)$ is defined. Since $D>\fthree(\lambda)$, we have $Q(\lambda)+Q(D)<1/2$
by \ref{Q-def},
so that $\xi_{k-2}(Q(\lambda)+Q(D))$ is defined and strictly positive (see \ref{reformulated big radius corollary from cusp}). We set
$$
\waschi_k(\alpha, \lambda,\delta,D,\mu)
=\Wone(\alpha,\delta,D)+\Wthree(\xi_{k-2}(Q(\lambda)+Q(D)), \lambda/2,\mu).$$
\EndNotation

\Lemma\label{stoppeth two}
Let $k>2$ be an integer, and let $M$ be a closed, orientable hyperbolic $3$-manifold such that
$\pi_1(M)$ is $k$-free. Let $\mu$ be a Margulis number for $M$. Let $p$ be  a point of $ M$, set $\lambda=\nextone(p)$, and let
$\alpha$ be a number with $0<\alpha\le\lambda$. 
Let $\delta$ be a number
strictly greater than $\alpha/4$,
and suppose that 
$\ell_M\ge\delta$ (in the notation of \ref{nbhd}). 
Then there is a number 
$D$ such that 
$D\ge\delta$ and $D>\f_1(\lambda)$
(so that 
$\waschi_k(\alpha, \lambda,\delta,D,\mu)$
is defined by \ref{chi def}), with
$
\vol M\ge 
\waschi_k(\alpha, \lambda,\delta,D,\mu)
$.
\EndLemma

\Proof
By \ref{short and next} we may fix a short maximal cyclic subgroup $C$ of $\pi_1(M,p)$. We choose a generator $g$ of $C$, and we define $D$ to be the minimal length of a loop based at $p$ that represents  $g$. The hypotheses of the present proposition imply that  $\alpha<4\delta$, and that the conjugacy class of $ g$ is represented by a closed geodesic in $M$ having length
greater than $\delta$. Hence it follows from \ref{reformulated surgeon general} that 
$(\alpha,\delta,D)\in\frakW $ (i.e. $D\ge\delta$), 
and that $\vol \nbhd_{\alpha/2}( p)\ge\Wone(\alpha,\delta,D)$.

Next,  note that according to 
Lemma
\ref{when m is two} we have $Q(\lambda)+Q(\shortone(p))<1/2$
(so that $\xi_{k-2}(Q(\lambda)+Q(\shortone(p)))$ is defined), and
there is a $\mu$-thick point $p'\in M$ such that
$\dist (p,p')\ge\xi_{k-2}(Q(\lambda)+Q(\shortone(p)))$. The definition of $\shortone(p)$ implies that $D\ge\shortone(p)$; since $Q$ is strictly monotone decreasing, and $\xi_{k-2}$ is strictly monotone 
increasing 
on $(0,1/2)$, it follows that $Q(\lambda)+Q(D)<1/2$ (so that $\xi_{k-2}(Q(\lambda)+Q(D))$ is defined), and $\dist (p,p')\ge\xi_{k-2}(Q(\lambda)+Q(D))$.
The inequality $Q(\lambda)+Q(D)<1/2$ is equivalent to $D>\fthree(\lambda)$
by \ref{Q-def};
as we have already observed that $D\ge\delta$, it now follows that 
$(\alpha, \lambda,\delta,D,\mu)$ 
is in the domain of $\waschi_k$.
Since $\lambda=\nextone(p)$ and $p'\in\Mthick(\mu)$, we may apply  \ref
{new distant thunder} with $p$, $p'$, $\xi_{k-2}(Q(\lambda)+Q(D))$  and
$\lambda/2$
playing the respective roles of $p_0$, $p_1$, $\rho$ and $R$, to deduce that $\vol (M-\nbhd_{\lambda/2}(p))\ge \Wthree(\xi_{k-2}(Q(\lambda)+Q(D)),\lambda/2,\mu)$. 
Since $\lambda\ge\alpha$, it now follows that
$$
\begin{aligned}
\vol M&=\vol\nbhd\nolimits_{\lambda/2}(p)+\vol (M-\nbhd\nolimits_{\lambda/2}(p))\\
&\ge \vol
(\nbhd\nolimits_{\alpha/2}(p))
+\vol (M-\nbhd\nolimits_{\lambda/2}(p))\\
&\ge \Wone(\alpha,\delta,D)+\Wthree(\xi_{k-2}(Q(\lambda)+Q(D)),\lambda/2,\mu)\\
&
=
\waschi_k(\alpha, \lambda,\delta,D,\mu).
\end{aligned}
$$
\EndProof

\NotationRemarks\label{before stoppeth one}
Let $k>2$ 
be
an integer, and let $l$, $h$  and $\mu$ be positive real numbers. Since $\fthree(l)+h>\fthree(l)$, it follows from \ref{Q-def} that
$Q(\fthree(l)+h)+Q(l)<1/2$, and hence $\xi_{k-2}(Q(\fthree(l)+h)+Q(l))$ is a well-defined positive number. We set 
$$\waspsi_k( h,l,\mu)=\Wtwo(\fthree(l)+h,l)+\Wthree(\xi_{k-2}(Q(\fthree(l)+h)+Q(l)),(\fthree(l)+h)/2,\mu)$$
(where the functions $\Wtwo$ and $\Wthree$ are defined by 
\ref{short geo def} 
and \ref{is that you wthree?}). 
\EndNotationRemarks

%y

\Lemma\label{stoppeth one}
 Let $k>2$ be an integer, and let $M$ be a closed, orientable hyperbolic $3$-manifold such that
$\pi_1(M)$ is $k$-free. Let $\mu$ be a Margulis number for $M$. Let $c$ be a closed geodesic in $M$, let $l$ denote the length of $c$, and
 suppose that
$l<\log3$. Then for any  point  $p$  of $|c|$,
we have $\nextone(p)>\f_1(l)$, and
$\vol M
\ge \waspsi_k(\nextone(p)-\f_1(l),l,\mu)$. In particular,
there is  a positive real number 
$h$ such that
%$Q(x)+Q(l)<1/2$ (so that $
$\vol M
\ge \waspsi_k(h,l,\mu)$.
\EndLemma

Note that the first assertion of the Lemma \ref{stoppeth one}, that $\nextone(p)>\f_1(l)$, is very slightly stronger than the conclusion of Corollary \ref{nebulizer}, which depends only on the assumption that $\pi_1(M)$ is $2$-free.

\Proof[Proof of Lemma \ref{stoppeth one}]
Since $\pi_1(M)$ is in particular $2$-free, $\log3$ is a Margulis number for $M$ by 
\cite[Corollary 4.2]{acs-surgery}. 
Let us fix a closed geodesic $c$ of length $l$, and a point $p$ of $|c|$. 
We set $x=\nextone(p)$. We apply 
Lemma \ref{surgeon special},
letting $\log3$ play the role of the Margulis number that is denoted $\mu$ in 
that lemma. This gives 
$\shortone(p)=l$ and $\vol\nbhd_{x/2}(p)=\Wtwo(x,l)
$.

Next, note that according to 
Lemma
\ref{when m is two} we have $Q(x)+Q(l)=Q(\nextone(p))+Q(\shortone(p))<1/2$ (so that $\xi_{k-2}(Q(x)+Q(l))$ is defined), and there is a $\mu$-thick point $p'\in M$ such that
 $\dist (p,p')\ge\xi_{k-2}(Q(x)+Q(l))$. Since $x=\nextone(p)$, we may apply Lemma \ref
{new distant thunder} with $p$, $p'$, $\xi_{k-2}(Q(x)+Q(l))$ and $x/2$ playing the respective roles of $p_0$, $p_1$, $\rho$ and $R$, to deduce that $\vol (M-\nbhd_{x/2}(p))\ge \Wthree(\xi_{k-2}(Q(x)+Q(l)),x/2,\mu)$. It now follows that
\Equation\label{unfrench fries}
\begin{aligned}
\vol M&=
\vol(\nbhd\nolimits_{x/2}(p))
+\vol (M-\nbhd\nolimits_{x/2}(p))\\
&\ge \Wtwo(x,l)+\Wthree(\xi_{k-2}(Q(x)+Q(l)),x/2,\mu).
\end{aligned}
\EndEquation

Now since $Q(x)+Q(l)<1/2$, we  have 
%$Q(x)<1/2-Q(l)=Q(\fthree(l))$; since $Q$ is strictly monotone decreasing on $(0,\infty)$, this gives 
$x>\fthree(l)$
by \ref{Q-def}. 
We  therefore have  $h:=x-\fthree(l)>0$.
The definition of $\waspsi$ now gives $\waspsi_k( h,l,\mu)=\Wtwo(x,l)+\Wthree(\xi_{k-2}(Q(x)+Q(l)),x/2,\mu)$, which with (\ref{unfrench fries}) gives $\vol M\ge\waspsi_k(h,l.\mu)$.
\EndProof

\section{Volumes, diameters, and Margulis numbers}\label{Margulis section}

\NotationRemarks\label{what's tVbor}
For any positive numbers $R$ and $\rho$ we set
$$\tVbor(R,\rho)=\begin{cases}
\begin{aligned}
\Vbor(R,\rho)\qquad&\text{ if }&\rho> h_3(R)\\
B(R)/\density(R)\qquad&\text{ if }& \rho\le h_3(R).
\end{aligned}
\end{cases}
$$

Note that the strict monotonicity of $\Vbor$ in its second argument (cf. \ref{Boroczky number}), together with the inequality $\Vbor(R,\rho)\ge B(R)/\density(R)$ (cf. \ref{Boroczky number}), implies that $\tVbor$ is monotone increasing, in the weak sense, in its second argument. 
\EndNotationRemarks

\Proposition\label{all about tVbor}
Let $p$ be a point  of a hyperbolic $3$-manifold $M$, and let $R$ and $\rho$ be positive numbers. Suppose that $R\le\shortone(p)/2$, and that there is a point of $M$ whose distance from $p$ is at least  $\rho$. Then $
\vol\nbhd_{ h_3(R)}(p)
\ge \tVbor(R,\rho)$.
\EndProposition

\Proof
Since $R\le\shortone(p)/2$, there is a hyperbolic ball in $M$ having radius $R$ and center $p$ (see \ref{short and next}). 
If $\rho\le h_3(R)$, the conclusion now follows from
 the inequality $\vol
\nbhd_{ h_3(R)}(p)
\ge  B(R)/\density(R)$, which was reviewed in \ref{Boroczky number}. Now suppose that $\rho<h_3(R)$. Since $M$ is connected, and some point of $M$ has distance at least $\rho$ from $p$, there is a point of $M$ whose distance from $p$ is exactly $\rho$. Hence \cite[Proposition 9.7]{fourfree}, which was also reviewed in \ref{Boroczky number}, gives $
\vol\nbhd_{ h_3(R)}(p)
\ge \Vbor(R,\rho)$, which in this case is equivalent to the conclusion of the present proposition.
\EndProof

\Proposition\label{better diameter from thick volume}
Let $M$ be a closed, orientable hyperbolic $3$-manifold. Let $\mu$ be a Margulis number for $M$
(so that $\Mthick(\mu)\ne\emptyset$ by Proposition \ref{my summer vocation}). 
Let 
$\Delta  $ denote the extrinsic diameter 
(see \ref{nbhd})
of $\Mthick(\mu)$, and  assume that $\Delta  \ge\mu$.
Let $R$ be a positive number such that 
$R\le\max_{p\in M}(\shortone(p)/2)$.
Then 
$$\vol M\ge
\tVbor(R,\Delta  /2)+\min(B(\mu/2),2B(\max(0,\Delta  /2-h_3(R)))$$
(where the 
function $B$ is defined as in \ref{bee of eggs} and the functions $h_3$ is defined as in \ref{Boroczky number}).
\EndProposition

\Proof
According to the 
hypothesis,
we may fix a point $p_0\in M$ such that $\shortone(p_0)\ge2R$. 

According to the definition of $\Delta  $, we may fix points $p_1,p_2\in \Mthick(\mu)$ such that $\dist(p_1,p_2)=\Delta  $. We may suppose the $p_i$ to be indexed so that $\dist(p_0,p_1)\ge\dist(p_0,p_2)$. Since $\Delta  \le \dist(p_0,p_1)+\dist(p_0,p_2)$,  it follows that $\dist(p_0,p_1)\ge \Delta  /2$. Now since 
$\shortone(p_0)\ge2R$,
%$\nbhd_R(p_0)$  is intrinsically isometric to a ball of radius $R$ in $\HH^3$,
it follows from 
Proposition \ref{all about tVbor}
that   $N:=\nbhd_{h_3(R)}(p_0)$ has volume at least $\tVbor(R,\Delta  /2)$.

For $i=1,2$, set $d_i=\max(\dist(p_i,p_0),h_3(R))$. Since $d_i\ge h_3(R)$ we have $e_i:=d_i-h_3(R)\ge0$. Since $\dist(p_0,p_1)\ge\dist(p_0,p_2)$, we have $d_1\ge d_2$ and hence $e_1\ge e_2$. On the other hand, we have $\dist(p_1,p_2)\le\dist(p_1,p_0)+\dist(p_0,p_2)\le d_1+d_2$.

For $i=1,2$, set $U_i=\nbhd_{e_i}(p_i)$ (so that in particular $U_i=\emptyset$ if $e_i=0$). We claim that
\Equation\label{fibber mcgee}
U_i\cap N=\emptyset
\EndEquation
for $i=1,2$. 
To prove (\ref{fibber mcgee}) for a given $i$, first note that if $e_i=0$, the assertion is immediate since $U_i=\emptyset$. Now suppose that $e_i>0$.
Then we have $d_i>h_3(R)$, which in view of the definition of $d_i$ means that $\dist(p_i,p_0)>h_3(R)$ and that $d_i=\dist(p_i,p_0)$. Hence 
$\dist(p_i,p_0)=e_i+h_3(R)$, which by the triangle inequality (see \ref{nbhd}) implies that $\nbhd_{e_i}(p_i)\cap\nbhd_{h_3(R)}(p_0)=\emptyset$, which establishes 
\ref{fibber mcgee}.

Consider the case in which $e_1\ge\mu/2$. In this case we have $W:=\nbhd_{\mu/2}(p_1)\subset U_1$; in view of \ref{fibber mcgee}, it follows that $W\cap N=\emptyset$. Hence $\vol M\ge(\vol W)+(\vol N)$.
 We have seen that $\vol N\ge\tVbor(R,\Delta  /2)$. Since $p_1\in\Mthick(\mu)$, the set $W$ is intrinsically isometric to a ball of radius $\mu/2$ in $\HH^3$. Hence $
\vol W
=B(\mu/2)$. Thus in this case we have $\vol M\ge\tVbor(R,\Delta  /2)+B(\mu/2)$, which implies the conclusion of the proposition.

The rest of the proof will be devoted to the case in which $e_1<\mu/2$. Since $e_1\ge e_2$, we have $e_2<\mu/2$;
 and since $\mu\le\Delta$ by hypothesis, it follows that
$e_1+e_2<\mu\le \Delta  =\dist(p_1,p_2)$. 
Hence by
the triangle inequality (see \ref{nbhd}), we have
% it follows that
 $\nbhd_{p_1}(e_1)\cap \nbhd_{p_2}(e_2)=\emptyset$, i.e.
\Equation\label{golly}
U_1\cap U_2=\emptyset.
\EndEquation

According to (\ref{fibber mcgee}) and (\ref{golly}), the sets $U_1$, $U_2$ and $N$ are pairwise disjoint. Hence 
$\vol M\ge
(\vol U_1)+(\vol U_2)+(\vol N)$.
We have seen that $\vol N\ge\tVbor(R,\Delta  /2)$. For $i=1,2$, since $p_i\in\Mthick(\mu)$, and since $e_i<\mu/2$, the set $U_i$ is intrinsically isometric  to a ball in $\HH^3$ of radius $e_i$, and hence $
\vol U_i
=B(e_i)$. Thus we obtain $\vol M\ge\tVbor(R,\Delta  /2)+B(e_1)+B(e_2)$. But the function $B(x)=\pi(\sinh(2x)-2x)$ is convex for $x\ge0$, so that $(B(e_1)+B(e_2))/2\ge B((e_1+e_2)/2)$. We therefore have $\vol M\ge\tVbor(R,\Delta  /2)+2B((e_1+e_2)/2)$. But we have $e_1+e_2=d_1+d_2-2h_3(R)\ge\dist(p_1,p_0)+\dist(p_2,p_0)-2h_3(R)\ge\dist(p_1,p_2)-2h_3(R)=\Delta  -2h_3(R)$. As $e_1$ and $e_2$ are non-negative, we have $(e_1+e_2)/2\ge\max(0,\Delta  /2-h_3(R))$. Since $B(x)$ is also monotone increasing for $x\ge0$, it now follows that $\vol M\ge2B(\max(0,\Delta  /2-h_3(R)))$, which implies the conclusion of the proposition in this case. 
\EndProof
%UV

\Reformulation\label{reformulated better diameter from thick volume} 
For applications of Proposition \ref{better diameter from thick volume}, it will be convenient to define a function  $\Wfive$ on $(0,\infty)^3$ by 
$\Wfive(R,x,\mu)=\tVbor(R,x /2)+\min(B(\mu/2),2B(\max(0,x /2-h_3(R))))$. 
In terms of this definition, we may reformulate Proposition \ref{better diameter from thick volume} as follows: 
if $M$ is a closed, orientable hyperbolic $3$-manifold, if $\mu$ is a Margulis number for $M$,  if $\Delta $ denotes the extrinsic diameter of $\Mthick(\mu)$,   if $\Delta \ge\mu$,
and if $R$ is a positive number such that 
$R\le\max_{p\in M}(\shortone(p)/2)$, then 
$\vol M\ge \Wfive(R,\Delta,\mu)$. 
(The subscript ``D-R'' stands for ``diameter-radius,'' as this function allows one to estimate the volume of a hyperbolic $3$-manifold in terms of the diameter of its thick part and the radius of a certain hyperbolic ball in the manifold (cf. \ref{short and next}).)
\EndReformulation

Proposition \ref{better diameter from thick volume} 
(or its reformulation \ref{reformulated better diameter from thick volume})
will be applied in conjunction with the following result, which is stated in terms of the function $Q$ defined in \ref{Q-def}:

\Proposition\label{cover your face}
Let $k > 2$ be an integer and let $M$ be a closed, orientable hyperbolic $3$-manifold such that $\pi_1(M)$ is $k$-free. Let $\mu_0$ be a Margulis number for $M$, 
and let $\Delta$ denote the extrinsic diameter (see \ref{nbhd}) of $\Mthick(\mu_0)$ in $M$. Let $\mu$ be a positive
real number such that $Q(\mu)+Q(\mu_0)\ge1/2$ and
$2Q(\mu) +(k-2)Q(2\Delta)\ge1/2$.
Then $\mu$ is itself a Margulis number for $M$.
\EndProposition

\Proof
This is a paraphrase of \cite[Corollary 10.3]{fourfree}, using the function $Q$. The quantities denoted here by $\mu_0$ and $\mu$ are denoted by $\mu$ and $\lambda$ respectively in \cite[Corollary 10.3]{fourfree}.
\EndProof

\Notation\label{that's true too}
For each integer $k>2$, we define a function $\g_k$ on $(\log3,\infty)$ by 
$$\g_k(x)=\frac12 Q^{-1}\bigg(\frac{1/2-2Q(x)}{k-2}\bigg).$$ 
The function $\g_k$ is well defined for each $k>2$ since $0<Q(x)<1/4$ whenever $x>\log3$.
\EndNotation

\Proposition\label{new kite}
Let 
$k>2$
 be an integer, and let $V$, $R$,  and $\mu$ be positive real numbers. 
Suppose that $\mu>\log3$, that 
$\Wfive(R,\g_k(\mu), \fthree(\mu)
%\mu_0
)\ge V$,
and that $\g_k(\mu)\ge\fthree(\mu)$ (where $\fthree(\mu)$ is defined by \ref{Q-def}).
Let $M$ be a closed, orientable hyperbolic $3$-manifold such that $\pi_1(M)$ is $k$-free and  $\max_{p\in M}(\shortone(p)/2)\ge R$.
Then either $\vol M\ge V$, or $\mu$ is a Margulis number for $M$.
\EndProposition

\Remark
Proposition \ref{new kite} is the result, referred to in the Introduction, which gives strictly stronger information than Lemma 10.4 of \cite{fourfree}.
\EndRemark

\Proof[Proof of Proposition \ref{new kite}]
Set $\mu_0=\fthree(\mu)$ and $G=\g_k(\mu)$. According to the definitions of $\fthree$ and $\g_k$, we have $Q(\mu)+Q(\mu_0)=1/2$ and $2Q(\mu) +(k-2)Q(2G)=1/2 $.

Since $\mu>\log3$, we have $
\mu_0<\log3
$.
But since $\pi_1(M)$ is in particular $2$-free, it  follows from \cite[Corollary 4.2]{acs-surgery} that $\log3$ is a Margulis number for $M$; hence   $\mu_0$ is a Margulis number.
%If $\mu$ is a Margulis number, the conclusion of the lemma is true. We shall assume for the rest of the proof that $\mu_0$ is a Margulis number. 

We shall denote by $\Delta$ the extrinsic diameter of $\Mthick(\mu_0)$.

Consider the case in which $\Delta\ge G$. Since $G\ge\mu_0$ by hypothesis, we have $\Delta\ge\mu_0$. Since $\max_{p\in M}(\shortone(p)/2)\ge R$ by the hypothesis of the present lemma, we may apply  \ref{reformulated better diameter from thick volume}, with the Margulis number $\mu_0$ playing the role of $\mu$ in that proposition, to deduce that $\vol M\ge \Wfive (R,\Delta,\mu_0)$.

The hypothesis gives $\Wfive (
R,G,\mu_0
)\ge V$. On the other hand,
since $\tVbor$ is monotone increasing, in the weak sense, in its second argument (see \ref{what's tVbor}),  and $B$ is (strictly) monotone increasing, the function $\Wfive$ is monotone increasing, in the weak sense,
in its second argument.
Since  $\Delta\ge G$, we now have $\vol M\ge \Wfive (
R,\Delta,\mu_0)\ge \Wfive (R,G,\mu_0)
\ge V$, and the 
first alternative of the
conclusion of the lemma is true in this case.

There remains the case in which $\Delta<G$. Since $Q$ is monotone decreasing, 
in this case 
we have 
%the hypothesis of the present lemma gives 
$
2Q(
\mu
) +(k-2)Q(2\Delta)>2Q(\mu) +(k-2)Q(2G)=1/2 
$.
Since  $Q(\mu)+Q(\mu_0)=1/2$,
 it now follows from Proposition \ref{cover your face} that $\mu$ is a Margulis number for $M$, and 
the second alternative of the conclusion holds.
\EndProof
%\mu f\mu_0 H

\Notation\label{ikea}
For any integer $k>2$, we define an interval $I_k\subset\RR$ as follows: if $k<7$ we set 
$$I_k=\bigg(\log3,\log\bigg(\frac{20}{7-k}-1\bigg)\bigg),$$
and if $k\ge7$ we set $I_k=(\log3,\infty)$. It follows from the definitions of the functions $\g_k$ and $Q$ that $I_k$ is the set of all real numbers $x>\log3$ such that $g_k(x)>\log3$.
\EndNotation

In terms of the intervals $I_k$, we have the following corollary to Proposition \ref{new kite}:

\Corollary\label{marseillaise}
Let 
$k>2$
 be an integer, and let $V$, $R$,  and $\mu$ be positive real numbers. 
Suppose that $\mu\in I_k$ and that 
$\Wfive(R,\g_k(\mu), \fthree(\mu)
%\mu_0
)\ge V$.
Let $M$ be a closed, orientable hyperbolic $3$-manifold such that $\pi_1(M)$ is $k$-free and  $\max_{p\in M}(\shortone(p)/2)\ge R$.
Then either $\vol M\ge V$, or $\mu$ is a Margulis number for $M$.
\EndCorollary

\Proof
It suffices to show that the hypotheses of Proposition \ref{new kite} hold with the given choices of $k$, $V$, $R$, $\mu$ and $M$. Since $\mu\in I_k$, we have $\g_k(\mu)>\log3$ by \ref{ikea}. But, again since $\mu\in I_k$, we have $\mu>\log3$, and hence $\fthree(\mu)<\log3$. It now follows that $\g_k(\mu)>\fthree(\mu)$, which is one of the hypotheses of Proposition \ref{new kite}. The other hypotheses of Proposition \ref{new kite} are included in those of the present corollary. 
\EndProof

\section{Very short geodesics}\label{new very short section}

\Lemma\label{very short lemma}
Let $k\ge3$ be an integer, and let $\wasdelta$ be a positive real number less than $\min(0.7,\log(k-1)/2)$. Set 
$$C=\frac4{\cosh(\wasdelta/2)e^\wasdelta(e^\wasdelta+3)}.$$
Suppose that $M$ is a 
closed, orientable hyperbolic $3$-manifold such that $\pi_1(M)$ is $k$-free, and that 
the quantity $\ell_M$ (see \ref{nbhd}) is at most 
$\wasdelta$. Let $\mu$ be a Margulis number for $M$.
Then
$$\vol M\ge\pi C-\frac{\pi\wasdelta}2+B\bigg(\min\bigg(\frac\mu2,\frac{\log(k-1)}2-\wasdelta\bigg)\bigg),$$
where the function $B(x)$ is defined as in \ref{bee of eggs}.
\EndLemma

\Proof
We set $l=\ell_M\le\wasdelta$, and we fix a closed geodesic $c$ in $M$ of  length $l$. 
We denote by $R$ the  tube radius of $c$ (see \ref{nbhd}).
Since $\pi_1(M)$ is in particular $2$-free, it follows from \cite[Corollary 4.2]{acs-surgery} that $\log3$ is a strong Margulis number for $M$ in the sense defined in \cite{accs}; according to \cite[Proposition 10.1]{accs}, the fact that $\log3$ is a strong Margulis number implies that
\Equation\label{inadequate}
\cosh2R\ge\frac{e^{2l}+2e^l+5}{(\cosh(l/2))(e^l-1)(e^l+3)}>\frac{8}{(\cosh(l/2))(e^l-1)(e^l+3)}.
\EndEquation
(The fact that (\ref{inadequate}) holds when $\pi_1(M)$ is $2$-free was used in the proof of \cite[Lemma 13.4]{fourfree}, but the argument given there was incomplete.)
The mean value theorem gives $e^l-1=le^{l_0}$ for some $l_0\in(0,l)$. Since $l\le\wasdelta$, it follows that $e^l-1<le^\wasdelta$, and (\ref{inadequate}) then yields
\Equation\label{dogs}
\cosh 2R>\frac{8}{(\cosh(\wasdelta/2))(le^\wasdelta)(e^\wasdelta+3)}=\frac {2C}l.
\EndEquation

Since $\wasdelta\le0.7$, the definition of $C$ gives $C\ge0.373\ldots$. In particular we have $2C>0.7>\wasdelta\ge l$. Hence there is a unique 
strictly
positive number $R_0$ such that $\cosh2R_0=2C/l$. It follows from (\ref{dogs}) that $R_0<R$.
Since
$R_0>0$, 
the geodesic $c$ is simple 
and there is a
 tube $T$ having core $c$ and radius $R_0$
(see \ref{nbhd}). 
We have 
$$\vol T=\pi l\sinh^2R_0=
\pi l\bigg(\frac12\cosh2R_0-\frac12\bigg)=
\pi l\bigg(\frac{C}{l}-\frac12\bigg)\ge\pi C-\frac{\pi\wasdelta}2.
$$

Now fix a point 
$p\in|c|$.
Since the extrinsic diameter (see \ref{nbhd}) of 
$|c|$ 
is at most $l/2$, and every point of $T$ is at a distance at most $R_0$ from a point of 
$|c|$, 
we have $T\subset\nbhd_{R_0+l/2}(p)$. Hence
\Equation\label{cats}
\vol\nbhd\nolimits_{R_0+l/2}(p)\ge\vol T
\ge\pi C
-\frac{\pi\wasdelta}2.
\EndEquation

An orientation of $c$ defines a non-trivial element $u$ of the torsion-free group $\pi_1(M,p)$, which is represented by a loop of length $l$, and we may regard $\{u\}$ as a one-element independent set. Applying
Proposition \ref{big radius corollary from cusp} with $m=1$, we obtain
a $\mu$-thick point $p'\in M$ such that
 $\rho:=\dist (p,p')$ satisfies
$(k-1)Q(2\rho)+Q(l)
\le1/2$, i.e.
$$\frac{k-1}{1+e^{2\rho}}+\frac1{1+e^l}\le\frac12.$$
Solving the latter inequality for $e^{2\rho}$, we obtain
$$
\begin{aligned}
e^{2\rho}&\ge\frac{(2k-3)e^l+(2k-1)}{e^l-1}\\
&=(2k-3)+\frac{4k-4}{e^l-1}\\
&>\frac{4k-4}{e^l-1}>\frac{4k-4}{le^\wasdelta},
\end{aligned}
$$
so that 
\Equation\label{avenge}
\rho>\frac12\log(4k-4)+\frac12\log(1/l)-\wasdelta/2.
\EndEquation

In view of the definition of $R_0$, we have $e^{2R_0}\le2\cosh(2R_0)=4C/l$, and hence 
\Equation\label{my man yorick}
R_0\le\frac12\log(4C)+\frac12\log(1/l).
\EndEquation
Combining (\ref{avenge}) and (\ref{my man yorick}), 
we find
\Equation\label{i'll need it all right}
\rho-(R_0+l/2)>\frac12\log(4k-4)-\frac12\log(4C)-\wasdelta/2-l/2.
\EndEquation
But we have $l\le\wasdelta$, and the definition of $C$ immediately implies 
$C<1$. 
Hence (\ref{i'll need it all right}) gives $\rho-(R_0+l/2)>\log(k-1)/2-\wasdelta$, i.\,e.
\Equation\label{now here, still, is the thing}
\dist(p,p')=\rho>(R_0+l/2)+\bigg(\frac{\log(k-1)}2
-\wasdelta\bigg).
\EndEquation
By hypothesis we have $(\log(k-1))/2-\wasdelta>0$. In view of (\ref{now here, still, is the thing}) and the triangle inequality (cf. (see \ref{nbhd}), we obtain
$\nbhd\nolimits_{R_0+l/2}(p)\cap\nbhd\nolimits_{\log(k-1)/2-\wasdelta}(p')=\emptyset$, which implies
\Equation\label{pig and thistle}
\vol M\ge \vol(\nbhd\nolimits_{R_0+l/2}(p))+\vol(\nbhd\nolimits_{\log(k-1)/2-\wasdelta}(p')).
\EndEquation

Since $p'\in\Mthick(\mu)$, i.\,e.
 $\shortone(p')\ge\mu$,
there is a hyperbolic ball in $M$ having radius $\mu/2$ and center $p'$ (see \ref{short and next}). In particular, $p'$ is the center of a ball of radius $s:=\min(\mu/2,\log(k-1)/2-\wasdelta)$. We have $\nbhd_{\log(k-1)/2-\wasdelta}(p')\supset\nbhd_s(p')$, and hence
$
\vol\nbhd_{\log(k-1)/2-\wasdelta}(p')\ge\vol\nbhd_s(p')
=B(s)$; that is,
\Equation\label{not another paper!?}
\vol\nbhd\nolimits_{\log(k-1)/2-\wasdelta}(p')
\ge B\bigg(\min\bigg(\frac\mu2,\frac{\log(k-1)}2-\wasdelta\bigg)\bigg).
\EndEquation

The conclusion of the lemma follows from 
(\ref{pig and thistle}),
(\ref{cats}) and (\ref{not another paper!?}).
\EndProof

\Reformulation\label{reformulated very short lemma} 
For applications of Lemma \ref{very short lemma}, it will be convenient to define a %set $\frakW_4\subset\ZZ\times\RR^2$ and a 
function  $\Wfour$ on a subset of $\ZZ\times\RR^2$  
%$\frakW_4$ 
as follows. 
The domain of $\Wfour$ consists of all triples 
%We set $$\frakW_4=\{
$(k,\wasdelta,\mu)\in\ZZ\times\RR^2$ such that $k>2$, $\mu>0$, 
%\text{ 
and $0<\wasdelta<\log(k-1)/2$.
For each such triple $(k,\wasdelta,\mu)$
%\in\frakW_4$ 
we set
$$\Wfour (k,\wasdelta,\mu)=\frac{4\pi}{\cosh(\wasdelta/2)e^\wasdelta(e^\wasdelta+3)}-\frac{\pi\wasdelta}2+B\bigg(\min\bigg(\frac\mu2,\frac{\log(k-1)}2-\wasdelta\bigg)\bigg).$$
In terms of these definitions, we may reformulate Lemma \ref{very short lemma} as follows. Let  $k\ge3$ be an integer, let  $\wasdelta$ be a positive real number less than $\min(0.7,\log(k-1)/2)$, let $M$ be a closed, orientable hyperbolic $3$-manifold 
such that $\pi_1(M)$ is $k$-free  and $\ell_M\le\wasdelta$, 
and let $\mu$ be a Margulis number for $M$. Then $\vol M\ge\Wfour(k,\wasdelta,\mu)$.

(The subscript ``VSG'' stands for ``very short geodesic.'')

\EndReformulation

\section{
Sufficient conditions for being a lower bound} 
\label{my favorite monster}

The statement of the following proposition involves a number of real-valued functions, defined on certain intervals in $\RR$, that were introduced in earlier sections:
$\f_n$ (for $n\ge1$, Subsection \ref{Q-def}); $B$ (Subsection \ref{bee of eggs}); $\density$ (Subsection \ref{Boroczky number}); $\waschi_k$  (for $k>2$, Subsection \ref{chi def}); $\waspsi_k$
(for $k>2$, Subsection \ref{before stoppeth one}); $\Wfive$ (Subsection \ref{reformulated better diameter from thick volume}); and $\g_k$ (for $k>2$, Subsection \ref{that's true too}). The statement also involves the interval $I_k$, which was defined in Subsection \ref{ikea}, 
and the invariant $\ell_M$, which was defined in  Subsection \ref{nbhd}.

\Proposition\label{latest monster}
Let $k\ge4$ be an integer,  let $V_0>0$ be a real number, and let $T$ be a set of positive real numbers. Suppose that $\waslambdaminus$ is  a number with $\log7\le\waslambdaminus\le\log8$, and that
 $\scrM$ is a function defined on $[\waslambdaminus,\log8]$, and taking values in the interval $I_k$,
so that 
in particular, 
according to  \ref{ikea}, $\g_k(\scrM(\lambda))$ is defined whenever $\lambda\in[\waslambdaminus,\log8]$. 

Assume that the following conditions hold. 
\begin{enumerate}
\item $B(\fone(\waslambdaminus)/2)/\density(\fone(\waslambdaminus)/2)>V_0$.

\item
For each $\lambda\in[\waslambdaminus,\log8]$, there is an 
$R\in(0,\f_3(\lambda)/2]$ 
such that
$$\Wfive(R,\g_k(\scrM(\lambda)),\fthree(\scrM(\lambda)))
> V_0.$$

\item
For every $\lambda\in[\waslambdaminus,\log8]$, and every $\wasell\in T$,

either
\begin{enumerate}[(a)]
\item
there exist positive numbers $\alpha$ and $\delta$, with $\alpha\le\lambda$ and $\lambda/4<\delta<\wasell$,
such that for every 
number $D$ satisfying $D\ge\delta$ and $D>\fthree(\lambda)$,
we have
$\waschi_k(\alpha, \lambda,\delta,D,\scrM(\lambda))>V_0$, 
or
\item
we have $\wasell<\log3$, and for every $h>0$
there is a number $\mu\in(0, \scrM(\lambda)]$ such that
$\waspsi_k(h,\wasell, \mu)>V_0$.
\end{enumerate}
\end{enumerate}
Then for every closed, orientable hyperbolic $3$-manifold $M$ such
that $\pi_1(M)$ is $k$-free and $\ell_M\in T$ (see \ref{nbhd}),  we have $\vol M\ge V_0$.
\EndProposition

\Proof
Suppose that  $M$ is a closed, orientable hyperbolic $3$-manifold such
that $\pi_1(M)$ is $k$-free and $\ell_M\in T$. Since  $\pi_1(M)$ is $k$-free, one of the alternatives (i), (ii) of Corollary \ref{keyer corollary} holds (with the value of $\waslambdaminus$ given by the hypothesis of the present proposition). 

Consider first the case in which Alternative (i) of Corollary \ref{keyer corollary} holds, and fix a point $p\in M$ such that
$\shortone(p)\ge\fone(\waslambdaminus)$. Then $p$ is the center of a ball of
radius $\fone(\waslambdaminus)/2$ in $M$ (see \ref{short and next}).  By the
results reviewed in \ref{Boroczky number}, it follows that $\vol M\ge
B(\fone(\waslambdaminus)/2)/\density(\fone(\waslambdaminus)/2)$. But according to
Condition (1) of the present proposition, we have
$B(\fone(\waslambdaminus)/2)/\density(\fone(\waslambdaminus)/2>V_0$, and hence we
have $\vol M>V_0$  in this case.

We now turn to  the case in which Alternative (ii) of Corollary \ref{keyer corollary} holds.
We fix a point $p_0\in M$ with $\waslambdaminus\le\nextone(p_0)\le\log8$ and $\max_{p\in M}\shortone(p) \ge\fone(\nextone(p_0))$. We set $\lambda_0=\nextone(p_0)$, so that $\waslambdaminus\le\lambda_0\le\log8$ and
$\max_{p\in M}\shortone(p) \ge\fone(\lambda_0)$. 
We set $\mu_0=\scrM(\lambda_0)$, where the function $\scrM$ is given by the hypothesis of the present proposition. Since by hypothesis $\scrM$ takes its values in $I_k$, we have $\mu_0\in I_k$.
According to Condition (2) of the hypothesis, applied with $\lambda_0$ playing the role of $\lambda$, we may fix an $R_0\in(0,\f_3(\lambda_0)/2]$ 
such that
$\Wfive(R_0,\g_k(\mu_0),\fthree(\mu_0))
> V_0$. We have $\max_{p\in M}(\shortone(p)/2)\ge\fone(\lambda_0)/2\ge R_0$. 
Thus all the hypotheses of Corollary \ref{marseillaise} hold with $V_0$, $R_0$, and $\mu_0$ playing the respective roles of $V$,  $R$, and $\mu$, and with  $k$ given by the hypothesis of the present proposition.
It therefore follows from Corollary \ref{marseillaise} that either $\mu_0$ is  a Margulis number for $M$, or $\vol M\ge V_0$. Hence we may assume, for the remainder of the argument, that $\mu_0$ is a Margulis number for $M$.

We now apply Condition (3) of the hypothesis, with $\lambda_0\in[\waslambdaminus,\log8]$ and $\ell_M\in T$ playing the respective roles of $\lambda$ and $\wasell$. One of the alternatives (a), (b) of Condition (3) must hold with these choices of $\lambda$ and $\wasell$.

Consider the subcase in which Alternative (a) holds; in this subcase we may fix 
positive numbers $\alpha_0$ and $\delta_0$, with $\alpha_0\le\lambda_0$ and $\lambda_0/4<\delta_0< \ell_M$,
such that for every 
$D$ satisfying $D\ge\delta_0$ and $D>\fthree(\lambda_0)$,
we have 
\Equation\label{with these values}
\waschi_k(\alpha_0, \lambda_0,\delta_0,D,\mu_0)>V_0.
\EndEquation

We apply Lemma \ref{stoppeth two}, with $\lambda_0$, $\alpha_0$, $p_0$, $\mu_0$ and $\delta_0$  playing the respective roles of $\lambda$, $\alpha$, $p$, $\mu$ and $\delta$.
By definition we have $\lambda_0=\nextone(p_0)$. We have $0<\alpha_0\le\lambda_0$, and $\ell_M>\delta_0>\lambda_0/4\ge\alpha_0/4$; furthermore, $\mu_0$ is a Margulis number for $M$. Thus all the hypotheses of Lemma \ref{stoppeth two} hold, and that lemma gives a number
 $D_0$ satisfying  $D_0\ge\delta_0$ and $D_0>\fthree(\lambda_0)$,
such that 
\Equation\label{or those}
\vol M\ge
\waschi_k(\alpha_0, \lambda_0,\delta_0,D_0,\mu_0).
\EndEquation
Since 
$D_0\ge\delta_0$ and $D_0>\fthree(\lambda_0)$,
(\ref{with these values}) holds with $D=D_0$. Combining this with (\ref{or those}), we obtain the required inequality $\vol M>V_0$ in this subcase.

There remains the subcase in which Alternative (b) of Condition (3) of the present proposition holds, with   $\lambda_0\in[\waslambdaminus,\log8]$ and $\ell_M\in T$ playing the respective roles of $\lambda$ and $\wasell$. This means that 
$\ell_M<\log3$, and that there is a number $\mu_1\in(0,\mu_0]$ such that for every $h>0$ we have
\Equation\label{but on the other foot}
\waspsi_k(h,\ell_M, \mu_1)>V_0.
\EndEquation

By definition $\ell_M$ is the length of some closed geodesic in $M$. Since $0<\mu_1\le\mu_0$, and $\mu_0$ is a Margulis number for $M$, it follows from Proposition \ref{Margulis and s-one and s-two} that $\mu_1$ is also a Margulis number for $M$.
Since, in addition, $\pi_1(M)$ is  $k$-free and $\ell_M<\log3$,
it follows from Lemma \ref{stoppeth one} that 
there is  a positive number 
$h_0$ such that
$
\vol M
\ge \waspsi_k(h_0,\ell_M,\mu_1)$. But applying (\ref{but on the other foot}) with $h=h_0>0$, we obtain $\waspsi_k(h_0,\ell_M,\mu_1)>V_0$, and hence we have $\vol M>V_0$ in this subcase as well.
\EndProof

\Corollary\label{monster corollary}
Let $k\ge4$ be an integer, and let $V_0$ 
and $\wasdelta_0$ be  strictly positive numbers, with $\wasdelta_0<(\log8)/4$.
Suppose that $\waslambdaminus$ is  a number in $[\log7,\log8]$; that  $\wasdelta_1$ is a number satisfying $(\log8)/4<\wasdelta_1<\log3$; 
that $\scrD$ is a function defined on $[\waslambdaminus,\log8]$, and taking values in the interval $(\wasdelta_0,\wasdelta_1] \subset((\log8)/4,\log3]$; and that $\scrH$ is a (weakly) monotone increasing function defined on $(\wasdelta_0,\wasdelta_1] $ 
and taking values in the interval $I_k$
(so that, according to  \ref{ikea}, $\g_k(\scrH(\wasell))$ is defined whenever $\wasell\in (\wasdelta_0,\wasdelta_1]$).
Assume that the following conditions hold. 
\begin{enumerate}
\item $B(\fone(\waslambdaminus)/2)/\density(\fone(\waslambdaminus)/2)>V_0$.

\item
For each $\lambda\in[\waslambdaminus,\log8]$, there is an 
$R\in(0,\f_3(\lambda)/2]$ 
such that
$$\Wfive(R,\g_k(\scrH(\scrD(\lambda))),\fthree(\scrH(\scrD(\lambda))))
> V_0.$$

\item
For every $\lambda\in[\waslambdaminus,\log8]$, 
there is a  number $\alpha\in(0,\lambda]$ 
such that for every 
$D$ satisfying $D\ge\scrD(\lambda)$ and $D>\fthree(\lambda)$,
we have
$\waschi_k(\alpha, \lambda,\scrD(\lambda),D, \scrH(\scrD(\lambda)))>V_0$.

\item
For every $\wasell\in (\wasdelta_0,\wasdelta_1]$ 
there is a number $\mu\le \scrH(\wasell)$ such that for every $h>0$ we have
$\waspsi_k(h,\wasell, \mu)>V_0$.
\end{enumerate}

Then for every closed, orientable hyperbolic $3$-manifold $M$ such
that $\pi_1(M)$ is $k$-free and $\ell_M>\wasdelta_0$,  we have $\vol M\ge V_0$.

If in addition we assume that

\begin{itemize}
\item[{\it(5)}]
$0<\wasdelta_0<\min(0.7,\log(k-1)/2)$ (so that $\Wfour (k,\wasdelta_0, \scrH(\scrD(\log8)))$ is defined 
by \ref{reformulated very short lemma}),
 and $\Wfour (k,\wasdelta_0, \scrH(\scrD(\log8)))>V_0$, 

\end{itemize}
then  every closed, orientable hyperbolic $3$-manifold with $k$-free fundamental group has volume at least $V_0$.
\EndCorollary

\Proof
Suppose that
$k$, $V_0$, $\wasdelta_0$, $\waslambdaminus$, $\wasdelta_1$, $\scrD$, and $\scrH$ satisfy the hypotheses of the first assertion of the corollary. 
Set $T=(\wasdelta_0,\infty)$. Set $\scrM=\scrH\circ\scrD$, so that $\scrM$ is a well-defined function on $[\waslambdaminus,\log8]$, taking values in $I_k$. We will show that $k$, $V_0$, $T$, and $\scrM$ then satisfy Conditions (1)--(3) of Proposition \ref{latest monster}; the conclusion of the first assertion of the corollary 
will then follow immediately from Proposition \ref{latest monster}.

Conditions (1) and (2)  of Proposition \ref{latest monster} are immediate from Conditions (1) and (2), respectively, of the hypothesis of the present corollary.

To verify Condition (3)  of Proposition \ref{latest monster}, let $\lambda\in[\waslambdaminus,\log8]$ and $\wasell>\wasdelta_0$ be given. Consider first the case in which $\wasell>\scrD(\lambda)$. According to Condition (3) of the present corollary, we may fix a  number $\alpha\in(0,\lambda]$ 
such that for every 
$D$ satisfying $D\ge\scrD(\lambda)$ and $D>\fthree(\lambda)$ 
we have
\Equation\label{frumpin}
\waschi_k(\alpha, \lambda,\scrD(\lambda),D, \scrH(\scrD(\lambda)))>V_0.
\EndEquation
Set $\delta=\scrD(\lambda)$. Since $\scrD$ takes values in the interval $(\wasdelta_0,\wasdelta_1]\subset((\log8)/4,\wasdelta_1]$, and since we are in the case where $\wasell>\scrD(\lambda)=\delta$, we have $\lambda/4<\delta<\wasell$. 
As we have defined $\scrM$ to be $\scrH\circ\scrD$, the inequality (\ref{frumpin}), which holds for every 
$D$ satisfying $D\ge\scrD(\lambda)$ and $D>\fthree(\lambda)$, 
may be rewritten as
$\waschi_k(\alpha, \lambda,
\delta,
D,
\scrM(\lambda)) >V_0$.
Thus Alternative (a) of Condition (3) of Proposition \ref{latest monster} holds in this case.

Now consider the case in which $\wasell\le\scrD(\lambda)$. Since $\scrD$ takes values in $(\wasdelta_0,\wasdelta_1]\subset((\log8)/4,\wasdelta_1]$, we in particular have $\wasell\le\wasdelta_1$. This implies both that $\wasell<\log3$, and, since $\wasell$ was given to be greater than $\wasdelta_0$, that $\wasell$ lies in $(\wasdelta_0,\wasdelta_1]$, the domain of $\scrH$. 
According to Condition (4) of the present corollary, 
there is a number $\mu\le \scrH(\wasell)$ such that
 for every $h>0$, we have
$\waspsi_k(h,\wasell, \mu)>V_0$.
Since we are in
 the case  where $\wasell\le\scrD(\lambda)$, and since $\scrH$ is monotone increasing, we have $\mu\le \scrH(\wasell)\le\scrH(\scrD(\lambda))=\scrM(\lambda)$. 
Thus Alternative (b) of Condition (3) of Proposition \ref{latest monster} holds in this case. 
This completes the proof of the first assertion.

To prove the second assertion, suppose that the hypotheses of the first assertion hold, and that in addition 
Condition (5) holds. Let us set $\mu=\scrH(\scrD(\log8)))$, so that $\mu\in I_k$ by the general hypothesis 
of the present corollary; then Condition (5) asserts that
$0<\wasdelta_0<\min(0.7,\log(k-1)/2)$, so that $\Wfour (k,\wasdelta_0, \mu)$ is defined,
 and that $\Wfour (k,\wasdelta_0, \mu)>V_0$.

If we have $\ell_M>\wasdelta_0$,  it follows from the first assertion that $\vol M\ge V_0$. Hence we may assume that $\ell_M\le\wasdelta_0$. 

We apply Condition (2) with $\lambda=\log8$. Since $\fone(\log8)=\log5$, this gives an $R\in(0,(\log5)/2]$ such that $\Wfive(R,\g_k(\mu),\fthree(\mu))
> V_0$. Now since $\pi_1(M)$ is in particular $3$-free, it follows from \cite[Corollary 9.3]{acs-singular} that $M$ contains a hyperbolic ball of radius $(\log5)/2$; in particular, $M$ contains a hyperbolic ball of radius $R$.  It therefore follows from Corollary \ref{marseillaise} that either $\vol M\ge V_0$, or $\mu$ is a Margulis number for $M$. Thus we may assume that $\mu$ is a Margulis number for $M$. Now since $0<\wasdelta_0<\min(0.7,\log(k-1)/2)$, and since $\pi_1(M)$ is $k$-free  and $\ell_M\le\wasdelta_0$, it follows from \ref{reformulated very short lemma} that $\vol M$ is bounded below by $\Wfour(k,\wasdelta_0,\mu)$, which is in turn bounded below by $V_0$.
\EndProof

\section{
Lower bounds for the $4$-free and $5$-free cases
}\label{4 and 5}

In this section we will refer freely to  a number of functions defined in earlier sections, including $\Wone$ (Subsection \ref{reformulated surgeon general}), $\Wtwo$ (Subsection \ref{short geo def}), 
$\Wthree$ (Subsection \ref{is that you wthree?}), $\Wfive$ (Subsection \ref{reformulated better diameter from thick volume}), and $\Wfour$ (Subsection \ref{reformulated very short lemma}).

\NotationRemarks\label{ouchday}
We will need estimates for the functions $\omega$, $\barOmega$ and $\theta$ defined in \ref{all about Omega}.

If $x^-$, $x^+$, $y^- $, $y^+$ and $z^+$ are positive numbers with $x^-\le x^+$ and $y^-\le y^+$, we set
$$\omega^-(x^-,x^+,y^-,y^+,z^+)=(\coth x^+)(\coth y^+)-(\cosh z^+)(\cosech x^-)(\cosech y^-).$$ 
Comparing this definition with the definition of $\omega$ given in \ref{all about Omega}, we find that $\omega(x,y,z)\ge\omega^-(x^-,x^+,y^-,y^+,z^+)$ whenever $x^-\le x\le x^+$, $y^-\le y\le y^+$, and $0<z\le z^+$. Hence if we define
a function $\barOmega^+$ 
on
the same domain as $\omega^-$
by setting
$$\barOmega^+(x^-,x^+,y^-,y^+,z^+)=\arccos(\min(\max(\omega^-
(x^-,x^+,y^-,y^+,z^+),-1),1))\in[0,\pi]$$
for all $x^-,x^+,y^-,y^+,z^+ >0$, then we have 
\Equation\label{spooks}
\barOmega(x,y,z)\le\barOmega^+(x^-,x^+,y^-,y^+,z^+)
\EndEquation
 whenever 
$0<x^-\le x\le x^+$, $0<y^-\le y\le y^+ $,
and $0<z\le z^+$.

If $C^-$, $C^+$, $x^-$ and $x^+$ are numbers satisfying $0<C^-\le C^+$ and $0<x^-\le x^+$, we set
%We define a function $\theta^+$ on $(0,\infty)^4$ by setting
\Equation\label{tweedle}
\theta^+(C^-,C^+,x^-,x^+)=
\barOmega^+\bigg(\arccosh\bigg(\frac{\cosh x^-}{\cosh C^+}\bigg), \arccosh\bigg(\frac{\cosh x^+}{\cosh C^-}\bigg), x^-,x^+, C^+\bigg)
\EndEquation
if $x^->C^+$, 
and setting $\theta^+(C^-,C^+,x^-,x^+)=\pi/2$ if $x^-\le C^+$.
% and $x^+\ge C^-$,  and setting
%$\theta^- (C^-,C^+,x^-,x^+)=0$ if $x^+< C^-$. 

We then have 
\Equation\label{gadzooks}
\theta(C,x)\le\theta^+(C^-,C^+,x^-,x^+)
\EndEquation
whenever $C^-\le C\le C^+$ and $x^-\le x\le x^+$. Indeed, to establish (\ref{gadzooks}) in the case where  $x^-> C^+$, we note that  in this case we have $x>C$, so that $\theta(C,x)$ and $\theta^+(x^-,x^+,y^-,y^+,z^+)$ are defined by (\ref{threedle}) and (\ref{tweedle}) in this case. We have 
 $\arccosh((\cosh x^-)/(\cosh C^+))\le\arccosh((\cosh x)/(\cosh C))\le\arccosh((\cosh x^+)/(\cosh C^-))$, and 
(\ref{gadzooks}) is then simply the instance of  (\ref{spooks}) in which $\arccosh((\cosh x^-)/(\cosh C^+))$, $\arccosh((\cosh x)/(\cosh C))$, \linebreak
$\arccosh((\cosh x^+)/(\cosh C^-))$, $x^-$, $x$, $x^+$, $C$ and $C^+$ play the respective roles of $x^-$, $x$, $x^+$, $y^-$, $y$, $y^+$, $z$ and $z^+$.
In the case where $x^-\le C^+$, we have
$\theta^+(C^-,C^+,x^-,x^+)=\pi/2$,  and (\ref{gadzooks}) follows from the observation made in \ref{all about Omega} that $\theta(C,x)
\le
\pi/2$ for all positive $C$ and $x$.

If $\underline{a}=(C^-, C^+, a_{1,0}^-, a_{1,0}^+, a_{1,1}^-, a_{1,1}^+, a_{2,0}^-, a_{2,0}^+, a_{2,1}^- ,a_{2,1}^+)$ is a decuple of positive numbers such that $C^-\le C^+$ and
 $a_{i,0}^-\le a_{i,0}^+\le a_{i,1}^-\le a_{i,1}^+$ for $i=1,2$, we set 
$$A_1^+ (\underline{a})=\max_{(i, j)\in\{0,1\}\times\{0,1\}}\barOmega^+(a_{1,i}^-,a_{1,i}^+,a_{2,j}^-,a_{2,j}^+,C^+),$$
$$A_2^+(\underline{a})=\max_{(m,i)\in\{1,2\}\times\{0,1\}}\theta^+(C^-,C^+,a_{m,i}^-, a_{m,i}^+),$$ 
and 
$$
A^+(\underline{a})
=\max
(A_1^+(\underline{a}),
A_2^+(\underline{a})).
$$
Now given any $C\in[C^-,C^+]$, and given $a_{(m,i)}\in [a_{(m,i)}^-, a_{(m,i)}^+]$ for each ${(m,i)\in\{1,2\}\times\{0,1\}}$, it follows from 
(\ref{spooks}) and (\ref{gadzooks}), and from the definitions of $A_1$, $A_2$ and $A$ given in \ref{all about Omega}, that
 $A_t (C,a_{1,0},a_{1,1},a_{2,0},a_{2,1})\le A_t^+(\underline{a})$ for $t=1,2$, and hence 
\Equation\label{this is important}
A(C,a_{1,0},a_{1,1},a_{2,0},a_{2,1})\le A^+(\underline{a}).
\EndEquation

If $\delta$, $D^-$ and $D^+$ are numbers satisfying $0<\delta\le D^-\le D^+$, we set
$$\Lambda^+(\delta,D^-,D^+)=A^+( D^-,D^+,2\delta, 2\delta,\Phi_2(\delta,D^-),\Phi_2(\delta,D^+),3\delta, 3\delta, \Phi_3(\delta,D^-), \Phi_3(\delta,D^+)).$$
Given any $D\in[D^-,D^+]$, since $\Phi_2$ and $\Phi_3$ are monotone increasing in their second argument (cf. \ref{audrey junior}), we have $\Phi_n(\delta,D)\in[\Phi_n(\delta,D^-), \Phi_n(\delta,D^+)]$ for $n=2,3$. Hence (\ref{this is important}), with the definition of $\Lambda^+$ given above and the definition of $\Lambda$ given in \ref{whats Lambda}, implies:
\Equation\label{this is the point}
\Lambda(\delta,D)\le\Lambda^+(\delta,D^-,D^+)\text{ whenever }D\in[D^-,D^+].
\EndEquation
\EndNotationRemarks

\NotationRemarks\label{morstuf}
We set
$\frakW_0 =\{(\alpha,\delta,D^-,D^+)\in\RR^4: 0<\delta\le D^-\le D^+ \text{ and } \alpha>0\}$. 
Note that if $(\alpha,\delta,D^-,D^+)\in \frakW_0 $, then $\Phi_n(\delta, D^-)$ is defined for every $n \ge 1$. Furthermore, by 
\ref{appomatox} we have
$\Phi_2(\delta, D^-) \le 2D^-\le2D^+$, 
so that
$\Psi(D^+, 
\Phi_2(\delta, D^-))$ is defined (see \ref{moore}). For each 
$(\alpha,\delta,D^-,D^+)\in \frakW_0 $ we set
$$\Wminus(\alpha,\delta,D^-,D^+)=B(\alpha/2)-2\sigma(\alpha/2,D^-/2,\Phi_2(\delta,D^-)/2,\Psi(D^+, \Phi_2(\delta,D^-))).$$ 
Next we define $\frakW_0 '=\{(\alpha,\delta,D^-,D^+)\in \frakW_0 : 
D^+<\Phi_3(\delta,D^-)
\le \Phi_3(\delta,D^+)
<\alpha\}
$, and observe that by 
\ref{moore},
the quantities $\Theta(D^+/2, \alpha/2)$ and $\Theta(\Phi_3(\delta,D^-)/2 ,
\alpha/2)$ 
are defined for every $(\alpha,\delta,D^-,D^+)\in \frakW_0 '$. We define
$\frakW_0 ''$ to be the set of all points $(\alpha,\delta,D^-,D^+)\in \frakW_0 '$ such that
\Equation\label{wuddabout}
 \cos\bigg(\Theta\bigg(\frac {D^+}2,\frac\alpha2\bigg)-\Theta\bigg(\frac{\Phi_3(\delta,D^-)}2,\frac\alpha2\bigg)\bigg)<
\frac{\cosh D^-\cosh \Phi_3(\delta,D^-)-\cosh2D^+}{\sinh D^+\sinh \Phi_3(\delta,D^+)}.
\EndEquation
We then define the function $\Woneminus$ on $\frakW_0 $ by setting
$$\Woneminus(\alpha,\delta,D^-,D^+)=\W^-(\alpha,\delta,D^-,D^+)$$
if 
$(\alpha,\delta,D^-,D^+)\in\frakW_0 ''$, and setting
\Equation\label{the other case}
\begin{aligned}
\Woneminus(\alpha,\delta,D^-,D^+)&=
\W^-(\alpha,\delta,D^-,D^+)-2\kappa\bigg(\frac\alpha2,\frac{\Phi_3(\delta,D^-)
}2 \bigg)\\
&+2\iota\bigg(\frac\alpha2,\frac{\Phi_2(\delta,
D^+
)}2,\frac{\Phi_3(\delta,
D^+
)}2,\Lambda^+(\delta,D^-,D^+)\bigg)
\end{aligned}
\EndEquation
if 
$(\alpha,\delta,D^-,D^+)\in\frakW_0 -\frakW_0 ''$.

\Remarks\label{hbd}
Note that for any $(\alpha,\delta,D^-,D^+)\in \frakW_0$ and any $D\in[D^-,D^+]$, we have $(\alpha,\delta,D)\in\frakW$. It follows from the fact that $\Phi_2$ is increasing in its second argument  (see \ref{audrey junior}) that if $(\alpha,\delta,D^-,D^+)\in \frakW_0'$, then $(\alpha,\delta,D)\in \frakW'$ for any $D\in[D^-,D^+]$. 

If $(\alpha,\delta,D^-,D^+)\in \frakW_0''$, then $(\alpha,\delta,D)\in \frakW''$ for any $D\in[D^-,D^+]$. To see this, since we have already shown that $(\alpha,\delta,D)\in \frakW'$, it suffices to observe that the left- and right-hand sides of (\ref{erev birthday}) are respectively bounded above and below by the left- and right-hand sides of (\ref{wuddabout}); this observation in turn follows from the fact that $\Phi_3$ is  increasing in its second argument  (see \ref{audrey junior}), and the fact that $\Theta$ is decreasing in its first argument (see \ref{moore}).
\EndRemarks

\EndNotationRemarks

\Lemma\label{woney}
Let $\alpha$, $\delta$,  $D^-$ and $D^+$ be positive numbers such that $\delta\le D^-\le D^+$. Then for every $D$ with $D^-\le D\le D^+$, we have 
$\Wone(\alpha,\delta,D)\ge\Woneminus(\alpha,\delta,D^-,D^+)$. 
\EndLemma

\Proof
First note that since $\sigma$ is decreasing in its second and third arguments and increasing in its fourth argument (see \ref{what to call it}), while $\Phi_2$ and $\Phi_3$ are increasing in their second argument (see \ref{audrey junior})  and 
$\Psi$ is  increasing in its first argument and decreasing in its second (see \ref{moore}), it follows from the definition of $\W$ given in \ref{reformulated surgeon general} and the definition of $\Wminus$ given in \ref{morstuf} that
\Equation\label{almost like praying}
\W(\alpha,\delta,D)\ge\Wminus(\alpha,\delta,D^-,D^+).
\EndEquation
Likewise, since $\kappa$ is decreasing in its second argument (\ref{what to call it}), we find, again using that $\Phi_3$ is increasing in its second argument, that
\Equation\label{less noise}
\kappa\bigg(\frac\alpha2,\frac{\Phi_3(\delta,D)}2 \bigg)\le
\kappa\bigg(\frac\alpha2,\frac{\Phi_3(\delta,D^-)}2 \bigg).
\EndEquation
Next note that, according to (\ref{this is the point}), we have
$\Lambda(\delta,D)\le\Lambda^+(\delta,D^-,D^+)$. Using that $\iota$ is decreasing in its second, third and fourth arguments (see  \ref{what to call it}), and again using that $\Phi_2$ and $\Phi_3$ are increasing in their second argument, we deduce that
\Equation\label{no the other one}
\iota\bigg(\frac\alpha2,\frac{\Phi_2(\delta,D)}2,\frac{\Phi_3(\delta,D)}2,\Lambda(\delta,D)\bigg)\ge
\iota\bigg(\frac\alpha2,\frac{\Phi_2(\delta,
D^+
)}2,\frac{\Phi_3(\delta,
D^+
)}2,\Lambda^+(\delta,D^-,D^+)\bigg).
\EndEquation

We now distinguish two cases. The first is the case in which  $(\alpha,\delta,D^-,D^+)\in \frakW_0''$. It then follows from Remark \ref{hbd} that $(\alpha,\delta,D)\in \frakW''$. Hence, according to the definitions given in \ref{reformulated surgeon general}  and \ref{morstuf}, we have 
$\Woneminus(\alpha,\delta,D^-,D^+)=\W^-(\alpha,\delta,D^-,D^+)$ and 
$\Wone(\alpha,\delta,D)=\W(\alpha,\delta,D)$ in this case, and the conclusion of the lemma follows from
(\ref{almost like praying}).

Now consider the case in which  $(\alpha,\delta,D^-,D^+)\notin \frakW_0''$, so that $\Woneminus (\alpha,\delta,D^-,D^+)$ is given by (\ref{the other case}). According to Remark \ref{cases and cases}, $\Wone (\alpha,\delta,D)$ satisfies (\ref{inequality version}). But it follows from (\ref{almost like praying}), (\ref{less noise}) and (\ref{no the other one}) that the right-hand side of (\ref{inequality version}) is bounded below by the right-hand side of (\ref{the other case}). Hence the conclusion of the lemma is true in this case as well.
\EndProof

\Notation\label{chiless def}
Suppose that $k>2$ is an integer, that $\alpha$, $\lambda$, $\delta$ and $\mu$ are positive  numbers such that $\alpha\le\lambda$ and  $\delta>\f_1(\lambda)$, 
and that $E^-$ and $E^+$ are real numbers with $E^+>E^-\ge0$.
Since $\alpha>0$ and $E^++\delta>E^-+\delta\ge\delta>0$, the quantity $\Wone(\alpha,\delta,E^-+\delta,E^++\delta)$ is defined. Since $E^++\delta>\delta>\fthree(\alpha)$, 
we have $Q(\lambda)+Q(E^++\delta)<1/2$ by \ref{Q-def}, 
so that $\xi_{k-2}(Q(\lambda)+Q(E^++\delta))$ is defined and strictly positive. We set
$$\waschi^-_k(\alpha,\lambda,\delta,E^-,E^+,\mu)=\Woneminus(\alpha,\delta,E^-+\delta, E^++\delta)+\Wthree(\xi_{k-2}(Q(\lambda)+Q(E^++\delta)), \frac\lambda2,\mu).$$
\EndNotation

\Lemma\label{now a cow}
Let $k>2$ is an integer, let $\lambda^-$, $\lambda^+$, $E^-$ and $E^+$ are real numbers with $0<\lambda^-<\lambda^+$ and $0\le E^-<E^+$, let $\delta$ be a number greater than $\fthree(\lambda^-)$, and let $\mu$ be a positive number. Note that since $\fthree$ is monotone decreasing by 
\ref{Q-def}, 
we have $\delta>\fthree(\lambda^+)$, and hence $\waschi^-_k(\lambda^-,\lambda^+,\delta,E^-,E^+,\mu)$ is defined in view of  \ref{chiless def}. Let $E\in(E^-,E^+]$ and $\lambda\in[\lambda^-,\lambda^+]$ be given; note that  $E+\delta>\delta>\fthree(\lambda)$ (again by the monotonicity of $\fthree$), so that 
$\waschi_k(\lambda^-, \lambda,\delta,E+\delta,\mu)$ 
is defined in view of \ref{chi def}. 
Then we have 
$
\waschi_k(\lambda^-, \lambda,\delta,E+\delta,\mu)
\ge
\waschi^-_k(\lambda^-,\lambda^+,\delta,E^-,E^+,\mu)$.
\EndLemma

\Proof
Since $\delta\le E^-+\delta< E+\delta\le E^++\delta$, Lemma \ref{woney} gives
\Equation\label{watch it go by}
\Wone(\lambda^-,\delta,E+\delta)\ge\Woneminus(\lambda^-,\delta,E^-+\delta,E^++\delta).
\EndEquation
On the other hand, since $\Wthree$ is   increasing in its first argument and decreasing in its second
%and
by \ref{is that you wthree?}, and $\xi_{k-2}$ is increasing and $Q$ is decreasing (see \ref{reformulated big radius corollary from cusp} and \ref{Q-def}), we have 
\Equation\label{god still forbid}
\Wthree\bigg(\xi_{k-2}(Q(\lambda)
+Q(E+\delta)),\frac\lambda2,\mu\bigg)
\ge \Wthree\bigg(\xi_{k-2}(Q(\lambda^+)+Q(E^++\delta)),\frac{\lambda^+}2,\mu\bigg).
\EndEquation
Adding the inequalities (\ref{watch it go by}) and (\ref{god still forbid}), and using the definitions of the functions $\waschi_k$ and $\waschi^-_k$, we obtain the asserted inequality $
\waschi_k(\lambda^-, \lambda,\delta,E+\delta,\mu)
\ge
\waschi^-_k(\lambda^-,\lambda^+,\delta,E^-,E^+,\mu)$.
\EndProof

%\alpha

\Lemma\label{other woney}
Let $\maybecalLminus $, $\maybecalLplus$ and $\delta$ be positive numbers such that $\delta<\maybecalLminus \le\maybecalLplus$. Then for every $\alpha\in[\maybecalLminus ,\maybecalLplus]$ and every $D\in[\maybecalLplus,\infty)$, we have 
$$\Wone(\alpha,\delta,D)\ge B\bigg(\frac{\maybecalLminus }2\bigg)-2\kappa\bigg(\frac{\maybecalLplus}2, \frac{\Phi_2(\delta,\maybecalLplus)}2\bigg) -2\kappa\bigg(\frac{\maybecalLplus}2, \frac{\Phi_3(\delta,\maybecalLplus)}2\bigg).$$
\EndLemma

\Proof
It follows from Remark \ref{cases and cases} and the definition of the function $\W$ that
%$$
\Equation\label{dumber inequality version}
\begin{aligned}
\Wone(\alpha,\delta,D)&\ge
\W(\alpha,\delta,D)-2\kappa\bigg(\frac\alpha2,\frac{\Phi_3(\delta,D)}2 \bigg)\\
&=B(\alpha/2)-2\sigma\bigg(\frac\alpha2,\frac D2,\frac{\Phi_2(\delta,D)}2,\Psi(D, \Phi_2(\delta,D))\bigg)
-2\kappa\bigg(\frac\alpha2,\frac{\Phi_3(\delta,D)}2 \bigg).
\end{aligned}
%$$
\EndEquation
Since $D/2\ge\maybecalLplus/2\ge\alpha/2$, in the notation of \ref{what to call it} we have 
$K_N(\zeta , D/2)=\emptyset$, where
 $N\subset \HH^3$ is an open ball of radius $\alpha/2$ and
$\zeta$ is a point of $\partial\overline N$. The definitions in \ref{what to call it} then give  
$\sigma(\alpha/2, D/2,\Phi_2(\delta,D)/2,\Psi(D, \Phi_2(\delta,D)))=\kappa(\alpha/2, \Phi_2(\delta,D)/2)$. Thus
(\ref{dumber inequality version}) becomes
\Equation\label{inudda woids}
\Wone(\alpha,\delta,D)\ge
B(\alpha/2)-2\kappa\bigg(\frac\alpha2,\frac{\Phi_2(\delta,D)}2\bigg)
-2\kappa\bigg(\frac\alpha2,\frac{\Phi_3(\delta,D)}2 \bigg).
\EndEquation

Since $\alpha\ge\maybecalLminus $, we have $B(\alpha/2)\ge B(\maybecalLminus /2)$. On the other hand,
according to \ref{what to call it}, $\kappa$ is increasing in its first argument and decreasing in its second, while according to \ref{audrey junior}, $\Phi_2$ and $\Phi_3$ are increasing in their second argument; since $\alpha\le\maybecalLplus\le D$, it follows that for $n=2,3$ we have  $\kappa(\alpha/2,\Phi_n(\delta,D)/2)\le \kappa(\maybecalLplus/2,\Phi_n(\delta,\maybecalLplus)/2)$. The conclusion of the lemma therefore follows from (\ref{inudda woids}).
%$\kappa\bigg(\frac\lambda2,\frac{\Phi_3(\delta,D)}2 \bigg).
\EndProof

\Notation\label{hencricks}
Let  $k>2$ be an integer, 
%and let $h^-$, $h^+$, $l^-$, $l^+$ and $\mu$ be positive numbers such that  $h^-<h^+$ and $l^-<l^+$.
and let $h^-$, $h^+$, $l^-$, $l^+$ and $\mu$ be real numbers such that  $0\le h^-<h^+$ and $0<l^-<l^+$. 
 Since $\fthree$ is strictly  decreasing by  
\ref{Q-def},
we have $\fthree(l^-)+h^+>\fthree(l^-)>\fthree(l^+)$; hence by
\ref{Q-def},
we have
$Q(\fthree(l^-)+h^+)+Q(l^+)<1/2$, so that $\xi_{k-2}(Q(\fthree(l^-)+h^+)+Q(l^+))$ is defined.
We set
$$
\waspsi_k^-( h^-,h^+,l^-,l^+,\mu)=
\Wtwo(\fthree(l^+)+h^-,l^-)+
\Wthree\bigg(\xi_{k-2}(Q(\fthree(l^-)+h^+)+Q(l^+)),\frac{\fthree(l^-)+h^+}2,\mu
\bigg).
$$
\EndNotation

\Lemma\label{cranky-poo}
Let  $k>2$ be an integer, 
and let $h^-$, $h^+$, $l^-$, $l^+$ and $\mu$ be real numbers such that  $0\le h^-<h^+$ and $0<l^-<l^+$. 
%And $h\in[h^-,h^+]$ should be $h\in(h^-,h^+]$.
Then for every $h\in(h^-,h^+]$ and every $l\in[l^-,l^+]$, we have
$\waspsi_k( h,l,\mu)\ge\waspsi^-_k( h^-,h^+,l^-,l^+,\mu)$.
\EndLemma

\Proof
The function
$\xi_{k-2}$ is strictly increasing by \ref{reformulated big radius corollary from cusp}, 
while $\Wthree$ is increasing in its first argument and decreasing in its second (\ref{is that you wthree?}). According to
Lemma \ref{on account of},  the function $\Wtwo$
is monotonically increasing  in both its  arguments. Furthermore, $Q$ and $\fthree$ are strictly monotone decreasing by
\ref{Q-def}.
It follows that 
$$\Wtwo(\fthree(l)+h,l)\ge \Wtwo(\fthree(l^+)+h^-,l^-)$$
and that
$$\Wthree\bigg(\xi_{k-2}(Q(\fthree(l)+h)+Q(l)),\frac{\fthree(l)+h}2,\mu
\bigg)\ge
\Wthree\bigg(\xi_{k-2}(Q(\fthree(l^-)+h^+)+Q(l^+)),\frac{\fthree(l^-)+h^+}2,\mu
\bigg).
$$
Adding these two inequalities, and using the definitions of the functions $\waspsi_k$ and $\waspsi^-_k$, we obtain the asserted inequality $\waspsi_k( h,l,\mu)\ge\waspsi^-_k( h^-,h^+,l^-,l^+,\mu)$.
\EndProof

\Lemma\label{coffee}
Let $k\ge4$ be an integer, and let $V_0$  be a  strictly positive real numbers. Suppose that  $\wasdelta_0$ and $\wasdelta_1$ are numbers satisfying $\wasdelta_0< (\log8)/4<\wasdelta_1<\log3$; 
that $m$,  $n$, $p$ and $q$ are positive integers, and that $\lambda_0,\ldots,\lambda_m$, $\delta_0,\ldots,\delta_m$, $E_0,\ldots,E_n$, $h_0,\ldots,h_p$, and $\zeta_0,\ldots,\zeta_q$ are real numbers with $\log7\le\lambda_0<\cdots<\lambda_m=\log8$, $\wasdelta_1=\delta_0\ge\cdots\ge\delta_{m}\ge\wasdelta_0$, $0=E_0<\cdots<E_n=(\log8)-\delta_m$, $0=h_0<\cdots<h_p$, and $\wasdelta_0=\zeta_0<\cdots<\zeta_q=\wasdelta_1$; 
and that $\scrH$ is a (weakly) monotone increasing function defined on $(\wasdelta_0,\wasdelta_1] $ 
and taking values in the interval $I_k$
(so that, according to  \ref{ikea}, $\g_k(\scrH(\wasell))$ is defined whenever $\wasell\in (\wasdelta_0,\wasdelta_1]$).
Suppose that $\delta_m>\fthree(\lambda_0)$; since $\fthree$ is monotone decreasing, this gives that
$\delta_i\ge\delta_m>\fthree(\lambda_{0})>\fthree(\lambda_{i})$ for $i=1,\ldots,m$, so that
$\waschi^-_k(\lambda_{i-1},\lambda_i,\delta_i,E_{j-1},E_j,\scrH(\delta_i))$ 
is defined for $i=1,\ldots,m$ and for $j=1,\ldots,n$.
Note also that since $\wasdelta_0< (\log8)/4<\min(0.7,(\log3)/2)\le\min(0.7,\log(k-1)/2)$, the quantity $\Wfour (k,\wasdelta_0, \scrH(\delta_m))$ is defined 
by \ref{reformulated very short lemma}.

Assume that the following conditions hold.

\begin{itemize}

\item[(1)] $B(\fone(\lambda_0)/2)/\density(\fone(\lambda_0)/2)>V_0$.

\item[(2)]
For $i=1,\ldots,m$ we have 
$$\Wfive(\fone(\lambda_i)/2,\g_k(\scrH(\delta_i)),\fthree(\scrH(\delta_i)))
> V_0.$$

\item[(3a)]
For $i=1,\ldots,m$ and for $j=1,\ldots,n$, we have 
$\waschi^-_k(\lambda_{i-1},\lambda_i,\delta_i,E_{j-1},E_j,\scrH(\delta_i))>V_0$.

\item[(3b)]
$B(\lambda_0/2)-2\kappa((\log8)/2,\delta_m)-2\kappa((\log8)/2,3\delta_m/2)>V_0$.

\item[(4a)]
For $s=1,\ldots,q$ and for $t=1,\ldots,p$, we have
$\waspsi^-_k( h_{t-1},h_t,\zeta_{s-1},\zeta_s,\scrH(
\zeta_{s-1}))
>V_0$.

\item[(4b)]
For $s=1,\ldots,q$ we have
$\Wtwo(\fthree(\zeta_s)+h_p,\zeta_{s-1})
>V_0$.

\item[(5)]
$\Wfour (k,\wasdelta_0, \scrH(\delta_m))>V_0$. 
\end{itemize}

Then every closed, orientable hyperbolic $3$-manifold with $k$-free fundamental group has volume at least $V_0$.
\EndLemma

\Proof
Set $\waslambdaminus=\lambda_0$. Define a step function $\scrD$ on $[\waslambdaminus,\log8]$ by setting $\scrD(\waslambdaminus)=\delta_1$, and $\scrD(\lambda)=\delta_i$ whenever $i\in\{1,\ldots,m\}$ and $\lambda_{i-1}<\lambda\le\lambda_i$. We will show that Conditions (1)--(5) of the present lemma imply Conditions (1)--(5) of Corollary \ref{monster corollary}, with these definitions of $\waslambdaminus$, and with $k$, $V_0$, $\wasdelta_0$, $\wasdelta_1$ and $\scrH$ given by the hypothesis of the present lemma. The conclusion will then follow immediately from Corollary \ref{monster corollary}.

Since $\waslambdaminus=\lambda_0$ and $\scrD(\log8)=\delta_m$, Conditions (1) and (5) of Corollary \ref{monster corollary} are respectively identical with Conditions (1) and (5) of the present lemma.

To verify Condition (2) of Corollary \ref{monster corollary}, note that if $\lambda\in[\waslambdaminus,\log8]$ is given, then for some $i\in\{1,\ldots,m\}$ we have either $\lambda_{i-1}<\lambda\le\lambda_i$, or $i=1$ and $\lambda=\waslambdaminus$. Then $\scrD(\lambda)=\delta_i$. If we set $R=\fone(\lambda_i)/2$, we have $R\in(0,\fone(\lambda)/2]$ because $\fone$ is monotone decreasing (see \ref{Q-def}). Using  Condition (2) of the present lemma, we find $\Wfive(R,\g_k(\scrH(\scrD(\lambda))),\fthree(\scrH(\scrD(\lambda))))=
\Wfive(\fone(\lambda_i)/2,\g_k(\scrH(\delta_i)),\fthree(\scrH(\delta_i)))
> V_0$, as required for Condition (2) of Corollary \ref{monster corollary}.

To verify Condition (3) of Corollary \ref{monster corollary}, we again consider an arbitrary $\lambda\in[\waslambdaminus,\log8]$, and again choose an $i\in\{1,\ldots,m\}$ such that either $\lambda_{i-1}<\lambda\le\lambda_i$, or $i=1$ and $\lambda=\waslambdaminus$, so that we again have  $\scrD(\lambda)=\delta_i$. Set $\alpha=\lambda_{i-1}\in
(0,\lambda]$; we will verify Condition (3) of Corollary \ref{monster corollary} with this choice of $\alpha$. Since $\fthree$ is monotone decreasing, we have $\scrD(\lambda)=\delta_i\ge\delta_m>\fthree(\lambda_0))\ge\fthree(\lambda)$. Hence the pair of inequalities $D\ge\scrD(\lambda)$ and $D>\fthree(\lambda)$ which appear in Condition (3) of Corollary \ref{monster corollary} (and are needed to guarantee that 
$\waschi_k(\alpha, \lambda,\scrD(\lambda),D, \scrH(\scrD(\lambda)))$ 
is defined) are equivalent to the single inequality $D\ge\delta_i$. 
Thus we are required to show that for every $D\ge\delta_i$ we have
$\waschi_k(\lambda_{i-1}, \lambda,\delta_i,D, \scrH(\delta_i))>V_0$.

Consider first the case in which $D\le\log8$. Since $D\ge\delta_i$, we may write $D=E+\delta_i$ for some $E$ with $0\le E\le(\log8)-\delta_i\le(\log8)-\delta_m=E_n$. We then have $E_{j-1}\le E\le E_j$ for some $j\in\{1,\ldots,n\}$. Applying Lemma \ref{now a cow}, with $\lambda_{i-1}$, $\lambda_{i}$, $E_{j-1}$, $E_j$, $\delta_i$ and $\scrH(\delta_i)$ playing the respective roles of $\lambda^-$, $\lambda^+$, $E^-$, $E^+$, $\delta$ and $\mu$ in that lemma, and applying Condition (3a) of the present lemma, we find
$$\begin{aligned}
\waschi_k(\lambda_{i-1}, \lambda,\delta_i,D,
\scrH(\delta_i))&=\waschi_k(\lambda_{i-1}, \lambda,\delta_i,E+\delta_i,\scrH(\delta_i))\\&\ge
\waschi^-_k(\lambda_{i-1},\lambda_i,\delta_i,E_{j-1},E_j,\scrH(\delta_i))>V_0,
\end{aligned}
$$
as required.

Now consider the case in which $D>\log8$. The definition of $\waschi_k$ (see \ref{chi def}) immediately implies that
$\waschi_k(\lambda_{i-1}, \lambda,\delta_i,D,\scrH(\delta_i))\ge\Wone(\lambda_{i-1},\delta_i,D)$. 

We have $D>\log8=\lambda_m>\lambda_{i-1}$, and we may therefore apply Lemma 
\ref{other woney}, with $\waslambdaminus$, $\log8$, $\delta_i$ and $\lambda_{i-1}$ playing the respective roles of $\maybecalLminus$, $\maybecalLplus$, $\delta$ and $\alpha$ in that lemma, and with $D$ given as above, to deduce that
\Equation\label{other-dfoo}
\Wone(\lambda_{i-1},\delta_i,D)\ge B\bigg(\frac{\waslambdaminus}2\bigg)-2\kappa\bigg(\frac{\log8}2, \frac{\Phi_2(\delta_i,\log8)}2\bigg) -2\kappa\bigg(\frac{\log8}2, \frac{\Phi_3(\delta_i,\log8)}2\bigg).
\EndEquation
Recalling that $\waslambdaminus=\lambda_0$, that $\kappa$ is monotone decreasing in its second variable by 
\ref{what to call it}, 
and that $\Phi_2(\delta_i,\log8)\ge2\delta_i\ge2\delta_m$ and $\Phi_3(\delta_i,\log8)\ge3\delta_i\ge3\delta_m$ by 
\ref{wi-fi}, 
we find that the right hand side of (\ref{other-dfoo}) is bounded below by $B(\lambda_0/2)-2\kappa((\log8)/2,\delta_m)-2\kappa((\log8)/2,3\delta_m/2)$, which by Condition (3b) of the present lemma is in turn greater than $V_0$. This completes the verification of Condition (3) of Corollary \ref{monster corollary}.

It remains to verify Condition (4) of Corollary \ref{monster corollary}. Let $\wasell\in (\wasdelta_0,\wasdelta_1]$ be given. We may fix an index $s\in\{1,\ldots,q\}$ such that $\zeta_{s-1}<\wasell\le\zeta_s$.  
We set $\mu=\scrH(\delta_{s-1})$. 
Since $\scrH$ is monotone increasing and $\delta_{s-1}<\wasell$, we have $\mu\le\scrH(\wasell)$.
We must show that for every $h>0$ we have
$\waspsi_k(h,\wasell, \mu)>V_0$.

We first consider the case in which $h\le h_p$. In this case, there is an index $t\in\{1,\ldots,p\}$ such that $h_{t-1}\le h\le h_t$. We apply Lemma \ref{cranky-poo}, with $\zeta_{s-1}$, $\zeta_{s}$, $h_{t-1}$ and $h_{t}$ playing the respective roles of $l^-$, $l^+$, $h^-$ and $h^+$, and with $\mu$ defined as above. This shows that 
$\waspsi_k( h,\wasell,\mu)\ge\waspsi^-_k( h_{t-1},h_t,\zeta_{s-1},\zeta_s,\mu)=\waspsi^-_k( h_{t-1},h_t,\zeta_{s-1},\zeta_s,\scrH(\zeta_{s-1}))$. According to Condition (4a) of the present lemma, the right hand side of the latter inequality is strictly bounded below by $V_0$. Hence the required inequality $\waspsi_k( h,\wasell,\mu)>V_0$ holds in this case.

There remains the case in which $h> h_p$. In this case, we observe that it is immediate from the definition of $\waspsi_k( h,\wasell,\mu)$ (see \ref{before stoppeth one}) that 
$\waspsi_k( h,\wasell,\mu)\ge\Wtwo(\fthree(\wasell)+h,\wasell)$.
The function $\Wtwo$ is monotonically increasing in each of its  arguments by Lemma \ref {on account of}, and $\fthree$ is monotone decreasing (see \ref{Q-def}). Since $h> h_p$ and $\zeta_{s-1}\le\wasell\le\zeta_s$, it follows that 
$\waspsi_k( h,\wasell,\mu)\ge\Wtwo(\fthree(\zeta_s)+h_{p},\zeta_{s-1})$.
Since Condition (4b) of the present lemma asserts that the right hand side of the latter inequality is strictly bounded below by
$V_0$, the required inequality $\waspsi_k( h,\wasell,\mu)>V_0$ is established in this case as well.
\EndProof

\Number\label{calculations}
In the applications of Lemma \ref{coffee} to the proofs of Propositions \ref{four-free case} and \ref{five-free case}, it will be necessary to calculate a finite number of values of the functions which appear as the left-hand sides of the inequalities in Conditions (1)---(5) of Lemma \ref{coffee}. These functions are defined in terms of the functions $\kappa$ and $\iota$ (see Subsection \ref{what to call it}), the function $\density$ (see Subsection \ref{Boroczky number}), and standard elementary functions. The calculation of  values of  $\kappa$ and $\iota$  is done by the method described in Subsection 14.7 of \cite{fourfree}. The calculation of values of the function $\density$ is done using Lobachevksy functions, as in \cite{fourfree}.

The proof of Proposition
\ref{good homology} will be based on Proposition \ref{latest monster} rather than on Lemma \ref{coffee}, but it will involve a closely related method that also involves calculating a finite number of values of the functions $\kappa$, $\iota$ and $\density$. These calculations will also be done by the methods described above.
\EndNumber

\Proposition\label{four-free case}
Every closed, orientable hyperbolic $3$-manifold with $4$-free fundamental group has volume greater than $3.57$.
\EndProposition

\Proof
We apply Lemma \ref{coffee}, taking 
$k=4$ and $V_0=3.570002$. We set $\wasdelta_0=0.033$ and $\wasdelta_1=0.5505$; $m=100$, and $\lambda_i=(\log 8)-0.0335 + 0.000335i$ for $0\le i\le100$. We define the $\delta_i$, for $i=1,\ldots,100$, by setting $\delta_i=0.5505$ for $0\le i\le39$; $\delta_i=0.545$ for $40\le i\le66$; $\delta_i=0.544$ for $67\le i\le73$; $\delta_i=0.543$ for $74\le i\le82$; $\delta_i=0.5423$ for $83\le i\le88$; $\delta_i=0.542$ for $89\le i\le91$; $\delta_i=0.5418$ for $92\le i\le94$; and setting $\delta_{95}=\delta_{96}=0.5416$, $\delta_{97}=0.5415$, $\delta_{98}=0.54145$, $\delta_{99}=0.54138$, and $\delta_{100}=0.5413$. We take $n=2000$ and set $E_j=
 ((\log8)-0.5413)j/2000$ for $j=0,\ldots,2000$. 
We set $p=1000$, and $h_t=t/2000$ for $t=0,\ldots,1000$. We take $q=1290$, and we set $\zeta_s=0.033+0.00255s$ for $0\le s\le190$; $\zeta_s=0.5175+0.00013(s-190)$ for $191\le s\le290$; and $\zeta_s=0.5305+0.00002(s-290)$ for $291\le s\le1290$. We define the (weakly) monotone increasing function $\scrH$ on $(0.033,0.5505]$ by $\scrH(l)=\max(1.1253,l+0.584)$; its range is $[1.1253,1.1345]\subset I_4=(\log3,\log(17/3))$. Note that we have $\wasdelta_0< (\log8)/4<\wasdelta_1<\log3$ and $\delta_{100}=0.5413>\f_1(\lambda_0)=0.46\ldots$, as required for the application of Lemma \ref{coffee}.
The verification of Conditions (1)--(5) of Lemma \ref{coffee} consists of calculating a 
total of $1,491,393$ function values
(using the method referred to in Subsection \ref{calculations})
and observing that they are all bounded below by $V_0$. The smallest of these is $3.5700023\ldots$, and is given by the left-hand side of the inequality in Condition (2) of Lemma \ref{coffee} when $i=82$. As verifying Conditions (3a) and (4a) involves by far the largest number of inequalities, we also record that the smallest value of the left-hand side of any of the inequalities needed for Condition (3a) is $3.5714082\ldots$, and is achieved when $i=100$ and $j=512$, and that the smallest value of the left-hand side of any of the inequalities needed for Condition (4a) is
$3.57006\ldots$, and is achieved when $s=831$ and $t=298$.

It now follows from Lemma \ref{coffee} that the volume of every closed, orientable hyperbolic $3$-manifold with $4$-free fundamental group is at least $V_0>3.57$.
\EndProof

\Proposition\label{five-free case}
Every closed, orientable hyperbolic $3$-manifold with $5$-free fundamental group has volume greater than $3.77$.
\EndProposition

\Proof
We apply Lemma \ref{coffee}, taking $k=5$ and 
$V_0=3.7700008$. 
We set $\epsilon_0=0.033$, $\epsilon_1=0.5643$, and $m=121$. To define $\lambda_i$ for $0\le i\le121$, we first  define $\iota_i$ by setting $\iota_i=i$ when $0\le i\le49$; $\iota_i=(i+49)/2$ for $49\le i\le90$; and $\iota_i=i-21$ for $91\le i\le 121$. We then set $\lambda_i=(\log8)-.0457+.000457\iota_i$ for $0\le121$. We define $\delta_i$ for $1\le i\le121$ by setting $\delta_i=0.5642$ for $1\le i\le24$; $\delta_i=0.5593$ for $25\le i\le32$; $\delta_i=0.5573$ for $33\le i\le37$; $\delta_i=0.5559$ for $38\le i\le41$;  $\delta_{42}=0.5556$; $\delta_{43}=\delta_{44}=0.5548$; 
$\delta_i=0.55414-0.00019(i-46)$ for $46\le i\le48$; 
$\delta_{49}=0.553525$;
$\delta_{50}=\delta_{51}=0.55325$; 
$\delta_i=0.55312-.00014(i-52)$ for $52\le i\le54$;
$\delta_{55}=0.55273$;
$\delta_{56}=0.5526$;
$\delta_{57}=0.55248$;
$\delta_i=0.55235-0.00012(i-58)$ for $58\le i\le68$; 
$\delta_{69}=0.55109$;
$\delta_i=0.5509-.00011(i-70)$ for $70\le i\le73$;
$\delta_i=0.55058-.0001(i-74)$ for $74\le i\le91$;
$\delta_i=0.54885-.00016(i-92)$ for $92\le i\le98$;
$\delta_{99}=\delta_{100}=0.554765$; 
$\delta_{101}=\delta_{102}=0.554738$;
$\delta_i=0.54698$ for $103\le i\le105$;
$\delta_i=0.54625$ for $106\le i\le111$;
and $\delta_i=0.5452$ for $112\le i\le121$.

We take $n=4000$ and set $E_j=((\log8)-.5452)j/4000$
% ((\log8)-0.5413)j/2000$ for 
$j=0,\ldots,4000$. 
We set $p=1080$, and $h_t=t/2000$ for $t=0,\ldots,1080$. 
We take $q=1980$, and we set $\zeta_s=0.033+0.00255s$ for $0\le s\le190$; $\zeta_s=0.5175+0.00013(s-190)$ for $191\le s\le290$; and $\zeta_s=0.5305+0.00002(s-290)$ for $291\le s\le1980$. We define the (weakly) monotone increasing function $\scrH$ on $(0.033,0.5643]$ by $\scrH(l)=\max(1.1319,l+0.5867)$; its range is $[1.1319,1.151]\subset I_5=(\log3,\log9)$. 
Note that we have $\wasdelta_0< (\log8)/4<\wasdelta_1<\log3$ and $\delta_{121}=0.5452>\f_1(\lambda_0)=0.47\ldots$, as required for the application of Lemma \ref{coffee}.
The verification of Conditions (1)--(5) of Lemma \ref{coffee} consists of calculating a finite number of function values
(using the method referred to in Subsection \ref{calculations})
and observing that they are all bounded below by $V_0$. 
The smallest of these is $3.77000089\ldots$, and is given by the left-hand side of the inequality in Condition (2) of Lemma \ref{coffee} when $i=24$. As verifying Conditions (3a) and (4a) involves by far the largest number of inequalities, we also record that the smallest value of the left-hand side of any of the inequalities needed for Condition (3a) is $3.770016\ldots$, and is achieved when $i=68$ and $j=23$, and that the smallest value of the left-hand side of any of the inequalities needed for Condition (4a) is $3.770012\ldots$, and is achieved when $s=1026$ and $t=368$.

It now follows from Lemma \ref{coffee} that the volume of every closed, orientable hyperbolic $3$-manifold with $5$-free fundamental group is at least $V_0>3.77$.
\EndProof

\section{Volume and homology}\label{homology section}

Let us recall some standard definitions in $3$-manifold theory. By a {\it surface} in a closed $3$-manifold $M$ we will mean a connected $2$-dimensional submanifold $S$ of $M$ which is tame, i.e. is smooth with respect to some smooth structure on $M$. (Some authors use the word ``embedded'' to emphasize that $S$ is a submanifold rather than a $2$-manifold equipped with an immersion in $M$.) An orientable $3$-manifold $M$ is said to be {\it irreducible} if 
$M$ is connected and
every surface in $M$ which is homeomorphic to $S^2$ is the boundary of a $3$-ball in $M$. An {\it incompressible surface} in an irreducible, closed, orientable $3$-manifold $M$ is an orientable surface $S$ which is not a $2$-sphere, but has the property that the inclusion homomorphism $\pi_1(S)\to\pi_1(M)$ is injective. If $M$ is hyperbolic, then $M$ is irreducible and every incompressible surface in $M$ has genus at least $2$.

Following \cite{CDS}, we will say that $M$ is {\it $(g,h)$-small,} where $g$ and $h$ are given positive integers, if every incompressible surface in $M$ has genus at least $h$, and every separating incompressible surface in $M$ has genus at least $g$.

\Lemma\label{small or tall}
Let $g$ be a positive integer, and $M$ be a closed, orientable hyperbolic $3$-manifold that contains an incompressible (embedded) surface of genus $g$. Then there exist an integer $g'\ge2$ and an incompressible surface $T\subset M$ of genus $g'$, such that either 
\Alternatives
\item $g'\le2g-2$, the surface $T$ separates $M$, and $M$ is $(g',g'/2+1)$-small, or
\item $g'\le g$, the surface $T$ does not separate $M$, and $M$ is $(2g'-1,g')$-small.
\EndAlternatives
\EndLemma

\Proof
Let $g_0\le g$ denote the smallest positive integer which occurs as the genus of an incompressible surface in $M$. Since $M$ is hyperbolic we have $g_0\ge2$. Consider first the case in which $M$ contains a separating incompressible surface of genus $g_0$; we choose such a surface and denote it by $T_0$. According to our choice of $g_0$, every  incompressible surface in $M$ has genus at least $g_0$; and since $g_0\ge2$ we have $g_0/2+1\le g_0$. Thus $M$ is $(g_0,g_0/2+1)$-small in this case. Furthermore, since $2\le g_0\le g$, we have in particular that $g_0\le2g-2$. Thus, in this case, Alternative (i) of the conclusion holds with $T=T_0$ and $g'=g_0$.

Now consider the case in which $M$ contains no separating incompressible surface of genus $g_0$. According to our choice of $g_0$, there is an incompressible surface $T_1\subset M$ whose genus is $g_0$, and in this case $T_1$ must be non-separating. 

We distinguish two subcases. The first is the subcase in which every separating  incompressible surface in $M$ has genus at least $2g_0-1$. Our choice of $g_0$ guarantees that every incompressible surface in $M$ has genus at least $g_0$. Thus, according to the definition, $M$ is $(2g_0-1,g_0)$-small in this subcase, and Alternative (ii) of the conclusion holds with $T=T_1$ and $g'=g_0$.

There remains the subcase in which $M$ contains an incompressible surface whose genus is strictly less than  $2g_0-1$. Let
$g_2$ denote the smallest positive integer which occurs as the genus of a separating incompressible surface in $M$; then  $g_2\le2g_0-2$ and $g_2\ge2$ by hyperbolicity. Let us fix a separating incompressible surface $T_2$ of genus $g_2$. According to our choice of $g_0$, every  incompressible surface in $M$ has genus at least $g_0\ge g_2/2+1$; and according to our choice of $g_2$, every separating  incompressible surface in $M$ has genus at least $g_2$. Thus, according to the definition, $M$ is $(g_2,g_2/2+1)$-small in this case. Since in addition we have $g_2\le2g_0-2\le2g-2$,  Alternative (i) of the conclusion holds with $T=T_2$ and $g'=g_2$.
\EndProof

We now review some more definitions from \cite{CDS}. The Euler characteristic of a finitely triangulable space $Y$ will be denoted $\chi(Y)$, and we will set $\chibar(Y)=-\chi(Y)$. If $S$ is an incompressible surface in a closed, irreducible $3$-manifold $M$, we denote by $M\setminus\setminus S$ the manifold with boundary obtained by splitting $M$ along $S$. Each component  of the manifold $M\setminus\setminus S$ is irreducible and boundary-irreducible in the sense of \cite{hempel}. Furthermore, each component of $M\setminus\setminus S$ is {\it strongly atoral} in the sense that its fundamental group has no rank-$2$ free abelian subgroup. 

Any compact, orientable $3$-manifold $K$ which is irreducible and boundary-irreducible has a well-defined relative characteristic submanifold $\Sigma_K$ in the sense of \cite{Jo} and \cite{JS}. (In the notation of \cite{JS}, $(\Sigma_K,\Sigma_K\cap\partial K)$ is the characteristic pair of $(K,\partial K)$.) 

If $B$ is a  compact, orientable $3$-manifold  whose components are irreducible and boundary-irreducible, we denote by $\Sigma_{B}\subset B$ the union of the submanifolds $\Sigma_K$, where $K$ ranges over the components of $B$. In the case where the components of $B$
are strongly atoral,
%=M\setminus\setminus S$ for some closed, orientable hyperbolic $3$-manifold $M$ and some incompressible surface $S\subset M$, each 
component $C$ of $\Sigma_B$ may be given the structure of an $I$-bundle over a compact $2$-manifold with boundary in such a way that $C\cap\partial B$ is the associated $\partial I$-bundle. We denote by $\kish(B)$ the union of all components of $\overline{B-\Sigma_B}$ that have (strictly) negative Euler characteristic. 

Let $B$ be a  compact, orientable $3$-manifold  whose components are irreducible, boundary-irreducible, and strongly atoral. To say that  $B$ is {\it acyclindrical} means that $\Sigma_B=\emptyset$; this is equivalent to saying that $\kish(B)=B$.

\Proposition\label{really from CDS}
Let $g$ be a positive integer, and $M$ be a closed, orientable hyperbolic $3$-manifold that contains an incompressible surface of genus $g$. Suppose that the Heegaard genus of $M$ is strictly greater than $2g+1$. 
Then there exist an integer $g'\ge2$ and an incompressible surface $S\subset M$ of genus $g'$, such that either 
\Alternatives
\item  $g'\le 2g-1$, the surface $S$ separates $M$, and $M\setminus\setminus S$ has an acylindrical component; or
\item $g'\le 2g-1$, the surface $S$ separates $M$, and for each component $B$ of $M\setminus\setminus S$ we have $\kish(B)\ne\emptyset$; or
\item  $g'\le g$, the surface $S$ does not separate $M$, and $\chibar(\kish(M\setminus\setminus S))\ge2g'-2$.
\EndAlternatives
\EndProposition

\Proof
Let $G\ge 2g+2$ denote the Heegaard genus of $M$.

The hypothesis of the present proposition include those of Lemma \ref{small or tall}. Hence there exist an integer $g'$ and an incompressible surface $T\subset M$ of genus $g'$, such that one of the alternatives (i) or (ii) of the conclusion of Lemma \ref{small or tall} holds.

Consider first the case in which Alternative (i) of Lemma \ref{small or tall} holds. Thus $g'\le2g-2$, the surface $T$ separates $M$, and $M$ is $(g',g'/2+1)$-small. We have $G\ge2g+2\ge g'+4$. We now apply Theorem 5.1 of \cite{CDS}, with $g'$ playing the role of $g$ in that theorem. The theorem asserts that if 
$M$ is a closed, orientable 3-manifold containing a
separating incompressible surface of some genus $g'$, and if $M$ has Heegaard genus at least $g'+4$ and is $(g',g'/2+1)$-small, then $M$ contains a separating incompressible surface $S$ of genus $g$ such that either 
$M\setminus\setminus S$ has at least one acylindrical  component, or $\kish(B)\ne\emptyset$  for each component $B$ of $M\setminus\setminus S$. Thus in this case, one of the alternatives (i), (ii) of the present proposition holds.

Now consider the case in which Alternative (ii) of Lemma \ref{small or tall} holds. Thus  $g'\le g$, the surface $T$ does not separate $M$, and $M$ is $(2g'-1,g')$-small. We have $G\ge2g+2\ge2g'+2$. In this case we set $S=T$, and we apply Theorem 3.1 of \cite{CDS}, with $g'$ playing the role of $g$ in that theorem. The theorem asserts that if  $M$ is a closed, orientable, hyperbolic 3-manifold containing a non-separating  incompressible surface $S$ of genus $g'$, and if $\chibar(\kish(M\setminus\setminus S)) < 2g'-2$, and if $M$ is $(2g'- 1,g')$-small, then the Heegaard genus $G$ of $M$ is at most $2g' + 1$. Since in the present situation we have $G\ge2g'+2$, we must have $\chibar(\kish(M\setminus\setminus S))\ge2g'-2$. Thus in this case,  Alternative (iii) of the present proposition holds.
\EndProof

%Then $M$ contains a closed, incompressible (embedded) surface $S$, whose genus is strictly positive and is at most $g$, such that 

\Proposition\label{from CDS}
Let $M$ be a closed, orientable hyperbolic $3$-manifold, and let $g$ be a positive integer. Suppose that $M$ contains an  incompressible surface of genus $g$, and that the Heegaard genus of $M$ is strictly greater than $2g+1$. Then $\vol M>6.45$. 
\EndProposition

\Proof
Since $g$ is the genus of some incompressible surface in the closed, orientable hyperbolic $3$-manifold $M$, we have $g\ge2$. Thus the hypothesis implies that the Heegaard genus of $M$ is at least $6$. According to \cite[Theorem 6.1]{CDS}, a closed, orientable hyperbolic $3$-manifold which contains an  incompressible surface of genus $2$, and has Heegaard genus at least $6$, must have volume greater than $6.45$. Thus the conclusion of the proposition is true if $M$ contains an  incompressible surface of genus $2$. For the rest of the proof, we will  assume that $M$ contains no such surface.

The hypotheses of the present proposiotion are the same as those of Proposition \ref{really from CDS}. Hence we may fix an integer $g'\ge2$ and an incompressible surface $S\subset M$ of genus $g'$, such that one of the alternatives (i)---(iii) of the conclusion of Proposition \ref{really from CDS} holds.
Since  we have assumed that $M$ contains no  incompressible surface of genus $2$, we have $g'\ge3$.

We claim:
\Equation\label{yes it is}
\chibar(\kish(M\setminus\setminus S))\ge2.
\EndEquation

To prove (\ref{yes it is}), we first consider the case in which Alternative (i) of Proposition \ref{really from CDS} holds. In this case we may  label the components of $M\setminus\setminus S$ as $B_1$ and $B_2$ in such a way that  $B_1$  is acylindrical. Then 
$\kish(B_1)=B_1$; and since by definition every component of $\kish(B_2)$ has negative Euler characteristic, we have  $\chibar(\kish(B_2))\ge0$. Hence
$\chibar(\kish(M\setminus\setminus S))=\chibar(\kish(B_1))+\chibar(\kish(B_2))\ge \chibar(\kish(B_1)) =\chibar(B_1)=\chibar(S)/2=g'-1$.
Since $g'\ge3$, (\ref{yes it is}) is established in this case. 

In the case where Alternative (ii) of Proposition \ref{really from CDS} holds, let $B_1$ and $B_2$ denote the components of $M\setminus\setminus S$. For $i=1,2$, since $\kish(B_i)\ne\emptyset$, and since by definition every component of $\kish(B_I)$ has negative Euler characteristic, we have $\chibar(\kish(B_i)\ge1$; hence $\chibar(\kish(M\setminus\setminus S))=\chibar(\kish(B_1))+\chibar(\kish(B_2)\ge2$, and (\ref{yes it is}) holds. In the case where Alternative (iii) of Proposition \ref{really from CDS} holds, we have $\chibar(\kish(M\setminus\setminus S))=\chibar(\kish(B_1))+\chibar(\kish(B_2)\ge2g'-2\ge2$; thus (\ref{yes it is}) is established in all cases.

Let $\voct=3.66\ldots$ denote the volume of a regular ideal hyperbolic octahedron. According to \cite[Theorem 9.1]{AST}, whenever $S$ is an incompressible surface in a closed, orientable hyperbolic $3$-manifold $M$, we have
$\vol M\ge\voct\chibar(\kish(M\setminus\setminus S))$. In the present situation, (\ref{yes it is}) then gives $\vol M\ge2\voct>6.45$.
\EndProof

\Proposition\label{3.77 homology}
If $M$ is any closed, orientable hyperbolic $3$-manifold with $\vol M\le 3.77$, we have $\dim H_1(M;\FF_2)\le 10$. 
\EndProposition

\Proof
Assume that $\dim H_1(M;\FF_2)\ge 11$. We apply \cite[Proposition 8.1]{singular-two}, which asserts that if $k\ge3$ is an integer, and $ M$ is a closed orientable hyperbolic $3$-manifold with  $\dim H_1(M;\FF_2)\ge\max(3k-4, 6)$, then either $\pi_1(M)$ is $k$-free, or $M$ contains an incompressible surface of genus at most $k-1$. (The actual statement of the result quoted is slightly stronger than this.) Taking $k=5$, we deduce that either $\pi_1(M)$ is $5$-free, or $M$ contains an incompressible surface of genus at most $4$. If 
$\pi_1(M)$ is $5$-free, then 
Proposition \ref{five-free case} 
above gives $\vol M>3.77$, a contradiction to the hypothesis. Now suppose that $M$ contains an incompressible surface of some genus $k\le4$. We have $\dim H_1(M;\FF_2)\ge 11>9\ge2k+1$. In particular, the Heegaard genus of $M$ is strictly greater than $2k+1$. Proposition \ref{from CDS} now gives $\vol
 M>6.45$, and again the hypothesis is contradicted.
\EndProof

\section{
Dehn Drilling}
\label{drilling section}

The goal of this section is to prove 
Proposition \ref{good homology}, 
which was discussed in the Introduction.

\Number\label{drilling construction}
If $c$ is a simple closed geodesic in a closed, orientable hyperbolic $3$-manifold $M$, then according to 
Proposition 4 of \cite{kojima-geodesic}, a result that appears to be due to T. Sakai (see \cite{sakai}),
the manifold $M-|c|$ is homeomorphic to a one-cusped finite-volume hyperbolic $3$-manifold, which by Mostow rigidity is unique up to isometry. This hyperbolic manifold will be denoted by $M_c$.
\EndNumber

We will use the following result due to Agol and Dunfield:

\Theorem\label{AD-theorem}
Let $M$ be a hyperbolic 3-manifold, and let $l$ and $R$ be positive real numbers. Let $c$ be a simple closed geodesic in $M$, and suppose that the length of $c$ is at most $l$ and that its tube radius (see \ref{nbhd}) is at least $R$. 
Then
$$\vol M_c\le(\coth ^3 (2R)) \bigg(\vol M+ \frac\pi2l\tanh R\tanh2R \bigg).$$
\EndTheorem

\Proof
This is Theorem 10.1 of \cite{AST}.
\EndProof

\Reformulation\label{AD reformulation}
We define a function $\VAD$ on $(0,\infty)^3$ by
$$\VAD(V,R,l)=     V\tanh ^3 (2R)-\frac\pi2 l\tanh R\tanh2R.$$ (The subscript ``AD'' stands for ``Agol-Dunfield.'')
In this notation,  we may reformulate Theorem \ref{AD-theorem} as follows. Let $M$ be a hyperbolic 3-manifold,  let $l$ and $R$ be positive real numbers, let $c$ be a simple closed geodesic in $M$, and suppose that $c$ has length at most $l$ and tube radius 
at least $R$.
Then $\vol M\ge\VAD(\vol M_c ,R,l)$.
\EndReformulation

\Lemma\label{topper stopper}
Let $M$ be a closed, orientable hyperbolic $3$-manifold such that $\dim H_1(M;\FF_2)\ge 6$, let $\delta$ be a positive number, and let $c$ be a closed geodesic in $M$ whose length $l$ is at most $\delta$. Then for every point $p\in|c|$ we have $\shortone(p) = l$ and $\nextone(p)\ge\f_1( l)$. Furthermore, if $\eta$ is any number such that $0\le \eta\le\nextone(p)-\f_1(l)$ for every $p\in|c|$, then
%\fand if for every point $p\in|c|$ we have $\nextone(p)\ge \f_1(l)+\eta$, then 
%$c$ is simple, and for every point $p$ of $|c|$ the quantity $h_p:=\nextone(p)-\f_1(l)$ is non-negative. Furthermore, there is a point $p_0$ of $|c|$  such that for every $\delta\ge l$ and every $h\le h_{p_0}$
%we have 
$$\vol M\ge
\VAD\bigg(5.06,\frac12\arccosh\bigg(\frac{\cosh(\f_1(\delta) +\eta)}{\cosh(\delta /2)}\bigg),\delta \bigg).$$
(Here the quantity $\arccosh(\cosh(\f_1(\delta) +\eta)/\cosh(\delta /2))$ is defined because (\ref{Q-def}) gives $\f_1(\delta)+\eta\ge\f_1(\delta)>\log3>\delta>\delta/2$.)

\EndLemma

\Proof
%To prove the first assertion, note that 
Since we have in particular that  $\dim H_1(M;\FF_2)\ge 4$, it follows from \cite[Corollary 1.9]{sw} that $\pi_1(M)$ is $2$-free. It then follows from
\cite[Corollary 4.2]{acs-surgery} that
$\log3$ is a Margulis number for $M$. If $p$ is any point of $|c|$, then since $l<\delta<\log3$, we have $\shortone(p) = l$ by the first assertion of Lemma \ref{surgeon special}. Since $\pi_1(M)$ is $2$-free, it now follows from Corollary \ref{nebulizer} that  
$\nextone(p)\ge\f_1(\shortone(p))=\f_1(l)$; this establishes the first sentence of the conclusion of the lemma.

Now suppose that $\eta$ is a non-negative number such that $ \eta\le\nextone(p)-\f_1(l)$ for every $p\in|c|$. Let us write $M=\HH^3/\Gamma$, where $\Gamma$ is a cocompact subgroup of $\isomplus(\HH^3)$. Then the pre-image of $|c|$ under the quotient map $\HH^3\to M$ is a line $A\subset\HH^3$, whose stabilizer in $\Gamma$ is a maximal cyclic subgroup $C$ of the \iccg\  $\Gamma$ (see \ref{i see cg}).
Fix  a generator $\gamma$  of $C$.

Let $r$ denote the tube radius of $c$. (A priori we have $r>0$ if and only if $c$ is  simple, see \ref{nbhd}.) According to 
\ref{nbhd},
there is an element $\alpha$ of $\Gamma-C$ such that the minimum distance between the lines $A$ and $\alpha\cdot A$ is equal to $2r$. 

Now let $L$ denote the common perpendicular to the lines $A$ and $\alpha\cdot A$, and let $ X $ and $ Y $ 
denote the respective points of intersection of $L$ with the lines $A$ and $\alpha\cdot A$. We have $\dist( X , Y )=2r$ (so that in particular $ X = Y $ if and only if $c$ is not simple). The points $ X $ and $\alpha^{-1}\cdot  Y $ lie on the line $A$. Since $\gamma$ acts on $A$ by a translation of length $l$, the orbit of $\alpha^{-1}\cdot  Y $ under $C$ contains a point $P_0$ whose distance from $ X $ is at most $l/2$. Thus for some $n\in\ZZ$ we have $P_0=\gamma^n\alpha^{-1}\cdot  Y $. The points $P_0$, $ X $ and $ Y $ are the vertices of a hyperbolic triangle with a right angle at $ X $, and the hyperbolic Pythagorean theorem gives
\Equation\label{purple dinkeler}
\cosh(\dist( Y ,P_0))=\cosh(2r)\cosh(\dist( X ,P_0)).
\EndEquation

We let $p_0\in|c|$ denote the image of $P_0$ under the quotient map.

Since $c$ has length $l$, and $P_0$ lies on the axis $A$ of $\gamma$, we have $d(\gamma,P_0)=l$.

We have observed that 
%Since we have  in particular that
%$\dim H_1(M;\FF_2)\ge 4$, it follows from 
%\cite[Corollary 1.9]{sw} that $\Gamma\cong\pi_1(M)$ is $2$-free. It then follows from
%\cite[Corollary 4.2]{acs-surgery} that
%$\log3$ is a Margulis number for $M$. Since $l<\delta<\log3$, we have 
$\shortone(p_0) = l$. 
% by the first assertion of Lemma \ref{surgeon special}.
Thus $d(\gamma,P_0)=\shortone(p_0)$. According to the discussion in \ref{short and next}, this implies that
%The following argument should be replaced by a ref. to Lemma 
%\ref{surgeon special}, which implies that  $\shortone(p) = l$. Of course $\ell_M$ will disappear from the argument.
%Since $c$ has length $\ell_M$, and $P_0$ lies on the axis $A$ of $\gamma$, we have $d(\gamma,P_0)=\ell_M=\min_{Q\in\HH^3,1\ne x\in\Gamma}d(x,Q)$. In particular we have
%$d(\gamma,P_0)=\min_{1\ne x\in\Gamma}d(x,P_0)$. Thus if $p_0\in M$ denotes the image of $P_0$ under the quotient map, the discussion in \ref{short and next} gives that  $d(\gamma,P_0)=\shortone(p_0)$. As $C$ is the maximal cyclic subgroup of the \iccg\ $\Gamma$ containing $\gamma$, the discussion in  \ref{short and next}) then gives 
$\nextone(p_0)=\min_{x\in\Gamma-C}d(x,P_0)$.

Set $\beta=\alpha\gamma^{-n}\in\Gamma$, so that $ Y =\beta\cdot P_0$. Since $\alpha\notin C$, we have $\beta\notin C$, and hence $\nextone(p_0)\le d(\beta,P_0)=\dist( Y ,P_0)$. In particular it follows that $\dist( Y ,P_0)>0$, which by (\ref{purple dinkeler}) implies that $r>0$; thus $c$ is simple.

The hypothesis implies that $\nextone(p_0)\ge \f_1(l)+\eta$. This, together with the inequality $\nextone(p_0)\le \dist( Y 
,P_0)$, shows that the left-hand side of (\ref{purple dinkeler}) is bounded below by $\cosh(\f_1(l)+\eta)$. Since in addition we have $\dist( X ,P_0)\le l/2$, it now follows from (\ref{purple dinkeler}) that 
\Equation\label{pie oh my}
\cosh(\f_1(l)+\eta)\le\cosh(2r)\cosh(l/2).
\EndEquation

Let us define a function $F$ on $(0,\infty)$ by $F(x)=\cosh(\f_1(x)+\eta)/\cosh(x/2)$. Then (\ref{pie oh my}) asserts that $\cosh(2r)\ge F(l)$. We claim that $F$ is strictly monotone decreasing on $(0,\infty)$. To show this, first note that the definitions of the functions $Q$ and $\f_1$ 
(see \ref{Q-def})
imply that $e^{\f_1(x)}=(e^x+3)/(e^x-1)$ for $x>0$. Hence if we set $A=e^{\eta}$, we have
$$\cosh(\f_1(x)+h)=\frac12\bigg(A\frac{e^x+3}{e^x-1}+A^{-1}\frac{e^x-1}{e^x+3}\bigg)$$
for $x>0$. It follows that $F(x)=G(e^{x/2})$, where $G$ is the rational function defined by
\Equation\label{moofoo}
\begin{aligned}
G(x)&=\frac1{u+u^{-1}}\bigg(A\frac{u^2+3}{u^2-1}+A^{-1}\frac{u^2-1}{u^2+3}\bigg)\\
&=\frac {Au}{u^2-1}\cdot\frac
{u^2+3}{u^2+1}+\frac {A^{-1} u}{u^2+3}\cdot\frac{u^2-1}{u^2+1}.
\end{aligned}
\EndEquation
Each term on the last line of (\ref{moofoo}) is a product of two strictly decreasing functions on $(1,\infty)$; hence $G$ is strictly decreasing on $(0,\infty)$, so that $F$ is strictly decreasing on $(0,\infty)$ as claimed. 
Since $l<\delta$ and $\cosh(2r)\ge F(l)$, we have $\cosh(2r)>F(\delta)=\cosh(\f_1(\delta)+\eta)/\cosh(\delta/2)$. In view of the first assertion of the lemma, this may be rewritten as
\Equation\label{footsie}
r\ge\frac12\arccosh\bigg(\frac{\cosh(\f_1(\delta)+\eta)
)}{\cosh(\delta/2)}\bigg).
\EndEquation

We apply \ref{AD reformulation}, letting the right hand side of (\ref{footsie}) play the role of the lower bound $R$ for the tube radius required for \ref{AD reformulation}, and recalling that $c$ has length  $l\le \delta$. This gives
\Equation\label{fundry}
\vol M\ge\VAD\bigg(\vol M_c,\frac12\arccosh\bigg(\frac{\cosh(\f_1(\delta)+\eta)
)}{\cosh(\delta/2)}\bigg),\delta\bigg).
\EndEquation

On the other hand, since the one-cusped orientable hyperbolic $3$-manifold $M_c$ is homeomorphic to the complement of a simple closed curve in $M$, we have $\dim H_1(M_c;\FF_2)\ge \dim H_1(M;\FF_2)\ge 6$ in view of the hypothesis. According to \cite[Theorem 6.2]{CDS}, if $N$ is any one-cusped orientable hyperbolic $3$-manifold such that $\dim H_1(N;\FF_2)\ge 6$, we have $\vol N>5.06$. Hence $\vol M_c>5.06$. Since it is obvious from the definition of the function $\VAD$ that it is monotone increasing in its first argument, the inequality (\ref{fundry}) continues to hold if the expression $\vol M_c$ is replaced by its lower bound $5.06$. This gives the conclusion of the lemma.
\EndProof
%XYZUV

\Proposition\label{good homology}
If $M$ is any closed, orientable hyperbolic $3$-manifold with $\vol M\le 3.69$, we have $\dim H_1(M;\FF_2)\le7$.
\EndProposition

\Proof
Equivalently, we must prove that if $\dim H_1(M;\FF_2)\ge8$ then $\vol M>3.69$. 
We apply \cite[Proposition 8.1]{singular-two}, which implies that if $k\ge3$ is an integer, and $ M$ is a closed orientable hyperbolic $3$-manifold with  $\dim H_1(M;\FF_2)\ge\max(3k-4, 6)$, then either $\pi_1(M)$ is $k$-free, or $M$ contains an incompressible surface of genus at most $k-1$. Taking $k=4$, we deduce that if $\dim H_1(M;\FF_2)\ge8$ then either $\pi_1(M)$ is $4$-free, or $M$ contains an incompressible surface of genus at most $3$. 

Consider first the case in which $M$ contains an incompressible surface of some genus $k\le3$. We have $\dim H_1(M;\FF_2)\ge 8>7\ge2k+1$. In particular, the Heegaard genus of $M$ is strictly greater than $2k+1$. Proposition \ref{from CDS} now gives $\vol
 M>6.45$, which is stronger than the required conclusion.

The rest of the argument will be devoted to the case in which
$\pi_1(M)$ is $4$-free. We fix a closed geodesic $c$ in $M$ having length $\ell_M$ 
(see \ref{nbhd}).
We distinguish four subcases, according to which of the following conditions holds:
\begin{itemize}
\item[(a)]
$\ell_M>0.5912$;
\item[(b)]
$\ell_M<0.5637$;

\item[(c)]
$0.5637\le\ell_M\le0.5912$, and for every point $p$ of 
$|c|$
we have $\nextone(p)>\f_1(\ell_M)+0.08267$;

\item[(d)]
$0.5637\le\ell_M\le0.5912$, and there is a point $p$ of 
$|c|$ 
for which $\nextone(p)\le\f_1(\ell_M)+0.08267$.

\end{itemize}

To prove that the inequality $\vol M\ge3.69$ always holds in Subcase (a), we apply Proposition \ref{latest monster}. We take $k=4$, 
$V_0=3.690003$,
$T=(0.5912,\infty)$, and $\waslambdaminus=(\log8)-0.0409$.

To define the function $\scrM$ on $[\waslambdaminus,\log8]$ that appears in the hypothesis of Proposition \ref{latest monster}, we first define positive numbers $\lambda_0,\ldots,\lambda_{200}$ and $\mu_1,\ldots,\mu_{200}$, with $\waslambdaminus=\lambda_0<\cdots<\lambda_{200}=\log8$ and $\mu_1>\cdots>\mu_{200}$, as follows. We set $\lambda_i=
(\log8)-.0409+.0002045i$ for $i=0,\ldots,200$. We set $\mu_i=1.146$ for
$1\le i\le19$; $\mu_i=1.142$ for $20\le i\le27$; $\mu_i=1.1405$ for $28\le i\le31$; $\mu_i=1.1395$ for $32\le\mu_i\le34$; $\mu_i=1.1386$ for $35\le i\le 37$; $\mu_i=
1.1379-.0002(i-38)$ for $38\le i\le42$; $\mu_i=1.1368-.0002(i-43)$ for $43\le i\le46$; $\mu_{47}=1.13597$; $\mu_{48}=1.13573$; $\mu_i=1.13548-.00021(i-49)$ for $49\le i\le 52$; and $\mu_i=1.13461-.00019(i-53)$ for $53\le i\le200$. We then define $\scrM$ by setting $\scrM(\waslambdaminus)=\mu_1$, and $\scrM(\lambda)=\mu_i$ whenever $i\in\{1,\ldots,200\}$ and $\lambda_{i-1}<\lambda\le\lambda_i$.

Note that the $\mu_i$ lie in the interval $I_4=(\log3,\log(17/3))$, so that $\scrM$ takes its values in this interval. 
We must verify Conditions (1)--(3) of Proposition \ref{latest monster}.
The verification of Condition (1) is a direct computation; the left-hand side of the inequality is $3.6904\ldots$. 

To verify 
Condition (2), let $\lambda\in[\lambda^-,\lambda^+]$ br given, and choose $i\in\{1,\ldots,m\}$ so that either $\lambda_{i-1}<\lambda\le\lambda_i$, or $i=1$ and $\lambda=\waslambdaminus$. Then $\scrM(\lambda)=\mu_i$. If we set $R=\fone(\lambda_i)/2$, we have $R\in(0,\fone(\lambda)/2]$ because $\fone$ is monotone decreasing (see \ref{Q-def}). Thus we have $\Wfive(R,\g_4(\scrM(\lambda)),\fthree(\scrM(\lambda)))=
\Wfive(\fone(\lambda_i)/2,\g_k(\mu_i),\fthree(\mu_i))$. It therefore suffices to show that
$\Wfive(\fone(\lambda_i)/2,\g_k(\mu_i),\fthree(\mu_i))> V_0$ for $i=1,\ldots,200$. This is done by a finite number of calculations; 
the smallest value of $\Wfive(\fone(\lambda_i)/2,\g_k(\mu_i),\fthree(\mu_i))$ is $3.69002\ldots$, 
and is achieved when $i=57$.

We will verify Condition (3) of Proposition \ref{latest monster} by establishing a stronger assertion: for every $\lambda \in[\waslambdaminus,\log8]$, and every $\wasell\in T=(0.5912,\infty)$, Alternative (a) of
Condition (3) of Proposition \ref{latest monster} holds. If $\lambda \in[\waslambdaminus,\log8]$ and $\wasell>0.5912$ are given, we again choose $i\in\{1,\ldots,m\}$ so that either $\lambda_{i-1}<\lambda\le\lambda_i$, or $i=1$ and $\lambda=\waslambdaminus$; thus we again have $\scrM(\lambda)=\mu_i$. We set $\alpha=\lambda_{i-1}\le\lambda$, and we set $\delta=0.5912<\wasell$. Note that we have $\lambda/4\le(\log8)/4<\delta$, and $\delta>\f_1(\waslambdaminus)\ge\f_1(\lambda)$ since $\f_1$ is monotone decreasing. Since $\delta>\f_1(\lambda)$,  the pair of inequalities $D\ge\scrD(\lambda)$ and $D>\fthree(\lambda)$ which appear in Alternative (a) of
Condition (3) of Proposition \ref{latest monster} (and are needed to guarantee that
$\waschi_k(\alpha, \lambda,\delta,D,\scrM(\lambda))$ 
is defined) are equivalent to the single inequality $D\ge\delta$.
Thus to complete the verification of Condition (3), we must show that   with these choices of $\alpha$ and $\delta$ we have $\waschi_4(\alpha,\lambda,\delta,D, \scrM(\lambda))>V_0$, i.e. 
\Equation\label{barrel}
\waschi_4(\lambda_{i-1}, \lambda,\delta,D, \mu_i)>V_0, 
\EndEquation
for every $D\ge\delta$.

We first prove (\ref{barrel}) for $\delta\le D\le\log8$. 
For $j=1,\ldots,2000$ we set $E_j=
((\log8)-.5912)j/2000$.
%Then 
For any $D\in[\delta,\log8]$,
% there is a $j\in\{1,\dots,2000\}$ such that $\delta+E_{j-1}\le D\le\delta+E_j$. Since $D\ge\delta$, 
we may write $D=E+\delta$ for some $E$ with $0\le E\le(\log8)-\delta=E_n$. We then have $E_{j-1}\le E\le E_j$ for some $j\in\{1,\ldots,n\}$. Applying Lemma \ref{now a cow}, with $\lambda_{i-1}$, $\lambda_{i}$, $E_{j-1}$, $E_j$,  and $\mu_i$ playing the respective roles of $\lambda^-$, $\lambda^+$, $E^-$, $E^+$,  and $\mu$ in that lemma,  we find
$$
\waschi_4(\lambda_{i-1}, \lambda,\delta,D,\mu_i)
\ge
\waschi^-_4(\lambda_{i-1},\lambda_i,\delta,E_{j-1},E_j,\mu_i).
$$

It therefore suffices to show that we have $\waschi^-_4(\lambda_{i-1},\lambda_i,\delta,E_{j-1},E_j,\mu_i)>V_0$ for all $i\in\{1,\ldots,200\}$ and all $j\in\{1,\ldots,2000\}$. This is done by a finite number calculations, 
using the method described in Subsection \ref{calculations}; 
the smallest value of 
$\waschi^-_4(\lambda_{i-1},\lambda_i,\delta,E_{j-1},E_j,\mu_i)$
is 
$3.6900044\ldots$
and is achieved when 
$i=49$
and 
$j=1$.

We now turn to the proof of (\ref{barrel}) when $D>\log8$. The definition of $\waschi_k$ (see \ref{chi def}) immediately implies that 
$\waschi_4(\lambda_{i-1}, \lambda,\delta,D,\mu_i)\ge\Wone(\lambda_{i-1},\delta,D)$. 
We have $D>\log8=\lambda_m>\lambda_{i-1}$, and we may therefore apply Lemma 
\ref{other woney}, with $\waslambdaminus$, $\log8$,  and $\lambda_{i-1}$ playing the respective roles of $\maybecalLminus$, $\maybecalLplus$,  and $\alpha$ in that lemma,  to deduce that
\Equation\label{dfoo}
\Wone(\lambda_{i-1},\delta,D)\ge B\bigg(\frac{\waslambdaminus}2\bigg)-2\kappa\bigg(\frac{\log8}2, \frac{\Phi_2(\delta,\log8)}2\bigg) -2\kappa\bigg(\frac{\log8}2, \frac{\Phi_3(\delta,\log8)}2\bigg).
\EndEquation
Recalling that $\waslambdaminus=(\log8)-0.0409$ and that $\delta=0.5912$, we find by direct computation that the right hand side of (\ref{dfoo}) is equal to $5.264\ldots>V_0$.
This completes the verification of Condition (3) of Proposition \ref{latest monster}, and completes the proof of the assertion of the proposition in Subcase (a).

To prove the assertion in Subcase (b), we apply Lemma 
\ref{topper stopper}, taking $\eta=0$, $\delta=0.5637$, taking $M$ to be the manifold given by the hypothesis of the present proposition, and taking $c$ to be  a closed geodesic in $M$  whose length is $\ell_M$; thus we have $l=\ell_M$ in the notation of Lemma \ref{topper stopper}. For every point $p$ of $|c|$ we have $\nextone(p)\ge\f_1(\ell_M)=\f_1(\ell_M)+\eta$ by Lemma \ref{topper stopper}.
It now follows from Lemma \ref{topper stopper} that 
$$\vol M\ge
\VAD\bigg(5.06,\frac12\arccosh\bigg(\frac{\cosh(\f_1(0.5637))}{\cosh(0.5637/2)}\bigg),0.5637\bigg) =3.69019\ldots.$$

To prove the assertion in Subcase (c), we again apply Lemma 
\ref{topper stopper}, now taking $\eta=0.08267$ and $\delta=0.5912$. We again take $M$ to be the manifold given by the hypothesis of the present proposition, and take $c$ to be  a closed geodesic in $M$  whose length is $\ell_M$; thus again we have $l=\ell_M$ in the notation of Lemma \ref{topper stopper}. 
It follows from Lemma \ref{topper stopper} that 
$$\vol M\ge
\VAD\bigg(5.06,\frac12\arccosh\bigg(\frac{\cosh(\f_1(0.5912)+0.08267)}{\cosh(0.5912/2)}\bigg),0.5912\bigg) =3.690003\ldots.$$

It remains to establish the assertion in Subcase (d). In this subcase, we choose  a point $p_0$ of 
$|c|$ 
for which $h:=\nextone(p_0)-\f_1(\ell_M)\le0.08267$. 

Since $\pi_1(M)$ is in particular $3$-free, it follows from
\cite[Corollary 9.3]{acs-singular} that $M$ contains a hyperbolic ball of radius $(\log5)/2$; that is,
$\max_{p\in M}\shortone(p)\ge \log5$. If we set $\mu_0=1.12235$, we find by direct calculation that 
$\Wfive((\log5)/2,\g_4(\mu_0), \fthree(\mu_0))=3.6904\ldots>V_0$. Applying Corollary \ref{marseillaise},
with $k=4$ and with $V_0$, $\log5$  and $\mu_0$ playing the respective roles of $V$, $R$  and $\mu$, we deduce that either $\mu_0$ is a Margulis number for $M$ or $\vol M>V_0$. We may therefore assume that $\mu_0$ is a Margulis number for $M$. 

It now follows
from Lemma \ref{stoppeth one} 
that $h$ is strictly positive (so that $h\in(0,.08267]$) and that
$\vol M
\ge \waspsi_4(h,\ell_M,\mu_0)$.

Let us define positive numbers
$\zeta_0,\ldots,\zeta_{246}$
and  $h_0,\ldots,h_{2000}$ as follows. We set $\zeta_s=0.5637+0.0005s$ for $0\le s\le 40$; $\zeta=0.5837+0.00005(s-40)$ for $41\le s\le176$; and $\zeta_s=0.5905+0.00001(s-176)$ for $177\le s\le246$. We set 
$h_t=0.08267t/2000$ for $0\le t\le2000$. Thus we have $0.5637=\zeta_0<\cdots<\zeta_{246}=0.5912$ and $0=h_0<\cdots<h_{2000}=0.08267$.
Since $\ell_M\in[0.5637,0.5912]$ and $h\in[0,0.08267]$, we may fix indices $s\in\{0,\ldots,246\}$ and $t\in\{0,\ldots,2000\}$ such that $\zeta_{s-1}\le\ell_M\le\zeta_s$ and $h_{t-1}\le h\le h_t$. We apply Lemma \ref{cranky-poo}, with $\zeta_{s-1}$, $\zeta_{s}$, $h_{t-1}$ and $h_{t}$ playing the respective roles of $l^-$, $l^+$, $h^-$ and $h^+$, and with $\mu$ defined as above. This shows that 
$\waspsi_4( h,\ell_M,\mu_0)\ge\waspsi^-_4( h_{t-1},h_t,\zeta_{s-1},\zeta_s,\mu_0)$. It therefore suffices to show that $\waspsi^-_4( h_{t-1},h_t,\zeta_{s-1},\zeta_s,\mu_0)>V_0$ for $s=1,\ldots, 246$ and for $t=1,\ldots,2000$. This is done by a finite number calculations,
using the method described in Subsection \ref{calculations}; 
the smallest value of $\waspsi^-_4( h_{t-1},h_t,\zeta_{s-1},\zeta_s,\mu_0)$
is 
$3.690103\ldots$
and is achieved when 
$s=246$ 
and 
$t=2000$.
\EndProof

\bibliographystyle{plain}

\end{document}